\documentclass[pdflatex,sn-mathphys-num]{sn-jnl}% Math and Physical Sciences Numbered Reference Style
\usepackage{graphicx}%
\usepackage{multirow}%
\usepackage{amsmath,amssymb,amsfonts}%
\usepackage{amsthm}%
\usepackage{mathrsfs}%
\usepackage[title]{appendix}%
\usepackage{xcolor}%
\usepackage{textcomp}%
\usepackage{manyfoot}%
\usepackage{booktabs}%
\usepackage{algorithm}%
\usepackage{algorithmicx}%
\usepackage{algpseudocode}%
\usepackage{listings}%
\usepackage[utf8]{inputenc}
\usepackage{bbm}
\usepackage{hyperref}
\usepackage{cleveref}
\usepackage{mathtools}
\usepackage[dvipsnames]{xcolor}
\usepackage{caption}
\usepackage{placeins}
\usepackage{diagbox}
\usepackage{tabularx}
\usepackage{cancel}
\usepackage[symbol]{footmisc}
\usepackage[T1]{fontenc}

\usepackage{scalerel}
\usepackage{stackengine,wasysym}

\numberwithin{equation}{section}%

\newtheorem{lemma}{Lemma}[section]

\newcommand\reallywidetilde[1]{\ThisStyle{%
  \setbox0=\hbox{$\SavedStyle#1$}%
  \stackengine{-.1\LMpt}{$\SavedStyle#1$}{%
    \stretchto{\scaleto{\SavedStyle\mkern.2mu\AC}{.5150\wd0}}{.6\ht0}%
  }{O}{c}{F}{T}{S}%
}}

\newcommand{\E}{\mathbbm{E}}
\newcommand{\R}{\mathbbm{R}}
\newcommand{\bcdot}{\ensuremath{\boldsymbol{\cdot}}}
\newcommand{\Av}{A^{\scriptscriptstyle v}}
\newcommand{\AT}{A^{\scriptscriptstyle T}}
\newcommand{\AS}{A^{\scriptscriptstyle S}}

\newcommand{\tr}{{\scriptscriptstyle \rm T}}

\newcommand{\CQFD}{\begin{flushright}     $\square$ \end{flushright}}

\usepackage{amsmath,array}
\DeclareMathOperator*{\esssup}{ess\,sup}

\newcommand{\modif}[1]{\color{black} #1 \color{black}}

\theoremstyle{thmstyleone}%
\newtheorem{theorem}{Theorem}%  meant for continuous numbers
\theoremstyle{thmstyletwo}%
\newtheorem{remark}{Remark}%

\theoremstyle{thmstylethree}%

\begin{document}

\title[Stochastic interpretations of the oceanic primitive equations with relaxed hydrostatic assumptions]{Stochastic interpretations of the oceanic primitive equations with relaxed hydrostatic assumptions}

%%=============================================================%%
%% GivenName	-> \fnm{Joergen W.}
%% Particle	-> \spfx{van der} -> surname prefix
%% FamilyName	-> \sur{Ploeg}
%% Suffix	-> \sfx{IV}
%% \author*[1,2]{\fnm{Joergen W.} \spfx{van der} \sur{Ploeg} 
%%  \sfx{IV}}\email{iauthor@gmail.com}
%%=============================================================%%

\author[1]{\fnm{Arnaud} \sur{Debussche}} \email{arnaud.debussche@ens-rennes.fr}

\author[2]{\fnm{Étienne} \sur{Mémin}} \email{etienne.memin@inria.fr}
%\equalcont{These authors contributed equally to this work.}

\author*[2]{\fnm{Antoine} \sur{Moneyron}} \email{antoine.moneyron@inria.fr}
%\equalcont{These authors contributed equally to this work.}

\affil[1]{Univ Rennes, CNRS, IRMAR UMR 6625, F35000, Rennes, France.}

\affil[2]{Univ Rennes, INRIA, IRMAR UMR 6625, F35000, Rennes, France.}
%%==================================%%
%% Sample for unstructured abstract %%
%%==================================%%

\abstract{In this paper, we investigate how weakening the classical hydrostatic balance hypothesis impacts the well-posedness of the stochastic LU primitive equations. The models we consider are intermediate between the incompressible 3D LU Navier-Stokes equations and the LU primitive equations with standard hydrostatic balance. As such, they are expected to be numerically tractable, while accounting well for phenomena within the grey zone between hydrostatic balance and non-hydrostatic processes. Our main result is the well-posedness of a low-pass filtering-based stochastic interpretation of the LU primitive equations, with rigid-lid type boundary conditions, in the limit of ``quasi-barotropic'' flow. This assumption is linked to the structure assumption proposed in \cite{AHHS_2022, AHHS_2022_preprint}, which can be related to the dynamical regime where the primitive equations remain valid \cite{MHPA_1997}. Furthermore, we present and study two eddy-(hyper)viscosity-based models.}

%%\pacs[JEL Classification]{D8, H51}

%%\pacs[MSC Classification]{35A01, 65L10, 65L12, 65L20, 65L70}

\maketitle

\section*{Introduction}

Stochastic modelling for large-scale fluid flows and their dynamics representation is nowadays a key research topic. Geophysical flows cannot be numerically represented in their full complexity, since they are characterised by a fully developed turbulence and chaotic dynamical systems. Therefore, only approximated large-scale models can be computationally considered. In the last years, stochastic modelling has emerged as a powerful setting for such representation \cite{Berner-Al_2017, FO_2017, FOBWL_2015, MTV_1999}. In particular, climate change studies call for models able to handle approximations or errors together with their time evolution. The need for plausible realisations in probabilistic forecasting requires carrying out a set of physically relevant realisations, with efficient uncertainty quantification properties. Stochastic modelling and stochastic calculus offer practical methodological frameworks to deal with these issues. The first models of this type were based on phenomenological turbulence studies on backscattering energy across scales \cite{Leith_1990, MT_1992}, and typically involve multiplicative random forcing models or stochastic parametrisation \cite{BMP_1999, Shutts_2005}. However, the noise variance being uncontrolled \emph{a priori}, an eddy viscosity was added to balance the noise energy. The precise form of this additional viscosity is still unknown, but it often relies on the debatable Boussinesq turbulence assumption in practice \cite{Schmitt_2007}. Also, random forcing defined outside of any conservation principle may lead to a poor accuracy compared to the reference fine deterministic resolution \cite{CDMR_2018}, and to a lack of interpretability.

During the past decade, the location uncertainty approach (LU) has been developed and studied to propose physically consistent stochastic models \cite{Mémin_2014, TML_2023}. This approach relies on a stochastic version of the Reynolds transport theorem, applied to the mass, momentum and energy conservation \cite{Mémin_2014}. The LU formalism has been successfully applied to classical geophysical models \cite{BCCLM_2020, RMC_2017_Pt1, RMC_2017_Pt2, RMC_2017_Pt3}, stochastic reduced order models \cite{RMHC_2017, RPMC_2021, TCM_2021} and large eddy simulation models \cite{CMH_2020, CHLM_2018, HM_2017}. Its physical relevance has also been tested on prototypical flow models \cite{BCCLM_2020, BLBM_2021, CDMR_2018}. Additionally, in a recent paper \cite{DHM_2023}, the authors demonstrated that the stochastic version of the 2D Navier-Stokes equation under location uncertainty is well-posed, and that the 3D one admits a martingale solution. This is consistent with the results in the deterministic setting, as the well-posedness of the classical 3D Navier-Stokes equation remains an open problem. Another important result was the continuity of the stochastic model for a vanishing noise, which shows the consistency of the stochastic model compared to the associated deterministic one. The noise considered in the LU setting corresponds to the so-called transport noise. This type of noise has been subject to intensive research efforts in the mathematics community due to the need for well-posedness properties for fluid dynamics models \cite{AHHS_2022, AHHS_2022_preprint, BS_2021, CFH_2019, DHM_2023, FGL_2021, FGP_10, FL_2021, GCL_2023, LCM_2023, MR_2005}, and because of the emergence of enhanced dissipation and mixing \cite{FGL_2022, FR_2023}.

In the deterministic setting, the primitive equations are commonly used to model geophysical flows \cite{Vallis2017}. They are derived from the 3D Navier-Stokes equations, making the hypothesis whereby the vertical acceleration is negligible. This leads to the classical hydrostatic equilibrium on the vertical component, which relates the vertical derivative of the pressure to the density fluctuation. This balance is known to be physically valid in the ocean at large-scale. However, it breaks down outside of the shallow water regime, or when thermodynamic effects take place, such as deep convection. We call this assumption the \emph{strong hydrostatic hypothesis}. Remarkably, this deterministic model is known to be well-posed \cite{CT_2007} under suitable boundary conditions -- namely the rigid-lid ones -- the proof relying essentially on the study of barotropic and baroclinic velocity modes. Various authors have studied stochastic versions of these equations. The well-posedness of the stochastic primitive equations with multiplicative noise has been shown in \cite{DGHT_2011, DGHTZ_2012}, and more recently with a specific class of regular transport noise in \cite{BS_2021}. Moreover, in a recent paper, it has been shown that the stochastic primitive equations with transport noise that is similar to the LU one are well-posed, under the strong hydrostatic hypothesis and using water world-type boundary conditions \cite{AHHS_2022}. However, the horizontal noise was assumed to be independent of the vertical axis in this work, which makes the barotropic and baroclinic noises tractable.

As pointed out previously, the strong hydrostatic balance does not allow to represent processes with a non-negligible vertical acceleration such as deep convection phenomena, which exhibit strong up- or down-welling of water. Yet, modelling these phenomena is crucial to represent faithfully thermohaline circulation and deep currents -- such as the Atlantic meridional overturning circulation (AMOC) -- in the long-run and at coarse spatial resolution (i.e. at the climatic scale). The LU setting allows to relax easily the strong hydrostatic balance by considering the martingale terms of the vertical acceleration as deviation terms. This immediately yields a generalization of the primitive equations. Therefore, in this paper, we study the well-posedness of the LU primitive equations under a weaker hydrostatic equilibrium assumption: we account for the transport of the vertical velocity by the noise. However, the choice of the shape for the energy compensation dramatically impacts the theoretical properties of the system. This leads to two potential models of interest. Such models are intermediate between the stochastic primitive equations under strong hydrostatic hypothesis and the 3D Navier-Stokes equations -- with Boussinesq's assumption of weak compressibility. Also, they are expected to be numerically tractable, while accounting better for non-hydrostatic phenomena, such as solitons, internal waves or oceanic convection, to name but a few.

The paper is organised as follows. First, we introduce the LU framework and precise our assumptions to derive a suitable class of LU representations for the primitive equations. Then we describe the functional setting of this study and state the main results of the paper. After this, we show the existence of a global martingale solution for a class of models with low-pass filtered noise. Then we show that there exists a unique local pathwise solution for such models. Considering an approximation of this low-pass filtered model in the limit of ``quasi-barotropic'' flow, and improving the regularity of the filtering kernel, we prove eventually that there exists a unique global-in-time pathwise solution. Additionally, we show that the solution to this latter model is continuous with respect to the initial data and to the noise data in a well-chosen topology. Moreover, we present two non-filtered (hyper)diffusive models, for which there exist global martingale solutions. Sketches of proofs for their $L^2$-energy estimates are given in the appendix, the rest of the proofs being similar to the one for the filtered model.

\section{Primitive equations models in the LU framework}

The LU formulation is based on the following time-scale separation of the flow:
\begin{equation}
    dX_t = u(X_t,t) dt + \sigma(X_t,t) dW_t.
\end{equation}
In this decomposition, which should be understood in the It\={o} sense, $X$ denotes the Lagrangian displacement defined in a bounded cylindrical tridimensional domain $\mathcal{S} = \mathcal{S}_H \times [-h,0] \subset \R^3$, where $\mathcal{S}_H$ is a subset of $\R^2$ with smooth boundary. This formulation corresponds to a flat bottom assumption, even if we expect the results presented in this paper to be similar for a smooth enough non-flat bottom. The velocity component $u(X_t,t)$ denotes the large-scale Eulerian velocity (correlated in both space and time) and $\sigma(X_t,t) dW$ is a highly oscillating unresolved velocity (uncorrelated in time but correlated in space). We interpret this second component as a noise term in the following.

Let us define this noise term more precisely: consider a Wiener process $W$ on the space of square integrable functions $\mathcal{W} := L^2(\mathcal{S},\R^3)$. Thus, there exists a Hilbert orthonormal basis $(e_i)_{i \in \mathbbm{N}}$ of $\mathcal{W}$ and a sequence of independent standard Brownian motions $(\hat{\beta}^i)_{i \in \mathbbm{N}}$ on a filtered probability space $(\Omega, \mathcal{F}, (\mathcal{F}_t)_t, \mathbbm{P})$ such that,
\begin{equation*}
    W=\sum_{i \in \mathbbm{N}} \hat{\beta}^i e_i.
\end{equation*}
As such, $(\Omega, \mathcal{F}, (\mathcal{F}_t)_t, \mathbbm{P}, W)$ is a stochastic basis. Note that the sum $\sum_{i \in \mathbbm{N}} \hat{\beta}^i e_i$ does \emph{not} converge in $\mathcal{W}$. Hence, the previous identity only makes sense in a space $\mathcal{U}$ including $\mathcal{W}$, such that the embedding $\mathcal{W} \hookrightarrow \mathcal{U}$ is Hilbert-Schmidt. For instance, $\mathcal{U}$ can be the dual space of any reproducing kernel Hilbert subspace
of $\mathcal{W}$ for the inner product $(\bcdot, \bcdot)_{\mathcal{W}}$, e.g. $H^{-s}(\mathcal{S}, \R^3)$ with $s > \frac{3}{2}$. Then, we define the noise through a deterministic time dependent correlation operator $\sigma_t$: let $\hat{\sigma} : [0,T] \rightarrow L^2(\mathcal{S}^2,\R^3)$ be a bounded symmetric kernel and define
\begin{equation*}
    (\sigma_t f)(x) = \int_{\mathcal{S}} \hat{\sigma}(x,y,t) f(y) dy, \quad \forall f \in \mathcal{W}.
\end{equation*}
With this definition, $\sigma_t$ is a Hilbert-Schmidt operator mapping $\mathcal{W}$ into itself, so that the noise can be written as
\begin{equation*}
    \sigma_t W_t = \sum_{i \in \mathbbm{N}} \hat{\beta}_t^i \sigma_t e_i,
\end{equation*}
as this series converges in $\mathcal{W}$ almost surely, and in $L^p(\Omega, \mathcal{W})$ for all $p \in \mathbbm{N}$. Here we interpret $\mathcal{W}$ as the space carrying the Wiener process $W_t$, while the notation $L^2(\mathcal{S},\R^3)$ is kept for denoting the space of tridimensional velocities. Moreover, there exists a Hilbert basis $(\phi_n)_n$ consisting of eigenfunctions of the noise operator $\sigma_t$. Here, we scale these eigenfunctions by their corresponding eigenvalues. By a change of basis, there also exists a sequence of standard Brownian motions $(\beta_t^k)_k$, defined on the same filtered space, such that
\begin{equation*}
    \sigma_t W_t = \sum_{i \in \mathbbm{N}} \beta_t^k \phi_k.
\end{equation*}
Furthermore, we can associate a covariance tensor to the random field $\sigma W_t$: if $x,y \in \mathcal{S}$ are two space points, and $t,s \in \R^+$ are two time points, define $Q$ formally by 
\begin{equation*}
    Q(x,y,t,s) =  \E[(\sigma_t dW_t)(x) (\sigma_s dW_s)(y)] = \int_{\mathcal{S}} \hat{\sigma}(x,z,t) \hat{\sigma}(z,y,s) dz \delta(t-s).
\end{equation*}
The diagonal part of this covariance tensor is referred to as the variance tensor, and is denoted by
\begin{equation}
    a(x,t) = \int_{\mathcal{S}} \hat{\sigma}(x,y,t) \hat{\sigma}(y,x,t) dy = \sum_{k=0}^\infty \phi_k(x,t) \phi_k(x,t)^\tr \in \R^{3 \times 3}.
\end{equation}
Moreover, the variance tensor $a$ is assumed to fulfil $a \in L^1([0,T], L^2(\mathcal{S}, \R^{ 3\times 3}))$, which will be enforced below by the regularity conditions on $\sigma$ -- see equation \eqref{smoothness-noise}. Let us note that, in full generality -- that is when $\hat{\sigma}_t$ is itself a random function -- the random-matrix process $a$ is subject to an integrability condition,
\begin{equation*}
    \E \int_0^T \|a(\bcdot,t)\|_{L^2(\mathcal{S}, \R^{ 3\times 3})}^2 < \infty,
\end{equation*}
where $\| \bcdot \|_{L^2(\mathcal{S}, \R^{ 3\times 3}) }$ is the Hilbert norm associated to $L^2(\mathcal{S}, \R^{ 3\times 3})$, the matrix space $ \R^{ 3\times 3}$ being equipped with the Frobenius norm. As such, the integral $\int_0^t \sigma_s dW_s $ is a $\mathcal{W}$-valued Gaussian process with expectation zero and bounded variance: $\E\big[\|\int_0^t \sigma_s dW_s\|_{L^2}^2\big]<\infty$. The quadratic variation of $\int_0^t (\sigma_t d W_s)(x)$ is given by the bounded variation process $\int_0^t a(x,s) ds.$

Similarly to the classical derivation of the Navier-Stokes equations, we may derive the LU Navier-Stokes equations using a stochastic version of the Reynolds Transport Theorem (SRTT) \cite{Mémin_2014}. Let $q$ be a random scalar, within a volume $\mathcal{V}(t)$ transported by the flow. Then, for incompressible unresolved flows -- that is $\nabla \bcdot \sigma_t=0$ -- the SRTT reads
\begin{gather}
    d\Big( \int_{\mathcal{V}(t)} q(x,t) dx \Big) = \int_{\mathcal{V}(t)} \big(\mathbbm{D}_t q + q \nabla \bcdot (u-u_s) dt\big) dx,\\
    \mathbbm{D}_t q = d_t q + (u - u_s) \bcdot \nabla q \: dt + \sigma dW_t \bcdot \nabla q - \frac{1}{2} \nabla \bcdot (a \nabla q) dt, \label{transport-operator}
\end{gather}
where an additional drift $u_s = \frac{1}{2} \nabla \bcdot a$, coined as the It\={o}-Stokes drift in \cite{BCCLM_2020}, is involved. Here, $d_t q(x,t) = q(x,t+dt) - q(x,t)$ is the forward time increment at a fixed spatial point $x$, and $\mathbbm{D}_t q$ is a stochastic transport operator introduced in \cite{Mémin_2014, RMC_2017_Pt1}, which plays the role of the material derivative. The It\={o}-Stokes drift is directly related to the divergence of the variance tensor $a$, which represents the effects of noise inhomogeneity on the large-scale dynamics. Such advection terms are commonly added as corrective terms in ocean large-scale simulation to account for surface waves and Langmuir turbulence \cite{CL_1976, GOP_2004, MSM_1997}. As shown in \cite{BCCLM_2020}, the LU framework naturally exhibits similar features, generalizing the effects of the small-scale inhomogeneity on the large-scale flow.

The stochastic transport operator includes physically interpretable terms for large-scale representation of flows. The last term on the right-hand side of \eqref{transport-operator} is an inhomogeneous diffusion term, representing small-scale mixing. This stochastic diffusion is entirely defined by the variance of the noise, and can be interpreted as a matrix generalization of the Boussinesq eddy viscosity assumption. The third term on the right-hand side represents the transport of the large-scale quantity $q$ by the unresolved velocity. Remarkably, the energy associated with this backscattering term is exactly compensated by the stochastic diffusion term \cite{RMC_2017_Pt1}. This equilibrium can be interpreted as an immediate instance of the fluctuation-dissipation theorem.

By interpreting the stochastic diffusion term $- \frac{1}{2} \nabla \bcdot (a \nabla q) dt$ as the It\={o}-Stratonovitch correction of the back-scattering It\={o} noise term $\sigma dW_t \bcdot \nabla q$ in the context of ``pure'' transports -- i.e. $\mathbbm{D}_t q = 0$ -- we may define a Stratonovitch transport operator
$$\mathbbm{D}^{\circ}_t q := d_t q + (u - u_s) \bcdot \nabla q \: dt + \sigma dW_t \circ \nabla q.$$
Here, $\sigma dW_t \circ \nabla q$ denotes a Stratonovitch noise term. Assuming this noise term is well-defined, the two transport operators \emph{are} equivalent in the absence of forcing terms. By this, we mean that $(\mathbbm{D}_t q = 0)$ and $(\mathbbm{D}^{\circ}_t q = 0)$ are equivalent. The importance of distinguishing these transport operators will appear clearer when exploring weak hydrostatic assumptions.

\subsection{The LU primitive equations with strong hydrostatic assumption}

Let us now derive the LU primitive equations. For modelling purposes, we assume that the flow is isochoric with constant material density. In addition, we suppose that the noise is divergence-free, with a divergence-free corresponding It\={o}-Stokes drift, i.e.
\begin{equation}
    \nabla \bcdot u = \nabla \bcdot u_s =0, \quad \nabla \bcdot \sigma_t dW_t =0.
\end{equation}
From these two assumptions, one can deduce immediately that, for any conservative scalar quantity $q$,
$$\mathbbm{D}_t q =0.$$
A linear law of state, relating density, salinity and temperature, can be expressed through a Taylor expansion: write the density $\rho$ as
\begin{equation}
    \rho = \rho_0\Big(1 + \beta_T (T-T_r) + \beta_S (S-S_r)\Big),
\end{equation}
with $\rho_0$ the reference density of the ocean at a typical temperature $T_r$ and salinity $S_r$. We assume that the thermodynamic parameters $\beta_T := \frac{1}{\rho_0} \frac{\partial \rho}{\partial T}$ and $\beta_S := \frac{1}{\rho_0} \frac{\partial \rho}{\partial T}$ are constant.

Following the derivation of \cite{DHM_2023}, we derive the following stochastic equations of motion by applying the SRTT to the conservation of momentum principle in rotating frame,
\begin{equation}
    \mathbbm{D}_t u + f k \times (u \: dt + \sigma dW_t) = -\frac{1}{\rho_0} \nabla (p \: dt + dp_t^\sigma) - \Av (u \: dt + \sigma dW_t), \label{NS-u}
\end{equation}
where $A^v$ is a diffusion operator defined below -- see equation \eqref{eq-def-Ai}. Notice that we have introduced, in addition to the classical pressure term $p \: dt$, a martingale noise pressure term $dp_t^\sigma$ arising from the stochastic modelling. Importantly, this martingale pressure term must be interpreted in the It\={o} formalism, as follows
\begin{equation}
    dp_t^\sigma = \sum_{k=0}^\infty \pi_k \: d\beta_k,
\end{equation}
where $(\pi_k)_k$ are $\R$-valued functions, that depend implicitly on the \emph{semimartingale} $u$ through the divergence-free condition. Upon applying the SRTT to the conservation of energy and saline mass, we also obtain evolution equations on the temperature and salinity,
\begin{align}
    \mathbbm{D}_t T = - \AT T, \quad \mathbbm{D}_t S = - \AS S. \label{NS-TS}
\end{align}
In equations \eqref{NS-u} and \eqref{NS-TS}, we have employed the anisotropic diffusion operators $A^v$, $A^T$ and $A^S$. For $i \in \{v,T,S\}$, and given viscosities $\mu_i, \nu_i$ specified \emph{a priori}, we define
\begin{align}
    A^{\scriptscriptstyle i} &= -\mu_i (\partial_{xx} + \partial_{yy}) - \nu_i  \partial_{zz}. \label{eq-def-Ai}
\end{align}
Denote by $u^* = u - u_s$, and notice that by assumption $u^*$, $u$ and $u_s$ are divergence-free. Therefore, if we denote by $v^*$, $v$ and $v_s$ their respective horizontal components, we can express the corresponding vertical components using the integro-differential operator
\begin{equation}
    w(v) = \int_z^0 \nabla_H \bcdot v, \label{w-definition}
\end{equation}
under the hypothesis that $w(v)=0$ when $z = 0$. In particular, this boundary condition will be used to define the function spaces introduced in subsection \ref{subsec-def-spaces}. The horizontal gradient operator is denoted by $\nabla_H = (\partial_x \: \partial_y)^\tr$, and we define horizontal Laplace operator as $\Delta_H(\bcdot) = \nabla_H \bcdot (\nabla_H (\bcdot) )$. Additionally, let $\sigma^H dW_t$ and $\sigma^z dW_t$ denote the horizontal and vertical components of $\sigma dW_t$, respectively. Thus, the horizontal and vertical momentum equations write
\begin{gather}
    \mathbbm{D}_t v + \Gamma (v \: dt + \sigma^H dW_t) =  -\Av (v\:dt  + \sigma^H dW_t) - \frac{1}{\rho_0} \nabla_H (p \: dt + dp_t^\sigma), \label{NS-v}\\
    \mathbbm{D}_t w = -\Av (w\:dt  + \sigma^z dW_t) - \frac{1}{\rho_0} \partial_z (p \: dt + dp_t^\sigma) - \frac{\rho}{\rho_0} g dt, \label{NS-w}
\end{gather}
where $\Gamma((a \quad b)^\tr) = f (-b \quad a)^\tr$ stands for the horizontal projection of the Coriolis term.

As explored in a recent paper \cite{AHHS_2022}, a simplified system can be obtained considering an assumption similar to the classical deterministic hypothesis. Specifically, assuming that the vertical acceleration is negligible compared to gravitational one $g$, we have
$$\mathbbm{D}_t w + \Av (w\:dt  + \sigma^z dW_t) \ll g.$$
Under this hypothesis, the vertical momentum equation boils down to
\begin{equation}
    \partial_z p + \rho g = 0, \quad \text{and} \quad \partial_z dp_t^\sigma = 0. \label{strong-hydro}
\end{equation}
We refer to this assumption as the \emph{strong hydrostatic hypothesis}. The validity of the hydrostatic balance corresponds to a regime of small ratio $\frac{\epsilon^2}{Ri}$, where $\epsilon^2=h^2/L^2$ is the squared aspect ratio, with $h$ and $L$ denoting the vertical and horizontal length scales, respectively. The Richardson number is defined as $Ri = N^2/(\partial_z v)^2$, where $N^2= - \frac{g}{\rho_0} \partial_z \rho $ is the stratification factor given by the Brunt-Väsäilä frequency, and $(\partial_z v)^2$ stands for the squared vertical shear of the horizontal velocity \cite{MHPA_1997}. In the stochastic setting, the strong hydrostatic balance holds if the noise does not disrupt this regime (see also Remark \ref{big-remark}).

Gathering all the points described previously and assuming that the strong hydrostatic hypothesis holds, we eventually obtain the following problem
\begin{gather*}
    \mathbbm{D}_t v + \Gamma (v \: dt + \sigma^H dW_t) =  -\Av (v\:dt  + \sigma^H dW_t) - \frac{1}{\rho_0} \nabla_H (p \: dt + dp_t^\sigma),\\
    \mathbbm{D}_t T = - \AT T dt, \quad
    \mathbbm{D}_t S = - \AS S dt,\\
    \nabla_H \bcdot v + \partial_z w =0,\\
    \partial_z p + \rho g = 0, \quad 
    \partial_z dp_t^\sigma = 0,\\
    \rho = \rho_0(1 + \beta_T (T-T_r) + \beta_S (S-S_r)).    
\end{gather*}
Interpreting the stochastic transport operator as a material derivative, this system enjoys a similar structure as the deterministic primitive equations system. The model studied in \cite{AHHS_2022} essentially corresponds to the one above, with the stochastic diffusion term $\frac{1}{2} \nabla \bcdot (a \nabla(\bcdot))$ replaced by $\nu_\sigma \Delta (\bcdot)$, where $\nu_\sigma>0$ is a constant. The authors prove well-posedness results, assuming smooth enough initial conditions, with periodic horizontal boundary conditions and rigid-lid type vertical boundary conditions. This corresponds to a \emph{water world} configuration, also referred to as an \emph{aqua planet}. The main results in \cite{AHHS_2022} are local-in-time well-posedness of the model, and global-in-time well-posedness when the horizontal component of the noise is barotropic, i.e. independent of the vertical coordinate $z$.

With this in mind, we aim to prove similar well-posedness results for more general models, where we make a different assumption on the vertical momentum equation, retaining essentially more stochastic terms. This assumption is referred to as the \emph{weak hydrostatic hypothesis} in the following. These relaxed hydrostatic equilibria correspond to dynamical regimes at the limit of validity of the deterministic hydrostatic assumption, such as high-resolution, non-hydrostatic physical phenomena, including wind- or buoyancy-driven turbulence and deep oceanic convection, where strong enough noise disrupts the strong hydrostatic regime.

\subsection{Transitioning from strong to weak hydrostatic hypothesis} \label{section-strong-to-weak}

To derive our new models, we aim to relax the hydrostatic assumption \eqref{strong-hydro} and derive a weaker form. For this purpose, we begin with physical remarks on the scaling of the quantities involved in the Navier-Stokes vertical momentum equation.

\subsubsection{Scaling the vertical momentum equation}
Remind that the vertical velocity $w$ fulfils the following,
\begin{align}
    \mathbbm{D}_t w = -\Av (w\:dt  + \sigma^z dW_t) - \frac{1}{\rho_0} \partial_z (p \: dt + dp_t^\sigma) - \frac{\rho}{\rho_0} g dt. \label{momentum-w-reminder}
\end{align}
Writing the previous equation \eqref{momentum-w-reminder} in Stratonovitch form yields
\begin{multline}
    \mathbbm{D}_t^\circ w - \frac{1}{2 \rho_0} \sum_k \phi_k \bcdot \nabla \partial_z \pi_k \: dt - \frac{1}{2Re} \sum_k \phi_k \bcdot \nabla A^v \phi_k^z \: dt \\
    = -\Av (w\:dt  + \sigma^z dW_t) - \frac{1}{\rho_0} \partial_z (p \: dt + dp_t^\sigma) - \frac{\rho}{\rho_0} g dt \label{momentum-w-Strato},
\end{multline}
where $\frac{1}{2} \sum_k \phi_k \bcdot \nabla \partial_z \pi_k \: dt$ and $\frac{1}{2Re} \sum_k \phi_k \bcdot \nabla A^v \phi_k^z \: dt$ are additional It\={o}-Stratonovitch correction terms associated to covariation between the noise term and the martingale pressure and molecular dissipation, respectively. Remind that $(\phi_k)$ stands for the eigenfunctions of the noise operator $\sigma$, and that the functions $(\pi_k)_k$ define the martingale pressure term
$$dp_t^\sigma = \sum_{k=0}^\infty \pi_k d\beta_t^k.$$

\noindent
We denote by $U,W$ the typical horizontal and vertical large-scale velocities respectively. By this we mean that the horizontal and vertical velocities $(v,w)$ can be expressed as $(v,w) = (U \tilde{v}, W \tilde{w})$, where $\tilde{v}$ and $\tilde{w}$ have order one. In addition, we define $\tau$, the typical large-scale motion timescale, and $L,H$ the typical (horizontal and vertical) length scales for the dynamics of the ocean. Consequently,
$$U = \frac{L}{\tau} \text{ and } W = \frac{H}{\tau}.$$

\noindent
Moreover, we define $L_\sigma, H_\sigma$ the unresolved characteristic lengths of the small-scale dynamics. Here, $L_\sigma$ and $H_\sigma$ correspond to the mixing length of the unresolved small-scale dynamics. They are not directly related to the size of the large-scale length scales nor to the grid resolution. Define also $\Upsilon_H, \Upsilon_z$ the typical scales of the variance tensor components $a_H$ and $a_z$, respectively. Their physical units are $m^2/s$, that is the unit of a kinematic viscosity. Hence, since the unresolved small scale is modelled by the noise term $\sigma dW_t$, we deduce the following scaled lengths and velocities,
\begin{align}
    L_\sigma = (\Upsilon_H \,\tau)^{1/2} \text{ and }  H_\sigma = (\Upsilon_z \,\tau)^{1/2}. \label{eq:derivation-from-Ito-begin}
\end{align}
Denote also $U_S$ and $W_S$ the It\={o}-Stokes drift counterparts of $U$ and $W$. As the It\={o}-Stokes drift $u_s = \frac{1}{2} (\nabla \bcdot a)$ is associated to a bounded variation dynamics, its characteristic time is $\tau$, and
\begin{gather}
    U_S = \frac{\Upsilon_H}{L} = \frac{L_\sigma^2}{L^2} U, \quad W_S = \frac{\Upsilon_z}{H} = \frac{H_\sigma^2}{H^2} W. \label{IS-scaling}
\end{gather}
Thus we can define $L_S := U_S \,\tau$ and $H_S := W_S \,\tau$. Moreover, we define $\epsilon = \frac{H}{L}$ and $\epsilon_\sigma = \frac{H_\sigma}{L_\sigma}$ the large scale and small scale aspect ratios, respectively. Thus, the large-scale divergence-free conditions read $\frac{U}{L} = \frac{W}{H}$ and $\frac{U_S}{L} = \frac{W_S}{H}$, that is
\begin{equation}
    \frac{W}{U} = \frac{W_S}{U_S} = \epsilon.
\end{equation}
Using the definition of the It\={o}-Stokes drift scaling \eqref{IS-scaling}, we also infer $\frac{W_S}{U_S} = \epsilon_\sigma^2 \epsilon^{-1}$, that is $\epsilon_\sigma = \epsilon$. In particular, we may define
\begin{equation}
    r_\sigma := \frac{H_\sigma}{H} = \frac{L_\sigma}{L} = \frac{H_S}{H} = \frac{L_S}{L} = \frac{(\Upsilon_H \,\tau)^{1/2}}{L}.
\end{equation}
Notice that a direct scaling on the divergence-free condition of the correlation tensor leads coherently to $\Upsilon_z^{1/2}/\Upsilon_H^{1/2}= \epsilon$. With these elements, equation \eqref{momentum-w-reminder} reads
\begin{multline}
    \frac{W}{\tau} d_{\tilde{t}} \tilde{w} + \frac{W}{\tau} (\tilde{u} - r_\sigma^2 \tilde{u}_s) \bcdot \tilde{\nabla} \tilde{w} \: d\tilde{t} + \frac{W}{\tau} r_\sigma \Big(\tilde{\sigma} dW_{\tilde{t}} \bcdot \tilde{\nabla}\Big) \tilde{w} - \frac{W}{2 \tau} r_\sigma^2 \tilde{\nabla} \bcdot (\tilde{a} \tilde{\nabla} \tilde{w}) d\tilde{t} \\- \frac{\mu W}{L^2} (\Delta_{\tilde{H}} + \frac{\nu}{\mu \epsilon^2} \partial_{\tilde{z}\tilde{z}}) (\tilde{w} \: d\tilde{t} + r_\sigma \tilde{\sigma}_z dW_{\tilde{t}}) = - \frac{P}{\rho_0 H} \partial_{\tilde{z}} \tilde{p} \: d\tilde{t} - \frac{P_\sigma}{\rho_0 H} \partial_{\tilde{z}} d\tilde{p}_{\tilde{t}}^\sigma - g\tilde{\rho} d\tilde{t}, \label{scaled-vertical-equation}
\end{multline}
where the tilde notations stand for adimensioned variables and operators. Also, we write $P = \rho_0 g H$ using the classical hydrostatic assumption, and $P_\sigma$ the scaling of the martingale pressure term. Furthermore, the term $d_t w$ is expected to vanish in the primitive equations in the limit of small-aspect ratio, with rate $\epsilon^2$ -- see \cite{LT2019}. Thus, disregarding the terms $d_t w$, we infer that
\begin{equation*}
    \partial_{z} dp_t^\sigma = - \rho_0 \sigma dW_t \bcdot \nabla w  + (\mu \Delta_H + \nu \partial_{zz}) \sigma_z dW_t + O(\epsilon^2).
\end{equation*}
Consequently, $\partial_z dp_t^\sigma$ and $\sigma dW_t \bcdot \nabla w$ share the same order of magnitude, since the viscous term is dominated by the advection term in turbulent flows. Hence, dividing equation \eqref{scaled-vertical-equation} by $\frac{W}{\tau}$ yields
\begin{multline}
    d_{\tilde{t}} \tilde{w} + (\tilde{u} - r_\sigma^2 \tilde{u}_s) \bcdot \tilde{\nabla} \tilde{w} \: d\tilde{t} + r_\sigma \Big(\tilde{\sigma} dW_{\tilde{t}} \bcdot \tilde{\nabla}\Big) \tilde{w} - \frac{1}{2} r_\sigma^2 \tilde{\nabla} \bcdot (\tilde{a} \tilde{\nabla} \tilde{w}) d\tilde{t} \\- \frac{\mu}{LU} (\Delta_{\tilde{H}} + \frac{\nu}{\mu \epsilon^2} \partial_{\tilde{z}\tilde{z}}) (\tilde{w} \: d\tilde{t} + r_\sigma \tilde{\sigma}_z dW_{\tilde{t}}) =  - \frac{1}{\epsilon^2} \frac{gH}{U^2} (\partial_{\tilde{z}} \tilde{p} \: d\tilde{t} + \tilde{\rho} d\tilde{t}) - r_\sigma \partial_{\tilde{z}} d\tilde{p}_{\tilde{t}}^\sigma.
\end{multline}
Introducing the Reynolds and Froude numbers defined as $Re = \frac{LU}{\mu}$ and $Fr^2 = \frac{U^2}{gH}$, respectively, and assuming for simplicity that $\nu \sim \mu \epsilon^2$ (see \cite{LT2019}), we obtain
\begin{multline}
    d_{\tilde{t}} \tilde{w} + (\tilde{u} - r_\sigma^2 \tilde{u}_s) \bcdot \tilde{\nabla} \tilde{w} \: d\tilde{t} + r_\sigma \Big(\tilde{\sigma} dW_{\tilde{t}} \bcdot \tilde{\nabla}\Big) \tilde{w} - \frac{1}{2} r_\sigma^2 \tilde{\nabla} \bcdot (\tilde{a} \tilde{\nabla} \tilde{w}) d\tilde{t} \\- \frac{1}{Re} \tilde{\Delta} (\tilde{w} \: d\tilde{t} + r_\sigma \tilde{\sigma}_z dW_{\tilde{t}}) =  - \frac{Fr^{-2}}{\epsilon^2} (\partial_{\tilde{z}} \tilde{p} \: d\tilde{t} + \tilde{\rho} d\tilde{t}) - r_\sigma \partial_{\tilde{z}} d\tilde{p}_{\tilde{t}}^\sigma. \label{pre-hydro-asymptotic-Ito}
\end{multline}
\modif{Starting from equation \eqref{momentum-w-Strato}, we can adapt the reasoning from equations \eqref{eq:derivation-from-Ito-begin} to \eqref{pre-hydro-asymptotic-Ito}, which yields}
\begin{multline}
    d_{\tilde{t}} \tilde{w} + (\tilde{u} - r_\sigma^2 \tilde{u}_s) \bcdot \tilde{\nabla} \tilde{w} \: d\tilde{t} + r_\sigma \tilde{\sigma} dW_{\tilde{t}} \circ \tilde{\nabla} \tilde{w} + \frac{1}{2}  r_\sigma^2 \sum_k (\phi_k \bcdot \tilde{\nabla}) \partial_{\tilde{z}} \tilde{\pi}_k \: d\tilde{t} - \frac{1}{Re} \tilde{\Delta} (\tilde{w} \: d\tilde{t} + r_\sigma \tilde{\sigma}_z dW_{\tilde{t}})\\
    + \frac{1}{2 Re} r_\sigma^2 \sum_k (\phi_k \bcdot \tilde{\nabla})\tilde{\Delta} \tilde{\phi}_k^z \: dt = - \frac{Fr^{-2}}{\epsilon^2} (\partial_{\tilde{z}} \tilde{p} \: d\tilde{t} + \tilde{\rho} d\tilde{t}) -  r_\sigma \partial_{\tilde{z}} d\tilde{p}_{\tilde{t}}^\sigma. \label{pre-hydro-asymptotic-Strato}
\end{multline}
Hence, multiplying by $\epsilon^2$ and taking $\epsilon \rightarrow 0$ in equation \eqref{pre-hydro-asymptotic-Ito} -- It\={o} form equation -- yields,
\begin{multline}
    - r_\sigma \epsilon^2 \Bigg[ \tilde{\sigma} dW_{\tilde{t}} \bcdot \tilde{\nabla} \tilde{w} - \frac{1}{Re} \tilde{\Delta}(\tilde{\sigma}_z dW_{\tilde{t}})\Bigg] + r_\sigma^2 \epsilon^2 \Bigg[\frac{1}{2} \tilde{\nabla} (\tilde{a} \tilde{\nabla} \tilde{w}) + \tilde{u}_s \bcdot \tilde{\nabla} \tilde{w} \Bigg] d\tilde{t}\\
    = \Big[Fr^{-2}(\partial_{\tilde{z}} \tilde{p} + \tilde{\rho}) d\tilde{t} + r_\sigma \epsilon^2  \partial_{\tilde{z}} d\tilde{p}_{\tilde{t}}^\sigma\Big] + O(\epsilon^2). \label{hydro-asymptotic-Ito}
\end{multline}
\modif{Multiplying and passing to the limit similarly in equation \eqref{pre-hydro-asymptotic-Strato}} -- Stratonovitch form equation -- yields,
\begin{multline}
    - \epsilon^2 \Big(r_\sigma \tilde{\sigma} dW_{\tilde{t}} \circ \tilde{\nabla} \tilde{w}\Big) + r_\sigma \epsilon^2 \frac{1}{Re}\tilde{\Delta}(\tilde{\sigma}_z dW_{\tilde{t}}) + r_\sigma^2 \epsilon^2 \tilde{u}_s \bcdot \tilde{\nabla} \tilde{w} d\tilde{t} + r_\sigma^2 \epsilon^2  \frac{1}{2}\sum_k (\phi_k \bcdot \tilde{\nabla}) \partial_{\tilde{z}} \tilde{\pi}_k \: d\tilde{t} \\
    + \frac{1}{2Re} r_\sigma^2 \epsilon^2 \sum_k (\tilde{\phi}_k \bcdot \tilde{\nabla}) \tilde{\Delta} \tilde{\phi}_k^z \: d\tilde{t} = \Big[Fr^{-2} (\partial_{\tilde{z}} \tilde{p} + \tilde{\rho}) d\tilde{t} + r_\sigma \epsilon^2 \partial_{\tilde{z}} d\tilde{p}_{\tilde{t}}^\sigma\Big] + O(\epsilon^2). \label{hydro-asymptotic-Strato}
\end{multline}
Notice that two new correction terms emerge when transitioning to the Stratonovitch formulation. Additionally,  the evolution term $d_t w$ is implicitly included within the asymptotic notation $O(\epsilon^2)$. The previous computation will allow to derive two distinct expressions of the martingale pressure term $p^\sigma$ for a small yet non-zero large-scale aspect ratio $\epsilon$. For this purpose, let us discuss the scaling of $r_\sigma$. Reminding that $r_\sigma = \frac{(\Upsilon_H \,\tau)^{1/2}}{L}$, the condition $r_\sigma \gg 1$ corresponds to the case of a much greater mixing length than the large scale characteristic length, with the implicit assumption that the unresolved small-scale processes correspond to fast phenomena. Another interpretation is possible by introducing a mesoscopic Reynolds number $Re_\sigma := \frac{UL}{\Upsilon_H}$. As $r_\sigma^2 = \frac{1}{Re_\sigma}$, the condition $r_\sigma^2 \gg 1$ is equivalent to $Re_\sigma \ll 1$, which implies a larger variance tensor scaling than the intrinsic large-scale eddy viscosity $UL$. 

In such case, the terms of order $r_\sigma \epsilon^2$ \emph{are not negligible} compared to those of order $\epsilon^2$. Consequently, disregarding the terms order $O(\epsilon^2)$, equations \eqref{hydro-asymptotic-Ito} and \eqref{hydro-asymptotic-Strato} yield two distinct relations, respectively
\begin{multline}
    - r_\sigma \epsilon^2 \Bigg[ (\tilde{\sigma} dW_{\tilde{t}} \bcdot \tilde{\nabla} \tilde{w}) - \frac{1}{Re}\tilde{\Delta}(\tilde{\sigma}_z dW_{\tilde{t}})\Bigg] + r_\sigma^2 \epsilon^2 \Bigg[\frac{1}{2} \tilde{\nabla} \bcdot (\tilde{a} \tilde{\nabla} \tilde{w}) + \tilde{u}_s \bcdot \tilde{\nabla} \tilde{w} \Bigg] d\tilde{t}\\
    = \Big[ Fr^{-2}(\partial_{\tilde{z}} \tilde{p} + \tilde{\rho}) d\tilde{t} + r_\sigma \epsilon^2  \partial_{\tilde{z}} d\tilde{p}_{\tilde{t}}^\sigma \Big], \label{weak-hydro-eps-order1-Ito}
\end{multline}
and
\begin{multline}
    - \epsilon^2 \Bigg[ (r_\sigma \tilde{\sigma} dW_{\tilde{t}} \circ \tilde{\nabla} \tilde{w}) - \frac{1}{Re}\tilde{\Delta}(\tilde{\sigma}_z dW_{\tilde{t}})\Bigg] + r_\sigma^2 \epsilon^2 \tilde{u}_s \bcdot \tilde{\nabla} \tilde{w} d\tilde{t} + r_\sigma^2 \epsilon^2 \frac{1}{2}\sum_k (\tilde{\phi}_k \bcdot \tilde{\nabla}) \partial_{\tilde{z}} \tilde{\pi}_k \: d\tilde{t} \\
    + \frac{1}{2Re} r_\sigma^2 \epsilon^2 \sum_k (\tilde{\phi}_k \bcdot \tilde{\nabla}) \tilde{\Delta} \tilde{\phi}_k^z \: d\tilde{t} = \Big[ Fr^{-2}(\partial_{\tilde{z}} \tilde{p} + \tilde{\rho}) d\tilde{t} + r_\sigma \epsilon^2  \partial_{\tilde{z}} d\tilde{p}_{\tilde{t}}^\sigma \Big]. \label{weak-hydro-eps-order1-Strato}
\end{multline}
As the evolution term $d_t w$ has been neglected, these two relations are no longer equivalent. On the one hand, equation \eqref{weak-hydro-eps-order1-Ito} provides explicit expression for both the bounded variation and martingale components of the pressure term. On the other hand, equation \eqref{weak-hydro-eps-order1-Strato} involves a Stratonovitch transport noise $\Big(r_\sigma \tilde{\sigma} dW_{\tilde{t}} \circ \tilde{\nabla} \tilde{w}\Big)$, which must be interpreted as a function of $v$. This is because $w = w(v) = \int_z^0 \nabla_H \bcdot v \: dz'$ as imposed by the divergence-free condition. Consequently, the two formulations differ due to the It\={o}-Stratonovitch correction terms, which are anticipated to be more intricate in the Stratonovitch formulation \eqref{weak-hydro-eps-order1-Strato}.

Furthermore, observe that the deterministic hydrostatic hypothesis, $\partial_{\tilde{z}} \tilde{p} + \tilde{\rho} = 0$, is, in fact, a \emph{zeroth-order} approximation in $\epsilon$ whenever $r_\sigma \epsilon^2 = o(1)$. This becomes evident as equations \eqref{hydro-asymptotic-Ito} and \eqref{hydro-asymptotic-Strato} independently yield
\begin{equation}
    (\partial_{\tilde{z}} \tilde{p} + \tilde{\rho}) d\tilde{t} = O(\epsilon^2, r_\sigma \epsilon^2).
\end{equation}
In addition, the assumption $\partial_{z} dp_t^\sigma = 0$ can be justified in the limit of infinitely small aspect ratios by the following argument. As mentioned earlier, by retaining only the martingale terms in either equation \eqref{weak-hydro-eps-order1-Ito} or \eqref{weak-hydro-eps-order1-Strato}, we infer
\begin{equation*}
    \partial_{z} dp_t^\sigma \sim - \rho_0 \sigma dW_t \bcdot \nabla w,
\end{equation*}
that is the back-scattering advection term of $w$ and the vertical pressure gradient have the same order of magnitude. Thus, we deduce that
\begin{equation*}
    dp_t^\sigma(x,y,z) \sim dp_t^{\sigma,S}(x,y) + \rho_0 \int_z^0 (\sigma dW_t (x,y,z) \bcdot \nabla) w(x,y,z) \: dz',
\end{equation*}
where $p^{\sigma,S}$ is a small-scale surface pressure, which is independent of the vertical coordinate $z$. Consequently, the pressure scaling can be expressed as
\begin{equation}
    P_\sigma \sim P_\sigma^S + \rho_0 r_\sigma H \frac{W}{\tau} = P_\sigma^S + \rho_0 r_\sigma \epsilon^2 U^2 = P_\sigma^S + O(r_\sigma \epsilon^2),\label{scaled-pressure-equation}
\end{equation}
and hence, 
\begin{equation}
    \nabla_{\tilde{H}} P_\sigma \sim \frac{P_\sigma^S}{L} + \rho_0 r_\sigma \epsilon^2 \frac{U^2}{L} = \frac{P_\sigma^S}{L} + O(r_\sigma \epsilon^2).
\end{equation}
This implies that the contribution of the vertical martingale pressure gradient term is negligible compared to its horizontal gradient counterpart. In the limit of infinitely small aspect ratios $\epsilon=0$, the horizontal pressure gradient affecting the horizontal dynamics reduces  to $\nabla_{\tilde{H}} P_\sigma \sim \frac{P_\sigma^S}{L}$, corresponding to neglecting completely $\partial_z dp_t^\sigma$.

To summarise, the approach presented in this subsection can be interpreted as a refinement of the aforementioned zeroth-order approximation of the Navier-Stokes vertical momentum equation in the context of stochastic flows, which is valid when $r_\sigma \gg 1$, or equivalently $Re_\sigma \ll 1$, with the mesoscopic Reynolds number representing the ratio between the intrinsic large-scale eddy viscosity and the variance tensor scaling. 

\subsubsection{Deriving a weak hydrostatic hypothesis formulation} \label{sec:weak-hydro-hyp-derivation}
Based on the previous remarks, we neglect only the large-scale contribution of the vertical acceleration \modif{in the It\={o} form of the vertical momentum equation \eqref{momentum-w-reminder}.} Consequently, this equation becomes
\begin{align}
    - \frac{1}{\rho_0} \partial_z &(p \: dt + dp_t^\sigma) - \frac{\rho}{\rho_0} g dt = \mathbbm{D}_t w + \Av(w\:dt + \sigma^z dW_t) \nonumber \\
    &\approx - u_s \bcdot \nabla w \: dt + \sigma dW_t \bcdot \nabla w - \frac{1}{2} \nabla \bcdot (a \nabla w) dt + \Av(\sigma^z dW_t).
\end{align}
This can be expressed as two separate equations,
\begin{equation}
    \frac{1}{\rho_0} \partial_z p = - \frac{\rho}{\rho_0} g + \frac{1}{2} \nabla \bcdot (a \nabla w) + u_s \bcdot \nabla w, \text{ and } \quad
    \frac{1}{\rho_0} \partial_z dp_t^\sigma = -\sigma dW_t \bcdot \nabla w - \Av(\sigma^z dW_t). \label{weak-Ito}
\end{equation}
We call these relations the \emph{It\={o} weak hydrostatic hypothesis}. Consequently, the \emph{finite variation} pressure gradients are expressed as,
\begin{equation}
    \frac{1}{\rho_0}\nabla_H p = - \nabla_H \int_z^0 \frac{\rho}{\rho_0} g \: dz'  -  \nabla_H \int_z^0 \frac{1}{2} \nabla \bcdot (a \nabla w) \: dz' - \nabla_H \int_z^0 u_s \bcdot \nabla w \: dz' + \frac{1}{\rho_0} \nabla_H p^s.
\end{equation}
Similarly, and the \emph{martingale} pressure gradients read
\begin{equation}
    \frac{1}{\rho_0}\nabla_H dp_t^\sigma =  \nabla_H \int_z^0 \sigma dW_t \bcdot \nabla w + \Av(\sigma^z dW_t) \: dz' + \frac{1}{\rho_0} \nabla_H dp^{\sigma,s}_t.
\end{equation}
Ultimately, the problem under weak hydrostatic balance $(\mathcal{P})$ is formulated as:
\begin{equation} (\mathcal{P}) \left\{\begin{array}{l}
     \mathbbm{D}_t v + \Gamma (v \: dt + \sigma^H dW_t) =  -\Av (v\:dt  + \sigma^H dW_t) - \frac{1}{\rho_0} \nabla_H (p \: dt + dp_t^\sigma),\\
        \mathbbm{D}_t T = - \AT T dt,\\
        \mathbbm{D}_t S = - \AS S dt,\\
        \nabla_H \bcdot v + \partial_z w =0,\\
        \frac{1}{\rho_0}\nabla_H p = \nabla_H \int_z^0 \Big[- \frac{\rho}{\rho_0} g - \frac{1}{2} \nabla \bcdot (a \nabla w) - u_s \bcdot \nabla w \Big] dz' + \frac{1}{\rho_0} \nabla_H p^s,\\
        \frac{1}{\rho_0} \nabla_H dp_t^\sigma = \nabla_H \int_z^0 \Big[ \sigma dW_t \bcdot \nabla w + \Av(\sigma^z dW_t) \Big] dz' + \frac{1}{\rho_0} \nabla_H dp^{\sigma,s}_t,\\
        \rho = \rho_0\Big(1 + \beta_T (T-T_r) + \beta_S (S-S_r)\Big).
\end{array}\right.
\end{equation}
Let us emphasise that the model above relies on an important modelling choice. The stochastic integrals in the transport operator \eqref{transport-operator} are expressed in It\={o} form. One could use the Stratonovitch formalism instead -- i.e. the operator $\mathbbm{D}_t^\circ$ -- leading to the following equation on $w$,
\begin{multline}
    \mathbbm{D}^{\circ}_t w = d_t w + (u - u_s) \bcdot \nabla w \: dt + \sigma dW_t \circ \nabla w = -\Av (w\:dt  + \sigma^z dW_t)\\
    - \frac{1}{\rho_0} \partial_z (p \: dt + dp_t^\sigma) - \frac{\rho}{\rho_0} g dt. \label{NS-w-strato}
\end{multline}
In this case, using the previous scaling argument, the \emph{semimartingale} pressure gradient fulfils
\begin{multline}
    \frac{1}{\rho_0}\nabla_H (p \: dt + dp_t^\sigma) + \nabla_H \int_z^0 \frac{\rho}{\rho_0} g \: dz' \: dt = -\nabla_H \int_z^0 u_s \bcdot \nabla w \: dz' \: dt \\+ \nabla_H \int_z^0 \sigma dW_t \circ \nabla w \: dz' + \nabla_H \int_z^0 \Av(\sigma^z dW_t) \: dz' + \frac{1}{\rho_0} \nabla_H (p^s dt + dp_t^{\sigma,s}).
\end{multline}
We refer to  this relation as the \emph{Stratonovitch weak hydrostatic hypothesis}. Thus, we define another general problem $(\mathcal{P}^{\circ})$, which reads
\begin{equation} (\mathcal{P}^\circ) \left\{\begin{array}{l}
    \mathbbm{D}_t^{\circ} v + \Gamma (v \: dt + \sigma^H dW_t) =  -\Av (v\:dt  + \sigma^H dW_t) - \frac{1}{\rho_0} \nabla_H (p \: dt + dp_t^\sigma),\\
    \mathbbm{D}_t^{\circ} T = - \AT T dt,\\
    \mathbbm{D}_t^{\circ} S = - \AS S dt,\\
    \nabla_H \bcdot v + \partial_z w =0,\\
    \frac{1}{\rho_0}\nabla_H (p \: dt + dp_t^\sigma) = - \nabla_H \int_z^0 \frac{\rho}{\rho_0} g \: dz' \: dt - \nabla_H \int_z^0 u_s \bcdot \nabla w \: dz' \: dt \\
    \qquad \quad + \frac{1}{\rho_0} \nabla_H (p^s dt + dp_t^{\sigma,s}) + \nabla_H \int_z^0 \Av(\sigma^z dW_t) \: dz' + \nabla_H \int_z^0 \sigma dW_t \circ \nabla w \: dz'\\
    \qquad \qquad \qquad - \frac{1}{2} \nabla_H \int_z^0 \sum_k (\phi_k \bcdot \nabla) \partial_z \pi_k \: dz' \: dt - \frac{1}{2} \nabla_H \int_z^0 \sum_k (\phi_k \bcdot \nabla) \Delta \phi_k^z \: dz' \: dt,\\
    \rho = \rho_0\Big(1 + \beta_T (T-T_r) + \beta_S (S-S_r)\Big).    
\end{array}\right.
\end{equation}
In particular, the Stratonovitch transport noise terms appearing in the horizontal velocity dynamics are the following, since the other noise terms are additive,
\begin{multline}
    \sigma dW_t \circ \nabla v + \nabla_H \int_z^0 \sigma dW_t \circ \nabla w \: dz' =  \sigma dW_t \bcdot \nabla v + \nabla_H \int_z^0 \sigma dW_t \bcdot \nabla w \: dz'\\ - \frac{1}{2} \Big[ \nabla \bcdot a \nabla v + \nabla \bcdot \check{a}[\nabla w] +\nabla_H \int_z^0 \nabla \bcdot \hat{a}[\nabla v] dz' + \nabla_H \int_z^0 \nabla \bcdot \hat{\hat{a}} [\nabla w] \: dz' + \mathcal{C}_\sigma \Big] dt,
\end{multline}
where, with double index summation convention,
\begin{align}
    &\text{for $f$ an $\R^3$-valued vector field},\quad  \check{a}f = \sum_{k=0}^\infty (\phi_k \otimes \nabla_H \int_z^0  (f^\tr \phi_k) dz')^\tr, (\in \R^{2 \times 3}) \label{eq:def-Ito-Strato-term} \\
    &\text{for $M$ an $\R^{3\times2}$-valued matrix field}, \quad\hat{a}M = \sum_{k=0}^\infty \phi_k \int_z^0 \nabla_H \bcdot (M^\tr \phi_k) \: dz',  (\in \R^3)\\
    &\text{for $f$ an $\R^3$-valued vector field}, \quad \hat{\hat{a}}f = \sum_{k=0}^\infty \phi_k \Delta_H \int_{z}^{0} \Big(\int_{z'}^{0} (f^\tr \phi_k) dz''\Big) dz', (\in \R^3)
\end{align}
and \begin{equation}
    \mathcal{C}_\sigma = \sum_{k=0}^\infty \phi_k \bcdot \nabla_H (\Gamma \phi_k^H + A^v \phi_k^H) + \nabla_H \int_z^0 \phi_k \bcdot \nabla (w(\Gamma \phi_k^H) + A^v \phi_k^z) dz'. \label{equation-C-sigma}
\end{equation} 
Consequently, different It\={o}-Stratonovitch correction terms arise from this approach, indicating that the two problems differ significantly. The former may result in physical energy imbalance, even though its derivation is ``more direct'' since it originates from the 3D LU Navier-Stokes equations. In contrast, the latter involves more complex terms and cannot be interpreted as a straightforward simplification of the 3D LU Navier-Stokes equations. Nevertheless, introducing the transport noise with its "true" It\={o}-Stratonovich correction ensures \emph{a priori} an energy balance. Both models are analyzed in the following, with a particular focus on $(\mathcal{P})$ rather than $(\mathcal{P}^{\circ})$.

In both cases, a significant difficulty arises when addressing such problems: the transport noise $\sigma dW_t \bcdot \nabla w$ and the It\={o}-Stokes drift advection $u_s \bcdot \nabla w$ in the stochastic pressure term lead to the following contributions in the horizontal velocity dynamics,
$$\nabla_H \int_z^0 \sigma dW_t \bcdot \nabla w \: dz', \quad \nabla_H \int_z^0 u_s \bcdot \nabla w \: dz' \: dt.$$
Establishing a suitable energy estimate for global pathwise existence and uniqueness with such terms remains a challenge, primarily because they involve three derivatives of the horizontal velocity $v$, given that $w(v) = \int_z^0 \nabla_H \bcdot v$.
        %-----------------------------------------
        \iffalse
        One may argue that we could rather consider the following noise and It\={o}-Stokes drift terms, according to our scaling argument,
        \begin{align*}
            \nabla_H \int_z^0 \sigma_z dW_t \partial_z w \: dz' &= - \nabla_H \int_z^0 \sigma_z dW_t (\nabla_H \bcdot v) \: dz',\\
            \nabla_H \int_z^0 w_s \partial_z w \: dz' \:&= - \nabla_H \int_z^0 w_s (\nabla_H \bcdot v) \: dz',
        \end{align*}
        since these involve now only two derivatives on the horizontal velocity $v$. Although these terms look \emph{a priori} simpler to deal with, we do not know how to make energy estimates in this case either. For this reason, we choose to keep the problem as general as possible.
        \fi
        %-----------------------------------------
In the following, we propose two techniques to regularize problem $(\mathcal{P})$, and one technique to regularize problem $(\mathcal{P}^{\circ})$.

\subsubsection{Weak low-pass filtered hydrostatic hypothesis} \label{sec-filtered-problem}

We study the problem $(\mathcal{P})$ first. To enforce greater regularity for the vertical transport noise, we define a regularizing convolution kernel $K$, and replace the It\={o} weak hydrostatic hypothesis \eqref{weak-Ito} by
\begin{gather}
    \frac{1}{\rho_0} \partial_z p = - \frac{\rho}{\rho_0} g + \frac{1}{2} \nabla \bcdot (a^K \nabla w) + K*[u_s \bcdot \nabla w] \: dt, \nonumber\\
    \frac{1}{\rho_0} \partial_z dp_t^\sigma = - K*[\sigma dW_t \bcdot \nabla w] - \Av(\sigma^z dW_t), \label{weak-filtered}
\end{gather}
with 
\begin{equation}
    a^K f = \sum_{k=0}^\infty \phi_k \mathcal{C}_K^* \mathcal{C}_K (\phi_k^\tr f), \label{def-aK}
\end{equation}
where $\mathcal{C}_K$ is the operator $f \longmapsto K*f$. Notice that the regularising kernel only affects the It\={o}-Stokes drift advection, the vertical transport noise and the associated diffusion -- and not possible vertical additive noises. Moreover, the stochastic diffusion operator $\frac{1}{2} \nabla \bcdot (a^K \nabla (\bcdot))$ is chosen to be the covariation correction term associated to $K*[\sigma dW_t \bcdot \nabla (\bcdot))]$. We refer to this assumption as the \emph{(It\={o}) weak (low-pass) filtered hydrostatic hypothesis}. This strategy involves filtering the transport noise of the vertical component and neglecting the vertical acceleration of the resolved vertical velocity. The noise terms, together with the stochastic diffusion term, represent deviations from a strong hydrostatic equilibrium. Convolving the vertical transport noise with $K$ effectively removes its highest frequencies. This new hypothesis should be seen as a relaxation of the strong hydrostatic balance, allowing for the consideration of more general stochastic pressures and extending the range of dynamical regimes compared to the strong hydrostatic case. Furthermore, it is worth noticing that applying a filtering kernel is common practice in defining numerical models for the primitive equations. This technique is also frequently employed to establish the well-posedness of specific (mesoscale) subgrid models, such as the Gent-McWilliams model \cite{KT_2023}.

\modif{
\begin{remark}
    Rather than considering the regularisation above for the pressure equation, one could regularise $w$ the vertical velocity only, that is
    \begin{gather*}
        \frac{1}{\rho_0} \partial_z p = - \frac{\rho}{\rho_0} g + \frac{1}{2} \mathcal{C}_K^*\nabla \bcdot (a \nabla \mathcal{C}_K w) + (u_s \bcdot \nabla) \mathcal{C}_K w] \: dt,\\
        \frac{1}{\rho_0} \partial_z dp_t^\sigma = - [(\sigma dW_t \bcdot \nabla) \mathcal{C}_K w] - \Av(\sigma^z dW_t).
    \end{gather*}
    This regularisation is equally valid, and would yield similar results.
\end{remark}
}

We can compare this approach with that of \cite{AHHS_2022_preprint}, where the authors introduced a temperature noise affecting the pressure equation. In their model, the pressure noise is of \emph{thermodynamic} origin. By contrast, our model considers transport noise in the vertical velocity component, arising from a \emph{mechanical} origin. This approach explicitly retains additional terms dependent on the vertical velocity $w$ and captures the influence of unresolved small-scale velocity (e.g., turbulence or submesoscale components) on the vertical large-scale velocity. As a result, the problem formulation is closer to a ``true'' three-dimensional problem, which formally justifies the use of additional regularization via a filtering kernel. Thus, assuming that this weak hydrostatic hypothesis holds instead of the strong one, the following problem is derived,
\begin{gather}
    \mathbbm{D}_t v + \Gamma (v \: dt + \sigma^H dW_t) =  -\Av (v\:dt  + \sigma^H dW_t) - \frac{1}{\rho_0} \nabla_H (p \: dt + dp_t^\sigma), \label{primitive-v}\\
    \mathbbm{D}_t T = - \AT T dt,\label{primitive-T}\\
    \mathbbm{D}_t S = - \AS S dt,\label{primitive-S}\\
    \nabla_H \bcdot v + \partial_z w =0,\\
    \frac{1}{\rho_0}\partial_z p + \frac{\rho}{\rho_0} g - \frac{1}{2} \nabla \bcdot (a^K \nabla w) - K*[u_s \bcdot \nabla w] = 0,\label{primitive-p}\\
    \frac{1}{\rho_0}\partial_z dp_t^\sigma +K*[\sigma dW_t \bcdot \nabla w ]+ \Av(\sigma^z dW_t) = 0,\label{primitive-dp}\\
    \rho = \rho_0\Big(1 + \beta_T (T-T_r) + \beta_S (S-S_r)\Big).    
\end{gather}
The weak filtered hydrostatic hypothesis impacts the horizontal momentum equation through the horizontal pressure gradients, as follows,
\begin{align*}
    \frac{1}{\rho_0}\nabla_H p &= - \nabla_H \int_z^0 \frac{\rho}{\rho_0} g \: dz'  +  \nabla_H \int_z^0 \frac{1}{2} \nabla \bcdot (a^K \nabla w) \: dz' + \nabla_H \int_z^0 K*[u_s \bcdot \nabla w] \: dz'\\
    & \qquad \qquad \qquad \qquad \qquad \qquad \qquad \qquad \qquad \qquad \qquad \qquad + \frac{1}{\rho_0} \nabla_H p^s,\\
    \frac{1}{\rho_0}\nabla_H dp_t^\sigma &=  -\nabla_H \int_z^0 K*[\sigma dW_t \bcdot \nabla w] - \Av(\sigma^z dW_t) \: dz' + \frac{1}{\rho_0} \nabla_H dp^{\sigma,s}_t.
\end{align*}
While new additive noise terms do not pose any issue for a sufficiently regular $\sigma dW_t$, regularizing the term $\nabla_H \int_z^0 \sigma dW_t \bcdot \nabla w \: dz'$ in the expression of $\frac{1}{\rho_0} \nabla_H dp^{\sigma}_t$ is essential for our theoretical analysis. This term represents the horizontal influence of the vertical transport noise $\sigma dW_t \bcdot \nabla w$. Applying a smoothing filter as a regularization method naturally increases its spatial scale, resulting in the vertical transport noise having a spatial scale larger than the resolution scale. Similarly, the term $\nabla_H \int_z^0 \nabla \bcdot (a \nabla w) \: dz'$ corresponds to the horizontal influence of the covariation correction arising from the LU Navier-Stokes equations, and has been modified accordingly. The well-posedness of this model, under appropriate regularity and structural conditions, constitutes our main result.

Additionally, we introduce two alternative methods to regularize the vertical dynamics, this time incorporating (hyper)viscosity terms. Such approaches are also commonly employed in numerical implementations of primitive equations. These viscosity-based regularized models are detailed below.

\subsubsection{Weak hydrostatic hypothesis with an additional eddy-(hyper)viscosity}

Another approach to regularize the vertical dynamics of $(\mathcal{P})$ is to identify the large-scale contributions of the transport and the molecular diffusion terms associated to $w$ in \eqref{NS-w} to a (hyper)diffusion term $\alpha (-\Delta)^{\gamma_r} w$, with $\alpha >0$ and $\gamma_r \geq 1$. This leads to the following identity,
\begin{align*}
    &\mathbbm{D}_t w + \Av(w\:dt + \sigma^z dW_t) \\
    &= \underbrace{d_t w + u \bcdot \nabla w \: dt + \Av(w\:dt)}_{\approx -\alpha (-\Delta)^{\gamma_r} w \: dt} - u_s \bcdot \nabla w dt + \sigma dW_t \bcdot \nabla w - \frac{1}{2} \nabla \bcdot (a \nabla w) dt + \Av(\sigma^z dW_t)\\
    &\approx - u_s \bcdot \nabla w \: dt + \sigma dW_t \bcdot \nabla w - \frac{1}{2} \nabla \bcdot (a \nabla w) dt + \Av(\sigma dW_t) + \alpha (-\Delta)^{\gamma_r} w \: dt.
\end{align*}
 Using \eqref{NS-w} a new relation between the vertical velocity and the pressure is derived,
\begin{multline*}
    - u_s \bcdot \nabla w \: dt + \sigma dW_t \bcdot \nabla w - \frac{1}{2} \nabla \bcdot (a \nabla w) dt + \Av(\sigma^z dW_t) +\alpha (-\Delta)^{\gamma_r} w \: dt\\
    = - \frac{1}{\rho_0} \partial_z (p \: dt + dp_t^\sigma) - \frac{\rho}{\rho_0} g dt,
\end{multline*}
which implies
\begin{gather}
    \frac{1}{\rho_0} \partial_z p = - \frac{\rho}{\rho_0} g + \frac{1}{2} \nabla \bcdot (a \nabla w) + u_s \bcdot \nabla w - \alpha (-\Delta)^{\gamma_r} w, \nonumber\\
    \frac{1}{\rho_0} \partial_z dp_t^\sigma = -\sigma dW_t \bcdot \nabla w - \Av(\sigma^z dW_t). \label{weak-eddy}
\end{gather}
We refer to this assumption as the \emph{(It\={o}) weak eddy-(hyper)viscosity hydrostatic hypothesis}. The introduction of a hyper-viscosity term is, again, common practice for improving the regularity and stability of ocean numerical models. In this paper, we establish the existence of martingale solutions for such systems when $\gamma > 2$.

However, it is important to remark that $\nabla_H \int_z^0 \nabla \bcdot (a \nabla w) \: dz'$ is \emph{not} the covariation correction of the horizontal momentum noise term $\nabla_H \int_z^0 \sigma dW_t \bcdot \nabla w \: dz'$, as mentioned in Subsection \ref{section-strong-to-weak}. For this reason, we investigate how this approach applies to the problem $(\mathcal{P}^{\circ})$.

\subsubsection{Weak hydrostatic hypothesis with energy balanced perturbation}

Doing the same computation as before in the Stratonovitch formalism yields
\begin{gather}
    \mathbbm{D}_t v + \Gamma (v \: dt + \sigma^H dW_t) =  -\Av (v\:dt  + \sigma^H dW_t) - \frac{1}{\rho_0} \nabla_H (p \: dt + dp_t^\sigma) + \frac{1}{2} \nabla \bcdot \check{a}[\nabla w] \: dt \nonumber\\
    \qquad \qquad \qquad \qquad \qquad \qquad \qquad \qquad \qquad \qquad \qquad \qquad + \frac{1}{2} \mathcal{C}_\sigma dt, \label{primitive-v-2}\\
    \mathbbm{D}_t T = - \AT T dt,\label{primitive-T-2}\\
    \mathbbm{D}_t S = - \AS S dt,\label{primitive-S-2}\\
    \nabla_H \bcdot v + \partial_z w =0,\\
    \frac{1}{\rho_0}\partial_z p + \frac{\rho}{\rho_0} g = u_s \bcdot \nabla w + \frac{1}{2} \nabla \bcdot \hat{a}[\nabla v] + \frac{1}{2} \nabla \bcdot \hat{\hat{a}} [\nabla w] - \alpha (-\Delta)^{\gamma_r} w \nonumber\\
    \qquad \qquad \qquad \qquad + \frac{1}{2}\sum_k \phi_k \bcdot \nabla \partial_z \pi_k + \frac{1}{2} \sum_k (\phi_k \bcdot \nabla) \Delta \phi_k^z, \label{primitive-p-2}\\
    \frac{1}{\rho_0}\partial_z dp_t^\sigma +\sigma dW_t \bcdot \nabla w + \Av(\sigma^z dW_t) = 0,\label{primitive-dp-2}\\
    \rho = \rho_0\Big(1 + \beta_T (T-T_r) + \beta_S (S-S_r)\Big),    
\end{gather}
where $\mathcal{C}_\sigma$ has been defined in equation \eqref{equation-C-sigma}. We refer to the equations \eqref{primitive-p-2} and \eqref{primitive-dp-2} as the \emph{weak eddy-viscosity energy-balanced hydrostatic hypothesis}. \modif{We remind that the main difference between this system and the aforementioned one with \emph{(It\={o}) weak eddy-(hyper)viscosity hydrostatic hypothesis} is that ``true'' It\={o}-Stratonovitch correction term is present on the right-hand side of equation \eqref{primitive-v-2}. The existence of two distinct possible formulations arises from the derivation we present in section \ref{sec:weak-hydro-hyp-derivation}. Thus, the model above is \emph{a priori} balanced in terms of the $L^2$-energy.} However, we keep an extra diffusion term $\alpha (-\Delta)^{\gamma_r} w$ on the vertical velocity $w$, which is used to compensate for residual energy terms arising from the use of It\={o}'s lemma. Our main result, that is Theorem \ref{theorem-EB}, is only proven when $\gamma_r>1$. However, establishing the existence of martingale solutions in the (Newtonian) diffusive case $\gamma_r=1$ remains a challenge without additional regularisation.

\subsubsection{Initial and boundary conditions}

In addition, we introduce the initial and boundary conditions we use in our study. Decompose the boundary as $\partial \mathcal{S} = \Gamma_u \cup \Gamma_b \cup \Gamma_l$ -- respectively the upper, bottom and lateral boundaries -- and equip this problem with the following free-slip rigid-lid  boundary conditions \cite{BS_2021,DGHT_2011},
\begin{align}
    \partial_z v  = 0, && w = 0, && \nu_T\partial_z T + \alpha_T T =0, && \partial_z S=0 && \text{on } \Gamma_u, \label{boundary-conditions}\\
    \partial_z v = 0, && w = 0, && \partial_z T =0, && \partial_z S=0 && \text{on } \Gamma_b, \nonumber \\
    \partial_{\bold{n}_H} v \times \bold{n}_H = 0, && v \bcdot \bold{n}_H = 0, && \partial_{\bold{n}_H} T =0, && \partial_{\bold{n}_H} S =0 && \text{on } \Gamma_l. \nonumber
\end{align}
Moreover, we consider an initial condition $U(t=0) = U_0 = (v_0, T_0, S_0)^\tr$ and an It\={o}-Stokes drift $U_s = (v_s, 0, 0)^\tr$ which fulfil \eqref{boundary-conditions}.

\section{Abstract framework and main results}

\subsection{Function spaces}\label{subsec-def-spaces}

\subsubsection{Sobolev spaces and noise regularity}

In this subsection, we define the Sobolev spaces used below. Remind that the spatial domain is denoted by $\mathcal{S} = \mathcal{S}_H \times [-h,0] \subset \R^3$. We define the following function spaces on $\Omega_S \in \{\mathcal{S}, \mathcal{S}_H\}$, for any real numbers $1 \leq q < \infty$ and $1 \leq p \leq \infty$, and any integers $d \geq 1$ and $m \geq 0$,
\begin{gather*}
    L^q(\Omega_S,\R^d) = \bigg\{u : \Omega_S \rightarrow \R^d \bigg| \int_{\Omega_S} \|u\|^q < +\infty\bigg\},\\
    L^\infty(\Omega_S,\R^d) = \bigg\{u : \Omega_S \rightarrow \R^d \bigg| \esssup_{\Omega_S} \|u\| < +\infty\bigg\},\\
    \modif{W^{m,p}(\Omega_S,\R^d) = \bigg\{u \in L^p(\Omega_S,\R^d) \bigg| D^\alpha u \in L^p(\Omega_S, \R^d), |\alpha|\leq m \bigg\}.}
\end{gather*}
We write $W^{m,2}(\Omega_S,\R^d) = H^m(\Omega_S,\R^d)$ for any non negative integer $m$. Additionally, we define the spaces $H^s(\Omega_S,\R^d)$ by interpolation, for any positive real number $s$. Furthermore, for any Banach space $\mathcal{B}$, and for $I$ an interval of $\R_+$, we denote by $C(I, \mathcal{B})$ the space of continuous functions from $I$ to $\mathcal{B}$, and the function spaces
\begin{gather*}
    L^p(I, \mathcal{B}) = \Big\{ f : I \rightarrow \mathcal{B} \: \Big| \: \int_I \|f\|_{\mathcal{B}}^p < \infty \Big\}, \\
    W^{s,p}(I, \mathcal{B}) = \Big\{ f \in L^p(I, \mathcal{B}) \: \Big| \: \int_I \int_I \frac{\|f(t) - f(t')\|^p_{\mathcal{B}}}{|t-t'|^{1+s p}} dt dt'  < \infty \Big\},
\end{gather*}
with
\begin{equation*}
    \|\bcdot\|_{L^p(I, \mathcal{B})} = \Big(\int_I \|\bcdot\|_{\mathcal{B}}^p dt\Big)^{1/p}\!\!\!\!\!\!, \|\bcdot\|_{W^{s,p}(I, \mathcal{B})} = \Big(\int_I \|\bcdot\|_{\mathcal{B}}^p dt + \int_I \int_I \frac{\|f(t) - f(t')\|^p_{\mathcal{B}}}{|t-t'|^{1+s p}} dt dt'\Big)^{1/p}\!\!\!\!\!\!,
\end{equation*}
(see \cite{FG_1995}). For any interval $I \subset \R_+$, we also define $W^{s,p}_{loc}(I,\mathcal{B}) := \Big\{ f : I \rightarrow \mathcal{B} \: \Big| \: \forall T >0, f \in W^{s,p}([0,T], \mathcal{B}) \Big\}$. In addition, for any Hilbert spaces $\mathcal{H}_1$ and $\mathcal{H}_2$, we define $\mathcal{L}_2( \mathcal{H}_1,\mathcal{H}_2)$ the space of Hilbert-Schmidt operators from $\mathcal{H}_1$ to $\mathcal{H}_2$, and $\|\bcdot\|_{\mathcal{L}_2( \mathcal{H}_1,\mathcal{H}_2)}$ its associated norm.

Furthermore, the noise is assumed to be regular enough in the following sense,
\begin{gather}
    \sup_{t\in [0,T]} \sum_{k=0}^\infty \|\phi_k\|_{H^4(\mathcal{S}, \R^3)}^2 < \infty, \quad u_s \in L^\infty\Big([0,T], H^4(\mathcal{S}, \R^3)\Big), \label{smoothness-noise}\\
    d_t u_s \in L^\infty\Big([0,T],H^3(\mathcal{S}, \R^3)\Big), \quad  a \nabla u_s \in L^\infty\Big([0,T],H^2(\mathcal{S},\R^3)\Big), \nonumber
\end{gather}
and 
\begin{equation*}
    a \in H^1\Big([0,T],H^3(\mathcal{S},\R^{3 \times 3})\Big).
\end{equation*}

Notice that condition \eqref{smoothness-noise} allows to give a meaning to the value of $\phi_k$ on the boundary. These regularity assumptions are not limiting in practice since most models consider spatially smooth noises for ocean models \cite{TML_2023}, as they are the physically observed ones.

In addition, we impose that the noise fulfils some boundary conditions. We consider two possible settings, for martingale and pathwise solutions respectively. To prove the existence of martingale solutions, the following non-penetration boundary condition is used,
\begin{equation}
    \phi_k \bcdot \bold{n} = 0 \text{ on $\partial \mathcal{S}$}. \label{eq-noiseBC-weak}
\end{equation}
In particular, the last assumption is made in the first point of Theorem \ref{theorem-filtered-weak}. Moreover, to prove the existence of a pathwise solution, we assume that the noise cancels on the lateral boundary as well -- that is, we assume the following \emph{stronger} assumption,
\begin{equation}
    \phi_k \bcdot \bold{n} = 0 \text{ on $\Gamma_u \cup \Gamma_b$} \label{eq-noiseBC-strong} \: \text{ and } \: \phi_k = 0 \text{ on $\Gamma_l$}.
\end{equation}
In the following, we often decompose $\phi_k$ as $\phi_k = (\phi_k^H \quad \phi_k^z)^\tr$, where $\phi_k^H$ is a 2D vector and $\phi_k^z$ a scalar.

\subsubsection{Rigid-lid boundary conditions spaces}

Let us define the function spaces associated to the rigid-lid boundary conditions \eqref{boundary-conditions}. We define first the following inner products
\begin{gather*}
    (v,v^\sharp)_{H_1} = (v,v^\sharp)_{L^2(\mathcal{S},\R^2)}, \quad (v,v^\sharp)_{V_1} = (\nabla v, \nabla v^\sharp)_{L^2(\mathcal{S},\R^2)}\\
    (T,T^\sharp)_{H_2} = (T,T^\sharp)_{L^2(\mathcal{S},\R)}, \quad (T,T^\sharp)_{V_2} = (\nabla T,\nabla T^\sharp)_{L^2} + \frac{\alpha_T}{\nu_T}(T,T^\sharp)_{L^2(\Gamma_u, \R)},\\
    (S,S^\sharp)_{H_3} = (S,S^\sharp)_{L^2(\mathcal{S},\R)}, \quad (S,S^\sharp)_{V_3} = (\nabla S,\nabla S^\sharp)_{L^2(\mathcal{S},\R)}.
\end{gather*}
Then, let
\begin{align}
    (U,U^\sharp)_H = (v,v^\sharp)_{H_1} + (T,T^\sharp)_{H_2} + (S,S^\sharp)_{H_3}, \nonumber\\
    (U,U^\sharp)_V = (v,v^\sharp)_{V_1} + (T,T^\sharp)_{V_2} + (S,S^\sharp)_{V_3},
\end{align}
for $U,U^\sharp \in L^2(\mathcal{S},\R^4)$ and $U,U^\sharp \in H^1(\mathcal{S},\R^4)$, respectively. Also we denote by $\|.\|_H, \|.\|_{H_i}$ and $\|.\|_V,\|.\|_{V_i}$ the associated norms. With a slight abuse of notation, we may write $\|.\|_H, \|.\|_V$ in place of $\|.\|_{H_i}, \|.\|_{V_i}$, respectively. Then, denote by $\mathcal{V}_1$ the space of functions of $C^\infty(\mathcal{S},\R^2)$, such that for all $v \in \mathcal{V}_1$, $\nabla_H \bcdot \int_{-h}^0 v dz = 0$ on $\mathcal{S}_H$ and $v \bcdot \bold{n} = 0$ on $\Gamma_l$. In addition, define $\mathcal{V}_2 = \mathcal{V}_3$ the space of functions of $C^\infty(\mathcal{S},\R)$ that average to zero over $\mathcal{S}$. Denote by $H_i$ the closure of $\mathcal{V}_i$ for the norm $\|.\|_{H_i}$, and $V_i$ its closure by $\|.\|_{V_i}$. Eventually, define $H=H_1 \times H_2 \times H_3$ and $V=V_1 \times V_2 \times V_3$, which are also the closures of $\mathcal{V}_1 \times \mathcal{V}_2 \times \mathcal{V}_3$ by $\|.\|_{H}$ and $\|.\|_{V}$, respectively. Often, by abuse of notation, we write $(\bcdot,\bcdot)_H$ instead of $(\bcdot,\bcdot)_{V' \times V}$. More generally, if $K$ is a subspace of $H$ and $K'$ its dual space, we write $(\bcdot,\bcdot)_H$ instead of $(\bcdot,\bcdot)_{K' \times K}$. Moreover, we define $\mathcal{D}(A) = V \cap H^2(\mathcal{S},\R^4)$, where the operator $A$ is defined below \eqref{def-operators}. As such, $A : \mathcal{D}(A) \rightarrow H$ is an unbounded operator.

Notice that, for consistency, we distinguished the spaces $H_2$ and $H_3$ even though they formally denote the same space. However, $V_2$ and $V_3$ \emph{are} different spaces since they are not equipped with the same inner products due to different boundary conditions on the temperature and salinity (Robin and Neumann respectively). This distinction allows to interpret $H_2$ and $V_2$ as temperature spaces, and $H_3$ and $V_3$ as salinity spaces. In addition, $H_1$ and $V_1$ are interpreted as horizontal velocity spaces ($\R^2$-valued processes). Using this formalism, the vertical velocity $w$ is written as a functional of the horizontal velocity $v$ through the continuity equation, namely $w(v) = \int_z^0 \nabla_H \cdot v \: dz'$.

Eventually, we define the barotropic and baroclinic projectors $\mathcal{A}_2 : \R^3 \rightarrow \R^2$, $\mathcal{A} : \R^3 \rightarrow \R^3$ and $\mathcal{R} : \R^3 \rightarrow \R^3$ of the velocity component as follows. For $v \in H$, $h$ being the depth of the ocean, let
\begin{equation}
    \mathcal{A}_2 [v](x,y) = \frac{1}{h} \int_{-h}^0 v(x,y,z') dz', \quad \mathcal{A} [v](x,y,z) = \mathcal{A}_2 [v](x,y), \quad \mathcal{R} [v] = v - \mathcal{A} [v].
\end{equation}
Remark that $\mathcal{A}$ and $\mathcal{R}$ are orthogonal projectors with respect to the inner product $(\bcdot, \bcdot )_{H}$. To simplify notations, we may use $\bar{v}$ in place of $\mathcal{A}_2 [v]$ or $\mathcal{A} [v]$, and $\tilde{v}$ in place of $\mathcal{R} [v]$.

\begin{remark}[Useful algebraic rules on barotropic and baroclinic modes]\hfill\\
    \noindent
    Let $f,g : \mathcal{S} \rightarrow \R$ two scalar functions. Then,
    $$\mathcal{A} (fg) = \bar{f} \bar{g} + \mathcal{A} (\tilde{f} \tilde{g}), \quad \mathcal{R} (fg) = \bar{f} \tilde{g} + \tilde{f} \bar{g} + \mathcal{R} (\tilde{f} \tilde{g}).$$
    In particular, if $\tilde{f}=0$ i.e. $f=\bar{f}$, then
    $$\mathcal{A} (fg) = \bar{f} \bar{g} = f \bar{g}, \quad \mathcal{R} (fg) = \bar{f} \tilde{g} = f \tilde{g},$$
    that is $\mathcal{A}$ and $\mathcal{R}$ commute with the functional $g \mapsto fg$ when $f$ is barotropic. Additionally, for all $i \in \{x,y\}$, we have
    $$\mathcal{A}\partial_i f = \partial_i \mathcal{A} f, \quad \mathcal{R} \partial_i f = \partial_i \mathcal{R} f.$$
    If in addition $f=0$ on $\Gamma_u \cup \Gamma_b$, this is also true for $i=z$.
\end{remark}

\subsection{Abstract formulation of the problems}

In this section, we aim to express, in abstract form, the problem under the weak (low-pass) filtered hydrostatic hypothesis \eqref{weak-filtered}. First, define the 4D vector $U$, representing the state of the system, and the correction $U^*$ of $U$ by the It\={o}-Stokes drift, as
\begin{equation}
    U = (v, T, S)^\tr, \quad U^* = (v^*, T, S)^\tr = (v - v_s, T, S)^\tr,
\end{equation}
and denote the advection operator by
\begin{equation}
    B(U^*,\bcdot) = B(v^*,\bcdot) := (v^* \bcdot \nabla_H)( \bcdot) + w(v^*) \partial_z (\bcdot),
\end{equation}
whenever it is defined. Here, we make an abuse of notation by identifying $B(U^*, \bcdot)$ to $B(v^*, \bcdot)$. This does not pose a problem, since the advection only depends on the velocity components of  $U^*$. Then, we define Leray type projectors $\mathbf{P}^v$ and $\mathbf{P}$ as follows (see \cite{BS_2021}),
\begin{equation}
    \mathbf{P}^v v = \mathbf{P}_{2D} \mathcal{A} v + \mathcal{R} v = \mathbf{P}_{2D} \bar{v} + \tilde{v}, \quad \mathbf{P} U = (\mathbf{P}^v v, T, S)^\tr,
\end{equation}
where $\mathbf{P}_{2D}$ is the 2D Leray projector, subject to the boundary condition $\bar{v} \bcdot \bold{n} =0$ on $\partial \mathcal{S}_H$, which is associated to the barotropic component $\bar{v}$. Notice that the baroclinic component $\tilde{v}$ is left unchanged by the projector $\mathbf{P}^v$, that is $\mathbf{P}$ only affects the barotropic component of the velocity. In addition, $\mathbf{P}$ is the identity over the temperature and salinity components. Also, for notational convenience we keep the same notations for the composition of the following operators with the Leray projector,
\begin{gather}
    AU = \mathbf{P} ( \Av v, \: \AT T, \: \AS S)^\tr, \quad B(U^*,U) = \mathbf{P} ( B(v^*,v), \: B(v^*,T), \: B(v^*,S))^\tr, \nonumber\\
    \Gamma U = \mathbf{P} ( \Gamma v, \: 0, \: 0)^\tr. \label{def-operators}
\end{gather}
With this definition, $A: \mathcal{D}(A) \rightarrow H$ is an unbounded operator with $\mathcal{D}(A) = V \cap H^2(\mathcal{S},\R^4)$. By abuse of notation again, we will often denote $Av$ for $\mathbf{P}^v A^v v$, which is the velocity component of $A U$. This extends to the temperature and salinity, since we write $A T$ and $A S$ in place of $A^T T$ and $A^S S$ respectively, them being the temperature and salinity components of $A U$.

\bigskip

\noindent
At this point we may state some useful facts: let $v^\sharp \in \mathcal{D}(A)$, so that
\begin{equation}
    \nabla \bcdot \bar{v}^\sharp =0, \quad \bar{v}^\sharp \bcdot \bold{n} = \frac{\partial\bar{v}^\sharp}{\partial\bold{n}} \times \bold{n}=0 \text{ on } \partial \mathcal{S}_H.
\end{equation}
In particular,
\begin{enumerate}
    \item Using the boundary conditions \eqref{boundary-conditions} on $\Gamma_u \cup \Gamma_b$, we remark that the anisotropic Laplace operator $A^v$ fulfils
    $$\mathcal{A} A^v v^\sharp = \frac{1}{h} \int_{-h}^0 \big[\mu_v (-\Delta_H) v^\sharp + \nu_v (-\partial_{zz}) v^\sharp \big] dz' = \mu_v (-\Delta_H) \bar{v}^\sharp = A^v \bar{v}^\sharp.$$
    Thus, $\mathcal{R} A^v v = A^v \tilde{v}^\sharp$.
    \item Following the remark of \cite{CT_2007}, since $v^\sharp \in \mathcal{D}(A)$, then $ \Delta_H \bar{v}^\sharp \in L^2(\mathcal{S}_H)$. Therefore, $\Delta_H \bar{v}^\sharp \bcdot \bold{n} \in H^{-1/2}(\partial \mathcal{S}_H)$ -- see \cite{Temam_2001} for instance. Also, using the boundary conditions on $\bar{v}^\sharp$, we deduce that $\Delta_H \bar{v}^\sharp \bcdot \bold{n} =0$ on $\partial \mathcal{S}_H$ \cite{Ziane_1998}. Hence, $A^v v^\sharp = \mathbf{P}^v A^v v^\sharp$. Consequently, the operator $A$ is equal to $(A^v, A^T, A^S)$.
\end{enumerate}
Combining the previous remarks, we infer as a result,
\begin{equation}
    \mathcal{A} A v^\sharp = A \bar{v}^\sharp = A^v \bar{v}^\sharp, \text{  and  } \mathcal{R} A v^\sharp = A \tilde{v}^\sharp = A^v \tilde{v}^\sharp.
\end{equation}
In the following, we will assume that $U_s = (v_s, 0, 0)^\tr \in V$.

\subsubsection{Low-pass filtered problem $(\mathcal{P}_K)$}

These remarks being made, we aim to derive an abstract formulation for the (low-pass) filtered problem. From subsection \ref{sec-filtered-problem}, we infer
\begin{multline}
    d_t U + \Bigg[AU + B(U^*,U) + \Gamma U + \frac{1}{\rho_0} \mathbf{P} \begin{pmatrix}
        \nabla_H p \\0 \\0
    \end{pmatrix} - \frac{1}{2} \mathbf{P} \nabla \bcdot(a\nabla U)\Bigg] dt = - \mathbf{P}( \sigma dW_t \bcdot \nabla U)\\
    - (A + \Gamma) (\sigma^H dW_t) - \frac{1}{\rho_0} \mathbf{P} \begin{pmatrix}
        \nabla_H dp_t^\sigma \\0 \\0
    \end{pmatrix}.
\end{multline}
Reminding the notations $U_s = (v_s, 0, 0)^\tr$ and $(v_s,w_s)^\tr = u_s = \frac{1}{2}(\nabla \bcdot a)$, we rewrite the previous equation in terms of $U^*=U - U_s$ with a change of variable, to get
\begin{multline}
    d_t U^* + [AU^* + B(U^*) + \Gamma U^* + \frac{1}{\rho_0} \mathbf{P} \begin{pmatrix}
        \nabla_H p \\0 \\0
    \end{pmatrix} + F_\sigma(U^*)]dt\\
    = G_\sigma(U^*)dW_t - \frac{1}{\rho_0} \mathbf{P} \begin{pmatrix}
        \nabla_H dp_t^\sigma \\0 \\0
    \end{pmatrix},
\end{multline}
where the operators $F_\sigma$ and $G_\sigma$ are defined as
\begin{align}
    F_\sigma(U^*) dt &= \mathbf{P} \Big[ d_t U_s + [ B(U^*,U_s) - \frac{1}{2} \nabla \bcdot (a \nabla U_s) + AU_s + \Gamma U_s - \frac{1}{2} \nabla \bcdot (a \nabla U^*) ]dt \Big],\\
    G_\sigma(U^*) dW_t &= \mathbf{P} \Big[ - (\sigma dW_t \bcdot \nabla) U^* - (\sigma dW_t \bcdot \nabla) U_s - A (\sigma^H dW_t) - \Gamma (\sigma^H dW_t) \Big].
\end{align}
Moreover, as $u_s$ is divergence-free, $w_s = w(v_s)$ follows from the definition of operator $w(v)$ \eqref{w-definition}. Then, we derive the following relations for the pressure terms, using equations \eqref{primitive-p} and \eqref{primitive-dp},
\begin{align}
    \frac{1}{\rho_0} \nabla_H p = &-g \:  \nabla_H \int_z^0 {(\beta_T T + \beta_S S) dz'} + \frac{1}{\rho_0}\nabla_H p^s \nonumber\\
    &+ \nabla_H \int_z^0 \Big[u_s \bcdot \nabla (w(v^*) + w_s) - \frac{1}{2} \nabla \bcdot (a^K \nabla (w(v^*) + w_s))\Big] dz',\\
    \frac{1}{\rho_0} \nabla_H (dp_t^\sigma) = &- \nabla_H \int_z^0 \Big[K* [\sigma dW_t \bcdot \nabla w(v^*)] + K*[\sigma dW_t \bcdot \nabla w_s] + A (\sigma^z dW_t) \Big]dz' \nonumber\\
    & \qquad \qquad \qquad \qquad \qquad \qquad \qquad \qquad \qquad \qquad + \frac{1}{\rho_0} \nabla_H dp_t^{\sigma,s}.
\end{align}
The quantities $p^s \: dt$ and $dp_t^{\sigma,s}$ are respectively the bounded variation and the martingale contributions to the surface pressure. As they are independent on the z-axis, we have for all $v^\sharp \in V_1$,
$$\Big(\frac{1}{\rho_0} \mathbf{P}^v \nabla_H p^s, v^\sharp \Big)_H = -\Big(\frac{1}{\rho_0} p^s, \nabla_H \bcdot \mathcal{A} v^\sharp \Big)_H =0,$$
using the boundary conditions on $\mathcal{A}v^\sharp$: $\mathcal{A}v^\sharp =0$ on $\mathcal{S}_H$ and $\mathcal{A}v^\sharp\bcdot\bold{n} = \frac{\partial}{\partial\bold{n}}\mathcal{A}v^\sharp \times \bold{n}=0$ on $\partial \mathcal{S}_H$. This implies $\mathbf{P}^v \nabla_H p^s = 0$. Similarly, we get $\mathbf{P}^v \nabla_H dp_t^{\sigma,s} = 0$. Therefore, we obtain the following relations,
\begin{align*}
    \frac{1}{\rho_0} \mathbf{P}^v [\nabla_H p] &=  \mathbf{P}^v \Bigg[ -g \:  \nabla_H \int_z^0 {(\beta_T T + \beta_S S) dz'} \Bigg] \\
    &- \mathbf{P}^v \Bigg[\nabla_H \int_z^0 \Big[K*[u_s \bcdot \nabla (w(v^*) + w_s)] + \frac{1}{2} \nabla \bcdot (a^K \nabla (w(v^*) + w_s))\Big] dz' \Bigg],
\end{align*}
and,
\begin{equation*}
    \frac{1}{\rho_0} \mathbf{P}^v[\nabla_H (dp_t^\sigma) ] = \mathbf{P}^v \Bigg[ \nabla_H \int_z^0 \Big[K*[\sigma dW_t \bcdot \nabla (w(v^*)+w_s) + A (\sigma^z dW_t) \Big]dz'\Bigg].
\end{equation*}
On the one hand, the bounded variation surface pressure $p^s$ can be interpreted as a Lagrange multiplyer associated to the constraint $\nabla_H \bcdot \mathcal{A} v=0$. This interpretation aligns with the proof proposed by Cao and Titi in the deterministic case \cite{CT_2007}, which highlights  that the barotropic mode follows a 2D Navier-Stokes equation, while the baroclinic mode evolves according to a 3D Burgers equation, up to some coupling terms. Remarkably, the bounded variation surface pressure does not affect the baroclinic dynamics, the barotropic mode being divergence-free using the vertical boundary conditions. On the other hand, the martingale surface pressure $p^{s,\sigma}$ arises from the stochastic modelling we proposed, and may be seen as a perturbation of the pressure $p^s$. This new term $p^{s,\sigma}$ affects \emph{a priori} both the barotropic and the baroclinic dynamics. However, in the strong hydrostatic hypothesis the martingale pressure equation simplifies to $\partial_z dp_t^\sigma = 0$, that is $dp_t^\sigma = dp_t^{s,\sigma}$ (equation \eqref{strong-hydro}). Consequently, the martingale pressure term becomes completely barotropic, affecting only the barotropic dynamics. This reduction allows the use of methods similar to those employed in the deterministic case.

\noindent
We are now in position to define $(\mathcal{P}_K)$, for $K \in H^3(\mathcal{S}, \R)$, the abstract primitive equations problem with weak (low-pass) filtered hydrostatic hypothesis, as follows,
\begin{equation} (\mathcal{P}_K) \left\{\begin{array}{l}
     d_t U^* + [AU^* + B(U^*) + \Gamma U^* + \frac{1}{\rho_0} \mathbf{P} \nabla_H p  + F_\sigma(U^*)]dt = G_\sigma(U^*)dW_t - \frac{1}{\rho_0} \mathbf{P} \nabla_H dp_t^\sigma,\\
     \begin{array}{rcl}
        \frac{1}{\rho_0} \mathbf{P} \nabla_H (dp_t^\sigma) & = & \mathbf{P} \Bigg[ \nabla_H \int_z^0 \big[K * [\sigma dW_t \bcdot \nabla (w(v^*) + w_s)] + A (\sigma^z dW_t) \big]dz' \Bigg],\\
        \frac{1}{\rho_0} \mathbf{P} \nabla_H p & = & \mathbf{P} \Bigg[ - \nabla_H \int_z^0 K*\big[u_s \bcdot \nabla (w(v^*) + w_s) \big] dz'\\
        & & - \nabla_H \int_z^0 \frac{1}{2} \nabla \bcdot (a^K \nabla (w(v^*) + w_s)) dz' - g \:  \nabla_H \int_z^0 {(\beta_T T + \beta_S S) dz'} \Bigg]. 
    \end{array}
\end{array}\right. \label{problem-PK}
\end{equation}
This problem is equipped with an initial condition to be specified later, and we make the noise regularity assumptions \eqref{smoothness-noise}. For notational convenience, we have written $\nabla_H p$ and $\nabla_H dp_t^{\sigma}$ in place of $(\nabla_H p, \: 0, \: 0)^\tr$ and $(\nabla_H dp_t^{\sigma}, \: 0, \: 0)^\tr$ respectively. In addition, we denote by $(\pi_k)_k$ and $(\pi_k^s)_k$ the bases of $L^2(\mathcal{S},\R)$ and $L^2(\mathcal{S}_H,\R)$ such that
\begin{align*}
    dp_t^{\sigma} = \sum_{k=0}^\infty \pi_k d\beta_t^k , \text{ and } dp_t^{\sigma,s} = \sum_{k=0}^\infty \pi_k^s d\beta_t^k,
\end{align*}
and we apply the same abuse of notation to the bases $(\phi_k)_k$, $(\pi_k)_k$ and $(\pi_k^s)_k$.

\subsubsection{Approximated low-pass filtered problem $(\hat{\mathcal{P}}_K)$} \label{approx-low-pass-filtered-problem}

We also propose an approximation of the problem $(\mathcal{P}_K)$ in the limit of ``quasi-barotropic'' flow, that we denote by $(\hat{\mathcal{P}}_K)$. Remind first that
\begin{equation}
    \sigma dW_t = \sum_k \phi_k d\beta^k,
\end{equation}
where $(\phi_k)$ is a family of functions in $L^2(\mathcal{S}, \R^3)$. In this subsection, we assume that each horizontal component $\phi_k^H$ is either barotropic or baroclinic -- i.e. independent of $z$ or averaging to $0$ along the $z$-axis, respectively. This is equivalent to considering that $\sigma^H dW_t$ is a sum of two independent noises, $\sigma_{BT}^H dW_t$ and $\sigma_{BC}^H dW_t$, the former one being barotropic and the latter being baroclinic. Due to the divergence-free and boundary conditions on the noise, we denote by
\begin{align}
    \phi_k^{BT} &:= \phi_k = \begin{pmatrix}
        \phi_k^H \\0
    \end{pmatrix} \text{ when $\phi_k^H$ is barotropic}, 0 \text{ otherwise,}\\
    \phi_k^{BC} &:= \phi_k = \begin{pmatrix}
        \phi_k^H \\ w(\phi_k^H)
    \end{pmatrix} \text{ when $\phi_k^H$ is baroclinic}, 0 \text{ otherwise.}
\end{align}
Denoting by $\beta^{k,BT}$ and $\beta^{k,BC}$ their respectively associated Brownian motions, we make
\begin{equation}
    \sigma_{BT} dW_t = \sum_k \phi_k^{BT} d\beta^{k,BT}, \quad \sigma_{BC} dW_t = \sum_k \phi_k^{BC} d\beta^{k,BC}.
\end{equation}
As a consequence, since $\beta^{k,BT}$ and $\beta^{p,BC}$ are independent for all $k$ and $p$, the variance tensor $a$ can be split into two terms $a_{BT}$ and $a_{BC}$, so that
\begin{equation}
    a = a_{BT} + a_{BC}, \quad a_{BT} = \sum_k \phi_k^{BT} (\phi_k^{BT})^\tr, \quad a_{BC} = \sum_k \phi_k^{BC} (\phi_k^{BC})^\tr.
\end{equation}
Eventually, we define the problem $(\hat{\mathcal{P}}_K)$ as follows,
\begin{equation} (\hat{\mathcal{P}}_K) \left\{\begin{array}{l}
     d_t U^* + [AU^* + B(U^*) + \Gamma U^* + \frac{1}{\rho_0} \mathbf{P} \nabla_H p  + \hat{F}_\sigma(U^*)]dt = \hat{G}_\sigma(U^*)dW_t - \frac{1}{\rho_0} \mathbf{P} \nabla_H dp_t^\sigma,\\
     \begin{array}{rcl}
        \frac{1}{\rho_0} \mathbf{P} \nabla_H (dp_t^\sigma) & = & \mathbf{P} \Bigg[ \nabla_H \int_z^0 \Big[K * [\sigma dW_t \bcdot \nabla (w(v^*) + w_s)] + A (\sigma^z dW_t) \Big]dz' \Bigg],\\
        \frac{1}{\rho_0} \mathbf{P} \nabla_H p & = & \mathbf{P} \Bigg[ - \nabla_H \int_z^0 K*\Big[u_s \bcdot \nabla (w(v^*) + w_s) \Big] dz'\\
        & & - \nabla_H \int_z^0 \frac{1}{2} \nabla \bcdot (a^K \nabla (w(v^*) + w_s)) dz' - g \:  \nabla_H \int_z^0 {(\beta_T T + \beta_S S) dz'} \Bigg], 
    \end{array}
\end{array}\right. \label{problem-approx-PK}
\end{equation}
where the operators $\hat{F}_\sigma$ and $\hat{G}_\sigma$ are defined by,
\begin{multline}
    \hat{F}_\sigma(U^*) dt = \mathbf{P} \Bigg[ d_t U_s + [ B(U^*,U_s) + AU_s + \Gamma U_s ]dt\\
    - \frac{1}{2} \begin{pmatrix}
        \nabla \bcdot (a_{BT} \nabla v^*) + \nabla \bcdot (\overline{a_{BC}} \nabla \overline{v}^*)\\
        \nabla \bcdot (a \nabla T)\\
        \nabla \bcdot (a \nabla S)
    \end{pmatrix} dt - \begin{pmatrix}
        \nabla \bcdot (a_{BT} \nabla v_s) + \nabla \bcdot (\overline{a_{BC}} \nabla \overline{v}_s)\\
        0\\
        0
    \end{pmatrix} dt \Bigg],
\end{multline}
and
\begin{multline}
    \hat{G}_\sigma(U^*) dW_t = \mathbf{P} \Bigg[ - A (\sigma^H dW_t) - \Gamma (\sigma^H dW_t) \\
    - \begin{pmatrix}
        (\sigma_{BT} dW_t \bcdot \nabla) v^* + (\sigma_{BC} dW_t \bcdot \nabla) \overline{v}^*\\
        (\sigma dW_t \bcdot \nabla) T\\
        (\sigma dW_t \bcdot \nabla) S
    \end{pmatrix} - \begin{pmatrix}
        (\sigma_{BT} dW_t \bcdot \nabla) v_s + (\sigma_{BC} dW_t \bcdot \nabla) \overline{v}_s\\
        0\\
        0
    \end{pmatrix}\Bigg].
\end{multline}
The problem is supplemented with an initial condition to be specified below, and we make the same noise regularity assumptions \eqref{smoothness-noise}.

\begin{remark}\hfill \label{big-remark}\\
\begin{enumerate}
    \item Let us precise in what sense $(\hat{\mathcal{P}}_K)$ is an approximation of $(\mathcal{P}_K)$. Regarding the horizontal momentum equation, we have that
    \begin{equation}
        \mathbbm{D}_t v + \Gamma (v \: dt + \sigma^H dW_t) =  -\Av (v\:dt  + \sigma^H dW_t) - \frac{1}{\rho_0} \nabla_H (p \: dt + dp_t^\sigma), \label{horizontal-mom-equation-for-approx}
    \end{equation}
    where the stochastic transport term $\mathbbm{D}_t v$ reads 
    \begin{equation}
        \mathbbm{D}_t v = d_t v + (u^* \bcdot \nabla_3) v + (\sigma dW_t \bcdot \nabla_3) v - \frac{1}{2} \nabla_3 \bcdot( a \nabla_3 v) dt. \label{transport-operator-approx}
    \end{equation}
    Denote by $r$ the typical ratio between the baroclinic and barotropic modes -- i.e. $|\tilde{v}| \sim r |\bar{v}|$ -- and assume that the unresolved scale enjoys the same ratio -- i.e. $\sigma_{BC} dW_t \sim r \sigma_{BT} dW_t$. It is physically expected that $r \ll 1$, since $\tilde{v}$ is much smaller than $\bar{v}$. On the one hand, expanding the third term on the RHS of \eqref{transport-operator-approx} yields,
    \begin{equation*}
        (\sigma dW_t \bcdot \nabla_3) v =  (\sigma_{BT} dW_t \bcdot \nabla_3) v + (\sigma_{BC} dW_t \bcdot \nabla_3) \overline{v} + (\sigma_{BC} dW_t \bcdot \nabla_3) \tilde{v}.
    \end{equation*}
    Hence, the physical energy associated to this martingale transport term reads
    \begin{multline}
        \sum_k \|(\phi_k^{BT} \bcdot \nabla_3) v\|_{L^2}^2 + \|(\phi_k^{BC} \bcdot \nabla_3) \overline{v} + (\phi_k^{BC} \bcdot \nabla_3) \tilde{v}\|_{L^2}^2 \\ = \sum_k \|(\phi_k^{BT} \bcdot \nabla_3) v\|_{L^2}^2 + \|(\phi_k^{BC} \bcdot \nabla_3) \overline{v}\|_{L^2}^2 + O(r^3). \label{martingale-energy-splitting}
    \end{multline} 
    On the other hand, splitting the fourth term on the RHS of \eqref{transport-operator-approx} yields,
    \begin{align*}
        \overline{\nabla_3 \bcdot( a \nabla_3 v)} &= \nabla_3 \bcdot( a_{BT} \nabla_3 \bar{v}) + \nabla_3 \bcdot( \overline{a}_{BC} \nabla_3 \bar{v}) + \overline{\nabla_3 \bcdot( \tilde{a}_{BC} \nabla_3 \tilde{v})}\\
        \reallywidetilde{\nabla_3 \bcdot( a \nabla_3 v)} &= \nabla_3 \bcdot( a_{BT} \nabla_3 \tilde{v}) + \nabla_3 \bcdot( \overline{a}_{BC} \nabla_3 \tilde{v}) + \nabla_3 \bcdot( \tilde{a}_{BC} \nabla_3 \bar{v}) +\reallywidetilde{\nabla_3 \bcdot( \tilde{a}_{BC} \nabla_3 \tilde{v})}.
    \end{align*}
    Therefore, by taking the inner product with $v$, the energy dissipated by the stochastic diffusion is expressed as
    \begin{align*}
        (v, &\nabla_3 \bcdot( a \nabla_3 v))_{L^2} = (\bar{v},\nabla_3 \bcdot( a_{BT} \nabla_3 \bar{v}) + \nabla_3 \bcdot( \overline{a}_{BC} \nabla_3 \bar{v}) + \overline{\nabla_3 \bcdot( \tilde{a}_{BC} \nabla_3 \tilde{v})})_{L^2}\\
        &+ (\tilde{v},\nabla_3 \bcdot( a_{BT} \nabla_3 \tilde{v}) + \nabla_3 \bcdot( \overline{a}_{BC} \nabla_3 \tilde{v}) + \nabla_3 \bcdot( \tilde{a}_{BC} \nabla_3 \bar{v}) +\reallywidetilde{\nabla_3 \bcdot( \tilde{a}_{BC} \nabla_3 \tilde{v})} )_{L^2}\\
        &= (\bar{v},\nabla_3 \bcdot( a_{BT} \nabla_3 \bar{v}))_{L^2} + (\bar{v},\nabla_3 \bcdot( \overline{a}_{BC} \nabla_3 \bar{v}))_{L^2} + (\tilde{v},\nabla_3 \bcdot( a_{BT} \nabla_3 \tilde{v}))_{L^2} + O(r^3).
    \end{align*}
    Retaining only the terms which contribute at order $r^2$ or more to the energy finally yields
    \begin{multline}
        \mathbbm{D}_t v \approx d_t v + (u^* \bcdot \nabla_3) v + (\sigma_{BT} dW_t \bcdot \nabla_3) v + (\sigma_{BC} dW_t \bcdot \nabla_3) \overline{v}\\
        - \frac{1}{2} \nabla_3 \bcdot( a_{BT} \nabla_3 v) dt - \frac{1}{2} \nabla_3 \bcdot( a_{BC} \nabla_3 \bar{v}) dt =: \mathbbm{D}_t^{approx} v. 
    \end{multline}
    Replacing $\mathbbm{D}_t$ by its approximation $\mathbbm{D}_t^{approx}$ in equation \eqref{horizontal-mom-equation-for-approx} and reasoning similarly as for deriving $(\mathcal{P}_K)$, we find the problem $(\hat{\mathcal{P}}_K)$.

    \item Notice that, when $\sigma_{BC} dW_t = 0$ (or equivalently $a_{BC} = 0$), $(\mathcal{P}_K)$ and $(\hat{\mathcal{P}}_K)$ are in fact equivalent. This stands as a relaxation of the assumption proposed in \cite{AHHS_2022, AHHS_2022_preprint}. There, the authors prove the global-in-time well-posedness of an interpretation of the primitive equations with a noise such that its horizontal component is independent of $z$ -- that is to say barotropic. Conversely, our approach yields an approximation of the horizontal momentum equation in the limit of small -- yet non-zero -- baroclinic noise components, provided they are decorrelated with the barotropic ones.
    
    \item Moreover, we draw attention to the fact that, if we assume the noise is divergence-free and has a barotropic horizontal component that fulfils the boundary condition \eqref{boundary-conditions}, then it is bidimensional in the following sense,
    $$\sigma dW_t = \begin{pmatrix}
        \sigma dW_t^H (x,y) \\
        0
    \end{pmatrix}, \quad \nabla_H \bcdot \sigma dW_t^H (x,y) = 0.$$
    Therefore, only compressible (divergent) tridimensional noises can fulfil the barotropic horizontal noise condition, and our results still hold true in this case -- see Remarks \ref{remark-energ-noise-dep-v} and \ref{remark-sup-noise} in subsection \ref{subsec-energy-1}. However, for such noises the expression of the It\={o}-Stokes drift becomes $u_s = \frac{1}{2} \nabla \bcdot a - \sum_k \phi_k (\nabla \bcdot \phi_k)$ -- see \cite{TCM_2021}. This drift velocity still fulfils the regularity conditions \eqref{smoothness-noise} for a regular enough noise $\sigma dW_t$.

    \item The aforementioned barotropic horizontal noise assumption can be related to the validity range of the primitive equations. The (deterministic) primitive equations are physically valid when the squared aspect ratio $\epsilon^2 := (h/L)^2$ is negligible compared to the Richardson number $Ri := \frac{N^2}{(\partial_z v)^2}$ \cite{MHPA_1997}. Here $v$ denotes the horizontal velocity, $h$ refers to the depth of the ocean, $L$ to the horizontal scale (e.g $\sqrt{|\mathcal{S}_H|}$), and $N^2= - \frac{g}{\rho_0} \partial_z \rho$. This condition reads
    \begin{align}
        \frac{\epsilon^2}{Ri} \ll 1, \quad \text{  or equivalently  } \quad (\partial_z v)^2 \ll \frac{N^2}{\epsilon^2}. \label{condition-Marshall}
    \end{align}
    In particular, the latter holds in the limit of small enough vertical shear of the horizontal component. In such a case, the horizontal component of the velocity is almost independent of z -- therefore it is, so to say, ``quasi-barotropic''. In the context of stochastic flows, the horizontal noise models a small scale velocity denoted by $\eta^{1/2} v'$, where $\eta^{1/2}$ is a scaling factor so that $v$ and $v'$ have the same order of magnitude. Thus, the condition \eqref{condition-Marshall} becomes
    \begin{align}
        (\partial_z (v + \eta^{1/2} v'))^2 \ll \frac{N^2}{\epsilon^2}, \quad \text{  which also reads  } \quad (\partial_z v)^2 , \: \eta (\partial_z v')^2 \ll \frac{N^2}{\epsilon^2}.
    \end{align}
    Therefore, the LU stochastic primitive equations are physically valid under the condition \eqref{condition-Marshall}, and when the horizontal noise modelling $\eta^{1/2} v'$ is either small enough ($\eta \rightarrow 0$) or quasi-barotropic ($(\partial_z v')^2 \rightarrow 0$). In this setting, a purely barotropic noise is equivalent to $\partial_z v' = 0$. The approximated problem $(\hat{\mathcal{P}}_K)$ corresponds to a setting where $\partial_z v$ and $\partial_z v'$ are both small, yet non-zero.

    \item Furthermore, assuming $\sigma_{BT} dW_t$ and $\sigma_{BC} dW_t$ are independent is crucial for the energy splitting argument to be valid in equation \eqref{martingale-energy-splitting}. Otherwise, more terms of order $r^2$ would emerge, for example $((\phi_k^{BT} \bcdot \nabla_3) v, (\phi_k^{BC} \bcdot \nabla_3) \tilde{v})_{L^2}$, the term $\sigma_{BC} dW_t \bcdot \nabla_3 \tilde{v}$ becoming non-negligible in the expression of $\mathbbm{D}_t^{approx} v$. Such assumption removes a direct \emph{martingale} dependence on $\nabla_3 \tilde{v}$ in the barotropic equation, which allows to use similar estimates as in \cite{AHHS_2022, AHHS_2022_preprint}. Not doing so would imply keeping the term $\sigma_{BC} dW_t \bcdot \nabla_3 \tilde{v}$ in the horizontal momentum equation, and estimate \emph{a priori} $\tilde{v}$ in $H^2$, while classically it is estimated in a (much larger) $L^p$-space.
\end{enumerate}
\end{remark}

\subsubsection{Eddy-viscosity problem $(\mathcal{P}^{\gamma_r}_{EV})$}

For $\gamma_r \geq 1$ and $\alpha > 0$, we define $(\mathcal{P}^{\gamma_r}_{EV})$, the abstract primitive equations problem with weak eddy-viscosity hydrostatic hypothesis, as follows,
\begin{equation}(\mathcal{P}^{\gamma_r}_{EV}) \left\{\begin{array}{l}
      d_t U^* + [AU^* + B(U^*) + \Gamma U^* + \frac{1}{\rho_0} \mathbf{P} \nabla_H p + F_\sigma(U^*)]dt = G_\sigma(U^*)dW_t - \frac{1}{\rho_0} \mathbf{P} \nabla_H dp_t^\sigma,\\
     \begin{array}{rcl}
        \frac{1}{\rho_0} \mathbf{P} \nabla_H (dp_t^\sigma) & = & \mathbf{P} \Bigg[ \nabla_H \int_z^0 \Big[\sigma dW_t \bcdot \nabla (w(v^*) + w_s) + A (\sigma^z dW_t) \Big]dz' \Bigg],\\
        \frac{1}{\rho_0} \mathbf{P} \nabla_H p & = & \mathbf{P} \Bigg[-g \:  \nabla_H \int_z^0 (\beta_T T + \beta_S S) dz' \Bigg] + \mathbf{P} \Bigg[ - \nabla_H \int_z^0 u_s \bcdot \nabla (w(v^*)+w_s) dz' \Bigg] \\
        & + & \mathbf{P} \Bigg[\nabla_H \int_z^0 \Big[- \frac{1}{2} \nabla \bcdot (a \nabla (w(v^*) + w_s)) - \alpha (-\Delta)^{\gamma_r} w \Big] dz' \Bigg]. 
    \end{array}
\end{array}\right. \label{problem-PEV}
\end{equation}
Again, the problem is supplemented with an initial condition to be specified later, and we make the same noise regularity assumptions \eqref{smoothness-noise}.

\subsubsection{(Eddy-viscosity) energy-balanced problem $(\mathcal{P}^{\gamma_r}_{EB})$}

In a similar fashion, for $\gamma_r \geq 1$ and $\alpha>0$, we define $(\mathcal{P}^{\gamma_r}_{EB})$, the abstract primitive equations problem with weak eddy-viscosity energy-balanced hydrostatic hypothesis, as follows,
\begin{equation}(\mathcal{P}^{\gamma_r}_{EB}) \left\{\begin{array}{l}
      d_t U^* + [AU^* + B(U^*) + \Gamma U^* + \frac{1}{\rho_0} \mathbf{P} \nabla_H p  + F_\sigma(U^*)]dt \\ \qquad \qquad \qquad \qquad \qquad \qquad = G_\sigma(U^*)dW_t - \frac{1}{\rho_0} \mathbf{P} \nabla_H dp_t^\sigma + \frac{1}{2} \mathbf{P} \nabla \bcdot \check{a}[\nabla w] dt + \frac{1}{2} \mathbf{P} \mathcal{C}_\sigma dt,\\
     \begin{array}{rcl}
        \frac{1}{\rho_0} \mathbf{P} \nabla_H (dp_t^\sigma) & = & \mathbf{P} \Bigg[ \nabla_H \int_z^0 \Big[\sigma dW_t \bcdot \nabla w(v^*) + \sigma dW_t \bcdot \nabla w_s + A (\sigma^z dW_t) \Big]dz' \Bigg],\\
        \frac{1}{\rho_0} \mathbf{P} \nabla_H p & = & \mathbf{P} \Bigg[-g \:  \nabla_H \int_z^0 {(\beta_T T + \beta_S S) dz'}\Bigg] + \mathbf{P} \Bigg[- \nabla_H \int_z^0 \Big[u_s \bcdot \nabla (w(v^*) + w_s)\Big]\Bigg]\\
        & + & \mathbf{P} \Bigg[ \nabla_H \int_z^0  - \frac{1}{2}\Big(\nabla \bcdot \hat{a}[\nabla v] + \nabla \bcdot \hat{\hat{a}} [\nabla w] \Big) dz' - \nabla_H \int_z^0 \alpha (-\Delta)^{\gamma_r} w(v^*) dz' \Bigg]\\
        & + & \frac{1}{2} \mathbf{P} \Bigg[\nabla_H \int_z^0 \sum_k \Big( \phi_k \bcdot \nabla \partial_z \pi_k + (\phi_k \bcdot \nabla) \Delta \phi_k^z\Big) dz'\Bigg],
    \end{array}
\end{array}\right. \label{problem-PEB}
\end{equation}
\modif{where we use the notations introduced in equations \eqref{eq:def-Ito-Strato-term} to \eqref{equation-C-sigma}.} Once again, the initial condition is specified below, and we make the same noise regularity assumptions \eqref{smoothness-noise}.

\subsection{Main results}

We are now in position to state our main results. We remind that, in the definition of martingale solutions, the stochastic basis -- that is the filtered probability space and the Wiener process -- is unknown \emph{a priori}. Our main result concerns the well-posedness of the filtered problem $(\mathcal{P}_K)$ and its approximation $(\hat{\mathcal{P}}_K)$,
\begin{theorem}\label{theorem-filtered-weak}
Suppose $K \in H^{3}(\mathcal{S}, \R)$. Then, the following propositions hold,
\begin{enumerate}
    \item Equip the problem $(\mathcal{P}_K)$ -- resp. $(\hat{\mathcal{P}}_K)$ -- with the initial condition $U_0 = (v_0, T_0, S_0)^\tr \in H$, and assume that the noise fulfils \eqref{smoothness-noise} and \eqref{eq-noiseBC-weak}. Then, $(\mathcal{P}_K)$ -- resp. $(\hat{\mathcal{P}}_K)$ -- admits at least one global-in-time martingale solution $(\mathcal{S}_0, U)$, the stochastic basis $\mathcal{S}_0$ not being fixed \emph{a priori}. In addition, for all $T>0$,
    $$U \in L^2\Big(\Omega_0,L^2\big([0,T],V\big)\Big) \cap L^2\Big(\Omega_0,L^\infty\big([0,T],H\big)\Big),$$
    where $\Omega_0$ is the sample space of associated to $\mathcal{S}_0$.
    \item Assume that $U_0 \in V$ and fix a stochastic basis \emph{a priori}. In addition, assume that the noise fulfils \eqref{smoothness-noise} and \eqref{eq-noiseBC-strong}. Then, there exists a stopping time $\tau>0$ such that $(\mathcal{P}_K)$ -- resp. $(\hat{\mathcal{P}}_K)$ -- admits a local-in-time pathwise solution $U^*$, which fulfils, for all $T>0$ and for all stopping time $0< \tau' < \tau$,
    $$U_{\tau' \wedge \bcdot}^* \in L^2\Big(\Omega, L^2\big([0,T],\mathcal{D}(A)\big)\Big) \cap L^2\Big(\Omega, C([0,T],V)\Big),$$
    where $\Omega$ denotes the sample space associated to the aforementioned stochastic basis. This solution is unique up to indistinguishability, that is for all solutions $U^*$ and $\hat{U}^*$ of $(\mathcal{P}_K)$ -- resp. $(\hat{\mathcal{P}}_K)$ -- associated to the stopping times $\tau, \hat{\tau}$ respectively, the following holds,
    $$\mathbbm{P}\Big(\sup_{[0,T]} \|U_{\tau' \wedge \bcdot}^* - \hat{U}_{\tau' \wedge \bcdot}^*\|_H^2 = 0 ;\quad \forall T>0 \Big)=1,$$
    for all stopping time $0< \tau' < \tau \wedge \hat{\tau}$.
\end{enumerate}
\end{theorem}
Moreover, under stronger assumptions, we may state a global-in-time well-posedness result.
\begin{theorem} \label{theorem-filtered-strong}
    Using the same notations as in Theorem \ref{theorem-filtered-weak}, assume that $U_0 \in V$, and that the noise fulfils \eqref{smoothness-noise} and \eqref{eq-noiseBC-weak}. Then, the following propositions hold,
    \begin{enumerate}
        \item $(\hat{\mathcal{P}}_K)$ admits a global-in-time pathwise solution, which is unique up to indistinguishability, and almost surely belongs to the space $$L_{loc}^2\Big([0,+\infty),\mathcal{D}(A)\Big) \cap C\Big([0,+\infty),V\Big).$$
        \item This solution is continuous in the following sense: let $T>0$ and define $\Sigma$, a space of noise operators, as follows,
            $$\Sigma = \Bigg\{ \sigma \in \mathcal{L}_2\Big(\mathcal{W}, H^4(\mathcal{S}, \R^3)\Big) \: \Bigg| \: a \in H^1\Big([0,T],H^3(\mathcal{S}, \R^{3 \times 3})\Big) \Bigg\}.$$
        Let $U_0^n \in V$ a sequence of initial conditions, and $\sigma^n \in \Sigma$ a sequence of noise operators, such that
        $$U_0^n \rightarrow U_0 \text{ in $V$, and } \sigma^n \rightarrow \sigma \text{ in $\Sigma$.}$$
        Denote by $U^n$ the solution to $(\hat{\mathcal{P}}_K)$, associated to the initial condition $U_0^n$ and the noise operator $\sigma^n$. Similarly, denote by $U$ the solution associated to $U_0$ and $\sigma$. Then,
        \begin{equation}
            U^n \rightarrow U \text{ in probability, in the space $L^2([0,T],\mathcal{D}(A)) \cap C([0,T],V)$.}
        \end{equation}
    \end{enumerate} 
\end{theorem}
\begin{remark}\hfill
    \begin{itemize}
        \item In particular, for a fixed initial data, if we denote by $U^{\sigma}$ the solution of $(\hat{\mathcal{P}}_K)$ associated to the noise data $\sigma \in \Sigma$, then, in probability,
        \begin{equation*}
            U^{\Upsilon^{1/2} \sigma}\underset{ \Upsilon \rightarrow 0}{\longrightarrow} U^0, \text{ in } L^2([0,T],\mathcal{D}(A)) \cap C([0,T],V).
        \end{equation*}
        Here, $U^0$ denotes the solution to the problem $(\hat{\mathcal{P}}_K)$ with noise zero -- i.e. $\sigma = 0$ -- that is to say the classical deterministic primitive equations.
        \item In particular again, assuming that the noise $\sigma dW_t$ is purely barotropic, Theorem 2 holds for $(\mathcal{P}_K)$, because in such case $(\hat{\mathcal{P}}_K)$ and $(\mathcal{P}_K)$ are equivalent. Using the formalism of subsection \ref{approx-low-pass-filtered-problem}, pathwise solutions of the problem $(\hat{\mathcal{P}}_K)$ are also continuous as $\sigma_{BC}$ varies in $\Sigma$. Denote by $\hat{U}^{\sigma_{BT}, \Upsilon^{1/2} \sigma_{BC}}$ and $U^{\sigma_{BT}}$ the unique pathwise solutions of $(\hat{\mathcal{P}}_K)$ and $(\mathcal{P}_K)$, with noise operators $\hat{\sigma} =\sigma_{BT} + \Upsilon^{1/2} \sigma_{BC}$ and $\sigma =\sigma_{BT}$, respectively. Consequently, the following convergence holds in probability whenever $\Upsilon_n \rightarrow 0$,
        \begin{equation*}
            \hat{U}^{\sigma_{BT}, \Upsilon_n^{1/2} \sigma_{BC}} \underset{ n \rightarrow \infty}{\longrightarrow} U^{\sigma_{BT}}.
        \end{equation*}
    \end{itemize}
\end{remark}
Moreover, the following weaker existence properties hold for the problems $(\mathcal{P}^{\gamma_r}_{EB})$ and $(\mathcal{P}^{\gamma_r}_{EV})$ -- equations \eqref{problem-PEB} and \eqref{problem-PEV} respectively,
\begin{theorem} \label{theorem-EB}
    Suppose that $\alpha >0$ and $\gamma_r > 1$. Then the problem $(\mathcal{P}^{\gamma_r}_{EB})$, equipped with the initial condition $U_0 \in H$, admits at least one global-in-time martingale solution, for all $T>0$, in the space $$L^2\Big(\Omega,L^2([0,T],V)\Big) \cap L^2\Big(\Omega,L^\infty([0,T],H)\Big).$$
\end{theorem}

\begin{theorem} \label{theorem-EV}
    Suppose that $\alpha >0$ and $\gamma_r > 2$. Then the problem $(\mathcal{P}^{\gamma_r}_{EV})$, equipped with the initial condition $U_0 \in H$, admits at least one global-in-time martingale solution, for all $T>0$, in the space $$L^2\Big(\Omega,L^2([0,T],V)\Big) \cap L^2\Big(\Omega,L^\infty([0,T],H)\Big).$$
\end{theorem}
To summarise, we show a well-posedness result that is similar to the work of \cite{AHHS_2022}, but with rigid-lid boundary conditions instead of the water world ones, and where we consider additional terms stemming from our discussion about the relaxation of the hydrostatic assumption in a stochastic context. We detail the proof of the first point of Theorem \ref{theorem-filtered-weak} in section \ref{section-martingale}, and of its second point in section \ref{section-loc-wpn}. Theorem \ref{theorem-filtered-strong} is shown in section \ref{section-globality}. For Theorems \ref{theorem-EB} and \ref{theorem-EV}, we only explicit the arguments where the proofs change, since, in $(\mathcal{P}^{\gamma_r}_{EV})$ and $(\mathcal{P}^{\gamma_r}_{EB})$, only the pressure terms differ from problem $(\mathcal{P}_K)$. This is sketched in the appendix.

To be more precise, the sketch of our proof for Theorem \ref{theorem-filtered-weak} is similar to the one used in \cite{BS_2021, DGHT_2011, DHM_2023}: we consider a Galerkin approximation of the problem, then, using energy estimates and tightness arguments, we show that a subsequence of its solutions converges toward a solution of the initial problem. Our scheme of proof is closer in spirit to \cite{BS_2021} than to \cite{AHHS_2022}. Moreover, our proofs adapt when the vertical acceleration $\mathbbm{D}_t w$ is totally neglected, which corresponds to choosing $K=0$ and neglecting additive noises in the pressure equations of the problem $(\mathcal{P}_K)$. In such case, the problem is equivalent to the one with strong hydrostatic assumption.

\section{Existence of martingale solutions} \label{section-martingale}

In this section, we aim to show the first point of Theorem \ref{theorem-filtered-weak} concerning $(\mathcal{P}_K)$ -- equation \eqref{problem-PK}. For $(\hat{\mathcal{P}}_K)$, the additional terms will be treated 
\begin{itemize}
    \item in Remarks \ref{F-hat-sigma-Lipschitz} and \ref{G-hat-sigma-Lipschitz}, concerning preliminary results,
    \item in Remarks \ref{energy-estimate-boundedVar-P-hat-K} and \ref{energy-estimate-martingale-P-hat-K}, concerning bounds for the energy estimates. 
\end{itemize}
Moreover, we quickly discuss the influence of non divergence-free noises on the energy estimates in Remarks \ref{remark-energ-noise-dep-v} and \ref{remark-sup-noise}. To be precise, considering such noises would affect the expression of the It\={o}-Stokes drift, which would become 
$$u_s = \sum_k \phi_k^\tr (\nabla \bcdot \phi_k) - \frac{1}{2} \nabla \bcdot (\phi_k \phi_k^\tr).$$
All other conditions being fulfilled, the following results hold with the noise regularity assumptions \eqref{smoothness-noise}. To emphasise the difference between $\nabla_H$ and $\nabla$, from now on we use the notation $\nabla_3$ instead of $\nabla$. For similar reasons, we use $\sigma_3$ and $a_3$ for $\sigma$ and $a$ respectively. In addition, for notational convenience, we may write the transport terms $f \bcdot \nabla_3 g$ instead of $(f \bcdot \nabla_3) g$ for two differentiable vectors $f$ and $g$. Moreover, we remind that $K$ has been defined as a regularising kernel of regularity $H^3$. In particular, the operator $\mathcal{C}_K (\bcdot) := K * (\bcdot)$ satisfies
\begin{equation}
    \|\mathcal{C}_K f \|_{H^k} \leq C (1+\|f\|_{L^2}), \quad \forall k \in \{0,1,2,3\}. \label{regularisation-CK}
\end{equation}
The same property holds for the operator $\mathcal{C}_{\check{K}}$.

As mentioned earlier, the sketch of our proof is similar to the approach used in \cite{BS_2021, DGHT_2011, DHM_2023}: first we consider a well-chosen range of Galerkin approximation problems, each of them admitting a unique solution. Then we show they fulfil a uniform energy estimate, and that their laws are tight. By application of the Prohorov and Skorohod theorems, and using limit theorems, this implies the existence of a martingale solution. Eventually, we show an additional regularity result on this solution.

\subsection{Preliminaries} \label{subsec-prelim}

Denote by $P^n$ the projection onto $H_n = Span(e_i)_{i\leq n} \subset H$, the space generated by the $n$ first eigenvectors $(e_i)_{i\leq n}$ of $A$ (see below for their precise definition), and define the following $H_n$-valued operators,
\begin{gather*}
    B^n = P^n B, \quad \Gamma^n = P^n \Gamma, \quad F_\sigma^n = P^n F_\sigma, \quad G_\sigma^n = P^n G_\sigma,\\
    \Big(\frac{1}{\rho_0} \mathbf{P} \nabla_H p\Big)^n = P^n \Big(\frac{1}{\rho_0} \mathbf{P} \nabla_H p\Big), \quad \Big(\frac{1}{\rho_0} \mathbf{P} \nabla_H dp_r^\sigma\Big)^n = P^n \Big(\frac{1}{\rho_0} \mathbf{P} \nabla_H dp_r^\sigma\Big).
\end{gather*}
Then, denote by $(\mathcal{P}_n)$ the projection of the problem $(\mathcal{P}_K)$ onto $H_n$, so that
\begin{multline}
    U_n^*  = P^n(U_0^*) - \int_0^t [A U_n^* + B^n (U_n^*) + \Gamma^n U_n^* + \Big(\frac{1}{\rho_0} \mathbf{P} \nabla_H p\Big)^n + F_\sigma^n (U_n^*)] dr\\ + \int_0^t [G_\sigma^n(U_n^*) dW_r - \Big(\frac{1}{\rho_0} \mathbf{P} \nabla_H dp_r^\sigma\Big)^n]. \: \: \: (\mathcal{P}_n)
\end{multline}
The problems $(\mathcal{P}_n)$ are equipped with the initial condition $P^nU_0 \in H_n$. Applying the Cauchy-Lipschitz theorem for stochastic differential equations in finite dimension, we know that each problem $(\mathcal{P}_n)$ admits a local-in-time solution $(U_n^*,t_n)$, where $U_n^* \in H_n$. The objective of this section is to prove continuity results on the non projected operators involved in the abstract problem $(\mathcal{P}_K)$, to be used in the next steps of the proof.

\paragraph{Diffusion, advection and Coriolis operators ($A, B$ and $\Gamma$)} \hfill

\noindent
Begin with the diffusion operator $A$. Remind first that $\mathcal{D}(A) = V \cap H^2(\mathcal{S},\R^4)$, then for all $U,U^\sharp \in V$, we have
\begin{align}
    (AU,U^\sharp)_H &= \int_\mathcal{S} D_\mu \partial_x U \bcdot \partial_x U^\sharp + \int_\mathcal{S} D_\mu \partial_y U \bcdot \partial_y U^\sharp + \int_\mathcal{S} D_\nu \partial_{z}U \bcdot \partial_{z}U^\sharp,
\end{align}
with $D_\mu = (\mu_v, \mu_T, \mu_S)^\tr$, $D_\nu = (\nu_v, \nu_T, \nu_S)^\tr$, since $A = (A^v, A^T, A^S)^\tr$. Consequently,
\begin{align}
    3\min_i \{\mu_i,\nu_i\} \|U\|_V^2 \leq |(AU,U)_H| \quad \text{ and} \quad |(AU,U^\sharp)_H| \leq 3\max_i \{\mu_i,\nu_i\} \|U\|_V \|U^\sharp\|_V, \label{bound-A}
\end{align}
so in particular $\|AU\|_{V'} \leq C \|U\|_V$. Since $V \hookrightarrow H$ is a compact embedding using Sobolev theorem, we deduce that the embedding $\mathcal{D}(A) \hookrightarrow H$ is also compact. Consequently $A : \mathcal{D}(A) \longrightarrow H$ is a closed positive self-adjoint operator with compact resolvent. By spectral theorem, there exists a family $(e_i)$ of eigenvectors of $A$, such that $(e_i)$ is a Hilbert basis of $H$. Also, each eigenvector $e_i$ is associated to an eigenvalue $\lambda_i\in \R$, the family $(\lambda_i)$ being increasing and unbounded. We also define, for $\alpha>0$,
\begin{align}
    \mathcal{D}(A^\alpha) = \{u \in H | \sum_{k=0}^\infty \lambda_k^{2\alpha} |(u,e_k)|^2 < \infty \},
\end{align}
and remark that $(e_i) \in \mathcal{D}(A^\alpha), \forall \alpha>0$. Now equip this space with the inner product
\begin{align}
    (u,v)_{\mathcal{D}(A^\alpha)} := \sum_{k=0}^\infty \lambda_k^{2\alpha} (u,e_k)(v,e_k), \quad \forall u,v \in \mathcal{D}(A^\alpha).
\end{align}
One can show that, if $\alpha > 0$, then the norms $\|\bcdot\|_{\mathcal{D}(A^\alpha)}$ and $\|\bcdot\|_{H^{2\alpha}}$ are equivalent on $\mathcal{D}(A)$ since $\mathcal{S}$ is smooth. As a consequence, $\mathcal{D}(A^\alpha)$ is a Hilbert space.

To study the operator $B$, we define, for $U,U^\flat \in V$ and $U^\sharp \in \mathcal{D}(A)$,
\begin{align}
    b(U,U^\sharp, U^\flat) = \begin{pmatrix}
        \int_\mathcal{S} v^\flat (v \bcdot \nabla) v^\sharp  + v^\flat w(v) \partial_z v^\sharp\\
        \int_\mathcal{S} T^\flat (v \bcdot \nabla) T^\sharp  + T^\flat w(v) \partial_z T^\sharp\\
        \int_\mathcal{S} S^\flat (v \bcdot \nabla) S^\sharp  + S^\flat w(v) \partial_z S^\sharp\\
    \end{pmatrix}
    =\int_\mathcal{S} U^\flat (v \bcdot \nabla) U^\sharp + \int_\mathcal{S} U^\flat w(v) \partial_z U^\sharp.
\end{align}
Using integration by parts and that $(v, w(v))^\tr$ is divergence-free, remark that $b(u,\bcdot, \bcdot)$ is anti-symmetric, that is
\begin{align}
    b(U,U^\sharp, U^\flat)=-b(U,U^\flat,U^\sharp) = \int_\mathcal{S} U^\flat (v \bcdot \nabla) U^\sharp + \int_\mathcal{S} U^\flat w(v) \partial_z U^\sharp.
\end{align}
To show the continuity of $b$, we split it into two terms $b_1$ and $b_2$ as follows,
\begin{align}
    b_1 := \int_\mathcal{S} U^\flat (v \bcdot \nabla) U^\sharp, \qquad  b_2 := \int_\mathcal{S} U^\flat w(v) \partial_z U^\sharp.
\end{align}
Regarding $b_1$, using Cauchy-Schwarz and Hölder inequality yields,
\begin{align*}
    |b_1(U,U^\sharp, U^\flat)| \leq \int_\mathcal{S} \sum_{i=1}^2\sum_{j=1}^4 |U_j^\flat v_i|  |\partial_i U_j^\sharp| &\leq \|U^\flat \otimes U\|_{L^2(\mathcal{S},\R^{2\times4})} \|\nabla U^\sharp\|_{L^2(\mathcal{S},\R^{2\times4})}\\ &\leq C \|U^\flat\|_{L^4(\mathcal{S},\R^4)} \|U\|_{L^4(\mathcal{S},\R^4)} \|U^\sharp\|_{V},
\end{align*}
and by Sobolev theorem $\|U'\|_{L^4(\mathcal{S},\R^4)} \leq C \|U'\|_V, \: \forall U' \in V.$ Thus
\begin{align*}
    |b_1(U,U^\sharp,U^\flat)| \leq C \|U^\flat\|_{V} \|U\|_{V} \|U^\sharp\|_{V}.
\end{align*}
Regarding $b_2$, 
\begin{align*}
    |b_2(U,U^\sharp,U^\flat)| \leq C\|w(U)\|_{L^2} \|\partial_z U^\sharp\|_{L^4(\mathcal{S},\R^4)} \|U^\flat\|_{L^4(\mathcal{S},\R^4)} \leq \|U\|_V \|U^\sharp\|_{\mathcal{D}(A)} \|U^\flat\|_V,
\end{align*}
using similar arguments. Consequently, we reach
\begin{align}
    |b(U,U^\sharp,U^\flat)| &\leq \|U\|_V \|U^\sharp\|_{\mathcal{D}(A)} \|U^\flat\|_{V}.
\end{align}
As demonstrated in \cite{PTZ08} (Lemma 3.1), we also have
\begin{align}
    |b(U,U^\sharp, U^\flat)| &\leq C \|U^\flat\|_H \| U^\sharp\|_V^{1/2} \| U^\sharp \|_{\mathcal{D}(A)}^{1/2} \| U\|_V^{1/2} \| U\|_{\mathcal{D}(A)}^{1/2}. \label{majoration PTZ}
\end{align}
Now, let $B(U,U^\sharp) = \mathbf{P} b(U,U^\sharp,\bcdot)$ and $B(U)=B(U,U)$, so that the notation is consistent with the one we proposed earlier. The previous results yield
\begin{align}
    |B(U,U^\sharp)|_{V'} &\leq C\|U^\sharp\|_{\mathcal{D}(A)} \|U\|_V, \text{ and } |B(U,U^\sharp)|_{\mathcal{D}(A^{-1})} \leq C\|U^\sharp\|_{V} \|U\|_V,
\end{align}
that is the operator $B$ is continuous from $V \times V$ to $\mathcal{D}(A^{-1})$, and from $V \times \mathcal{D}(A)$ to  $V'$. To extend this operator to $H\times H$, using the Sobolev theorem, we claim that the embedding
\begin{equation}
    H^\beta(\mathcal{S},\R^4) \hookrightarrow L^\infty(\mathcal{S},\R^4) \text{ is continuous}
\end{equation}
when $\beta > 3/2$. So the embedding
\begin{equation}
    \mathcal{D}(A^{\beta/2}) \hookrightarrow L^\infty(\mathcal{S},\R^4) \text{ is continuous as well,}
\end{equation}
with the same condition on $\beta$. As a result, if $U,U^\sharp \in V$ and $U^\flat \in \mathcal{D}(A^{\frac{\beta+2}{2}})$, then
\begin{align*}
    |b(U,U^\sharp,U^\flat)|=|b(U,U^\flat,U^\sharp)| &\leq C\|A U^\flat\|_{L^\infty(\mathcal{S},\R^4)} \|U^\sharp\|_H \|U\|_H\\
    &\leq C\|A U^\flat\|_{\mathcal{D}(A^{\beta/2})} \|U^\sharp\|_H \|U\|_H\\
    &\leq C\|U^\flat\|_{\mathcal{D}(A^{\frac{\beta+2}{2}})} \|U^\sharp\|_H \|U\|_H.
\end{align*}
Hence, $b$ can be extended by continuity to $H \times H \times \mathcal{D}(A^\gamma) \longrightarrow \R$, that is $B$ can be extended on
$H \times H \longrightarrow \mathcal{D}(A^{-\gamma})$, where $\gamma := \frac{\beta+2}{2} > \frac{5}{4}$, and satisfies 
\begin{align}
    \|B(U,U^\sharp)\|_{\mathcal{D}(A^{-\gamma})} \leq C\|U^\sharp\|_H \|U\|_H. \label{eq-extension-B}
\end{align}
We deduce that the operator $B$ is locally Lipschitz in $H \times H \longrightarrow \mathcal{D}(A^{-\gamma})$. Additionally, remark that
\begin{align}
    |(\Gamma U,U^\sharp)_H| &\leq \int_\mathcal{S} |f(k\times v) \bcdot v^\sharp| \leq \|f\|_{L^\infty(\mathcal{S}, \R)} \|v\|_{H} \|v^\sharp\|_{H}\leq C \|U\|_{H} \|U^\sharp\|_{H}, \:\:\: \forall U^\sharp \in H.
\end{align}

\paragraph{Bounded variation pressure operator ($\nabla_H p$)} \hfill

\noindent
Let $U \in H$ and $U^\sharp \in V$. For $\frac{1}{\rho_0} \mathbf{P} \nabla_H p$,
\begin{align}
    \Big(\frac{1}{\rho_0} \mathbf{P} \nabla_H p,U^\sharp \Big)_H = &-\Big( g \:  \nabla_H \int_z^0 {(\beta_T T + \beta_S S) dz'}, v^\sharp \Big)_{L^2} \nonumber\\
    &- \Big( \nabla_H \int_z^0 \Big[K*(u_s \bcdot \nabla_3 w(v^*)) + \frac{1}{2} \nabla_3 \bcdot (a_3^K \nabla_3 w(v^*))\Big] dz',v^\sharp \Big)_{L^2} \nonumber\\
    &- \Big( \nabla_H \int_z^0 \Big[K*(u_s \bcdot \nabla_3 w_s) + \frac{1}{2} \nabla_3 \bcdot (a_3^K \nabla_3 w_s) \Big] dz',v^\sharp \Big)_{L^2}.
\end{align}
Moreover,
\begin{align*}
    \Big|\Big(\nabla_H \int_z^0 {g(\beta_T T + \beta_S S) dz'}, v^\sharp \Big)_{L^2}\Big| = \Big|\Big(g(\beta_T T + \beta_S S), w(v^\sharp) \Big)_{L^2} \Big| &\leq C \|U\|_H \|U^\sharp\|_V,\\
    \Big|\Big( \nabla_H \int_z^0 \Big[ K*(u_s \bcdot \nabla_3 w_s) + \frac{1}{2} \nabla_3 \bcdot (a_3^K \nabla_3 w_s) \Big] dz',v^\sharp \Big)_{L^2}\Big| &\leq C \|U^\sharp\|_{L^2},
\end{align*}
where the latter inequality holds immediately by the noise regularity assumption \eqref{smoothness-noise}, since the left factor in the last inner product depends only on $v^\sharp$ and stochastic parameters -- that are $u_s = (v_s \quad  w_s)^\tr = \frac{1}{2} \nabla_3 \bcdot a$ and $a_3^K$. In addition, for $U^\sharp \in \mathcal{D}(A^{3/2})$,
\begin{align*}
    &\Big|\Big( \nabla_H \int_z^0 \Big[K*(u_s \bcdot \nabla_3 w(v^*)) + \frac{1}{2} \nabla_3 \bcdot (a_3^K \nabla_3 w(v^*))\Big] dz',v^\sharp \Big)_{L^2}\Big|\\
    &\leq \Big|\Big( K*(u_s \bcdot \nabla_3 w(v^*)) + \frac{1}{2} \nabla_3 \bcdot (a_3^K \nabla_3 w(v^*)), w(v^\sharp)\Big)_{L^2}\Big|\\
    &\leq \Big|\Big(u_s \bcdot \nabla_3 w(v^*), \check{K}*w(v^\sharp) \Big)_{L^2}\Big| + \frac{1}{2} \Big|\Big( w(v^*), \nabla_3 \bcdot (a_3^K \nabla_3 w(v^\sharp)) \Big)_{L^2}\Big| \\
    &\leq \underbrace{\Big|\Big(v^*, \nabla_H \int_z^0  \nabla_3 ((\check{K}*w(v^\sharp)) \otimes u_s) dz'\Big)_{L^2}\Big|}_{\leq C \|U^*\|_H \|u_s\|_{W^{2,\infty}(\mathcal{S})} \|\check{K}*w(v^\sharp)\|_{\mathcal{D}(A)}} + \frac{1}{2} \Big|\Big( v^*, \nabla_H \int_z^0 \nabla_3 \bcdot (a_3^K \nabla_3 w(v^\sharp)) dz' \Big)_{L^2} \Big|\\
    &\leq C\|U^*\|_H (\|\check{K}*w(v^\sharp)\|_{\mathcal{D}(A)} + \|U^\sharp\|_{\mathcal{D}(A^{3/2})}) \leq C\|U^*\|_H \|U^\sharp\|_{\mathcal{D}(A^{3/2})},
\end{align*}
where we used that the embedding $H^4(\mathcal{S}) \hookrightarrow W^{2,\infty}(\mathcal{S})$ is compact by Sobolev theorem. Therefore,
\begin{equation}
    \|\frac{1}{\rho_0} \mathbf{P} \nabla_H p\|_{\mathcal{D}(A^{-3/2})} \leq C (1 + \|U\|_H).
\end{equation}
Consequently $\frac{1}{\rho_0} \mathbf{P} \nabla_H p : H \rightarrow \mathcal{D}(A^{-3/2})$ is continuous and Lipschitz.

\paragraph{Noise related operators ($F_\sigma, G_\sigma$ and $\nabla_H p_t^\sigma$):}

Remind that $\phi_k(t)$ is a collection of eigenvectors of the operator $\sigma_3$, so that one can decompose $\sigma_3 dW_t$ and $a_3$ as follows,
\begin{align}
    \sigma_3 dW_t = \sum_{k=0}^\infty \phi_k d\beta_t^k, \quad a_3 = \sum_{k=0}^\infty \phi_k \phi_k^\tr.
\end{align}
Also, remind that the noise $\sigma_3 dW_t$ is regular in the following sense: for all $T \in \R_+$,
\begin{align*}
    && \sup_{t\in [0,T]} \sum_{k=0}^\infty \|\phi_k\|_{H^4(\mathcal{S},\R^3)}^2 < \infty,\\
    U_s \in L^\infty([0,T], H^4(\mathcal{S},\R^4)), && d_t U_s \in L^\infty([0,T],H^1(\mathcal{S},\R^4)), && a_3\nabla U_s \in L^\infty([0,T],H^2(\mathcal{S},\R^{4})). \nonumber
\end{align*}
Therefore, by Sobolev embedding theorem,
\begin{align}
    \sup_{t\in [0,T]} \sum_{k=0}^\infty \|\phi_k\|_{L^\infty(\mathcal{S},\R^3)}^2 &< \infty, \quad \sup_{t\in [0,T]} \|a_3\|_{H^1(\mathcal{S},\R^3)} + \|a_3\|_{L^\infty(\mathcal{S},\R^{3\times 3})} < \infty, \\
    \sup_{t\in [0,T]} \sum_{k=0}^\infty |B(\phi_k,v_s)|^2 &< \infty, \quad \sup_{t \in [0,T]} \|u_s\|_{H^4(\mathcal{S},\R^{3\times 3})} < \infty. \nonumber
\end{align}
Then we want to show that $F_\sigma$ is Lipschitz. Let $U \in H$ and $U^\sharp \in \mathcal{D}(A^\gamma)$, then
\begin{align}
    (F_\sigma (U) dt, U^\sharp)_H &= \Bigg(d_t U_s + \mathbf{P} \Big[ B(v,U_s) - \frac{1}{2} \nabla_3 \bcdot (a_3 \nabla_3 U_s) + AU_s + \Gamma U_s \Big] dt - \mathbf{P} \frac{1}{2} \nabla_3 \bcdot (a_3 \nabla_3 U^*) dt, U^\sharp\Bigg)_H \nonumber\\
    &\leq C (\|U\|_H+1) \|U^\sharp\|_{\mathcal{D}(A)} dt + ( B(v,U_s), U^\sharp)_H dt
\end{align}
Notice that
\begin{align*}
    |( B(v,U_s), U^\sharp)_H| \leq C \|U\|_H \|U_s\|_H \|U^\sharp\|_{\mathcal{D}(A^\gamma)} .
\end{align*}
Therefore,
\begin{equation}
    \|F_\sigma (U)\|_{\mathcal{D}(A^{-\gamma})} \leq C (\|U\|_H+1).
\end{equation}
So $F_\sigma : H \rightarrow \mathcal{D}(A^{-\gamma})$ is continuous and Lipschitz.
\begin{remark}\label{F-hat-sigma-Lipschitz}
    Considering the approximated problem $(\hat{\mathcal{P}}_K)$, we may split the operator $\hat{F}_\sigma(U^*)$ as
    \begin{equation}
        \hat{F}_\sigma(U^*) = F_\sigma^{BT}(U^*) - \frac{1}{2} \begin{pmatrix}
        \nabla \bcdot (\overline{a_{BC}} \nabla \overline{v}^*)\\
        0\\
        0
    \end{pmatrix} - \frac{1}{2} \begin{pmatrix}
        \nabla \bcdot (\overline{a_{BC}} \nabla \overline{v}_s)\\
        0\\
        0
    \end{pmatrix},
    \end{equation}
    where $F_\sigma^{BT}$ stands for the operator $F_\sigma$ of the problem $(\mathcal{P}_K)$ with $\sigma = \sigma_{BT}$ -- i.e. $a = a_{BT}$. By similar arguments, it is immediate that $\hat{F}_\sigma : H \rightarrow \mathcal{D}(A^{-\gamma})$ is continuous and Lipschitz as well.
\end{remark}
We show the same for the noise operators $G_\sigma$ and $p_t^\sigma$. Let $\Psi \in \mathcal{W}$, so that
\begin{align}
    (G_\sigma(U)\Psi,U^\sharp)_H = -(\sigma_3 \Psi \bcdot \nabla_3 U + \sigma_3 \Psi \bcdot \nabla_3 U_s + A (\sigma^H \Psi) + \Gamma (\sigma^H \Psi),U^\sharp)_H.
\end{align}
Thus,
\begin{align*}
    |(G_\sigma(U)\Psi,U^\sharp)_H| \leq C[\|\sigma \Psi\|_H \|U^\sharp\|_{\mathcal{D}(A)} + \|\sigma \Psi\|_H \|U\|_H \|U^\sharp\|_{\mathcal{D}(A^\gamma)} ],
\end{align*}
that is
\begin{equation}
    \|G_\sigma(U)\Psi\|_{\mathcal{D}(A^{-\gamma})} \leq C (1+ \|U\|_H)\|\sigma \Psi\|_H.
\end{equation}
Hence, $G_\sigma \Psi : H \rightarrow \mathcal{D}(A^{-\gamma})$ is Lipschitz, so $G_\sigma : H \rightarrow \mathcal{L}_2(\mathcal{W}, \mathcal{D}(A^{-\gamma}))$ is Lipschitz as well, since $G_\sigma$ is linear with respect to $\sigma$.

\begin{remark}\label{G-hat-sigma-Lipschitz}
    Considering the approximated problem $(\hat{\mathcal{P}}_K)$ again, we may split $\hat{G}_\sigma(U^*) dW_t$ as
    \begin{equation}
        \hat{G}_\sigma(U^*) dW_t = G_\sigma^{BT}(U^*)dW_t - \begin{pmatrix}
        (\sigma_{BC} dW_t \bcdot \nabla) \overline{v}^*\\
        0\\
        0
    \end{pmatrix} - \begin{pmatrix}
        (\sigma_{BC} dW_t \bcdot \nabla) \overline{v}_s\\
        0\\
        0
    \end{pmatrix},
    \end{equation}
     where $G_\sigma^{BT}$ stands for the operator $G_\sigma$ of the problem $(\mathcal{P}_K)$ with $\sigma = \sigma_{BT}$. Then, it is immediate that $\hat{G}_\sigma : H \rightarrow \mathcal{L}_2(\mathcal{W}, \mathcal{D}(A^{-\gamma}))$ is continuous and Lipschitz.
\end{remark}
Moreover, for $U^\sharp \in H$,
\begin{align}
    (&\frac{1}{\rho_0} \mathbf{P} \nabla_H p_t^\sigma (\Psi),v^\sharp) = (\nabla_H \int_z^0 \Big[K *[\sigma_3 \Psi \bcdot \nabla_3 w(v) + \sigma_3 \Psi \bcdot \nabla_3 w_s] + A (\sigma^z \Psi) \Big]dz',v^\sharp)_{L^2} \nonumber\\
    &= (K*[\sigma_3 \Psi \bcdot \nabla_3 w(v) + \sigma_3 \Psi \bcdot \nabla_3 w_s] + A (\sigma^z \Psi), w(v^\sharp))_{L^2}\nonumber\\
    &= -( v, \nabla_H \int_z^0 \nabla_3 [(\sigma_3 \Psi) \check{K}* w(v^\sharp)] )_{L^2} + (\sigma_3 \Psi \bcdot \nabla_3 w_s, \check{K} * w(v^\sharp))_{L^2} + (A (\sigma^z \Psi), w(v^\sharp))_{L^2} \nonumber\\
    &\leq C(\|U^*\|_H + 1)\|\sigma_3 \Psi\|_{H^4(\mathcal{S},\R^3)} \|U^\sharp\|_{H},
\end{align}
where we used again the boundary conditions on $\mathcal{A}v^\sharp$ to cancel the martingale surface pressure term. Therefore,
\begin{equation}
    \|\frac{1}{\rho_0} \mathbf{P} \nabla_H p_t^\sigma (\Psi)\|_{H} \leq C(\|U^*\|_H + 1)\|\sigma_3 \Psi\|_{H^4(\mathcal{S},\R^3)},
\end{equation}
so that $\frac{1}{\rho_0} \nabla_H p_t^\sigma \Psi : H \rightarrow H$ is continuous and Lipschitz. Thus so is $\frac{1}{\rho_0} \nabla_H p_t^\sigma : H \rightarrow \mathcal{L}_2(\mathcal{W}, H)$, as $\frac{1}{\rho_0} \nabla_H p_t^\sigma$ is linear with respect to $\sigma$.

\subsection{Energy estimates in $L^\infty([0,T],H) \cap L^2([0,T],V)$ for all $T>0$}\label{subsec-energy-1}

In the rest of the paper, we assume that the constant $\gamma$ used for defining the extension of the operator $B$ in equation \eqref{eq-extension-B} fulfils $\gamma = 2$. Moreover, for notational convenience and without confusion, from now on we may write $W^{s,p}$, $H^s$ and $L^p$ instead of $W^{s,p}(\mathcal{S},\R^d)$, $H^s(\mathcal{S},\R^d)$ and $L^p(\mathcal{S},\R^d)$ respectively, for any $p \in [1, \infty]$, $s \in \R_+$ and $d \in \mathbbm{N}$. In addition, we use the convention under which $C$ refers to a positive constant, which may differ from one line to the other. The aim of this subsection is to prove the following lemma, which allows in particular to conclude that the solution $U_n^*$ of $(\mathcal{P}_n)$ is global-in-time.
\begin{lemma} \label{lemma-basic-estimate}
    Let $T >0$. If $p \geq 2$, then there exists a constant $C$, independent of $n$, such that,
    \begin{align}
        \E[\sup_{0\leq t \leq T} \|U_n^*\|_H^p ] &\leq C, \text{ and } \E[\int_0^{T} \|U_n^*\|_H^{p-2} \|U_n^*\|_V^2 dt] \leq C. \label{estimation-énerg}
    \end{align}
\end{lemma}
\emph{Proof:} Let $T>0$ and $t \in [0,T]$, and apply It\={o}'s lemma to $\|\bcdot\|_H^p$,
\begin{align}
    d_t \|U_n^*\|_H^p = &- p \|U_n^*\|_H^{p-2} \Big(U_n^*, A U_n^* + B^n(U_n^*) + \Gamma^n U_n^* + \Big(\frac{1}{\rho_0} \mathbf{P} \nabla_H p\Big)^n + F_\sigma^n (U_n^*)\Big)_H dt \nonumber\\ 
    & + p \|U_n^*\|_H^{p-2} \Big(U_n^*, G_\sigma^n(U_n^*)dW_t\Big)_H + \frac{p(p-2)}{2} \Big\|\Big[G_\sigma^n(U_n^*)- \mathcal{M}^{\sigma,U_n^*}\Big]^*U_n^*\Big\|_{H}^2 \|U_n^*\|_H^{p-4} dt \nonumber\\
    &+ \frac{p}{2} \Big\|G_\sigma^n(U_n^*) - \mathcal{M}^{\sigma,U_n^*} \Big\|_{\mathcal{L}_2(\mathcal{W},H)}^2 \|U_n^*\|_H^{p-2} dt \label{Itô-énergie},
\end{align}
where we denote by $\mathcal{M}^{\sigma,U_n^*}$ the linear operator
\begin{equation}
    \mathcal{M}^{\sigma,U^*} \Psi = \frac{1}{\rho_0} \mathbf{P} \nabla_H \sum_{k=0}^\infty (\Psi,\pi_k^{\sigma,U_n^*})_H \pi_k^{\sigma,U_n^*}, \quad \forall \Psi \in \mathcal{W},
\end{equation}
with $\pi_k^{\sigma,U^*} = \nabla_H \int_z^0 \Big[K *[\phi_k \bcdot \nabla_3 w(v^*)] + A (\phi_k^z) \Big]dz'$. This implies in particular
\begin{equation}
    \mathcal{M}^{\sigma,U_n^*} dW_t =  \frac{1}{\rho_0} \mathbf{P} \nabla_H dp_t^{\sigma,U_n^*}.
\end{equation}
\emph{Step 1:} We compute the contribution of the bounded variation terms first. Remind that, from \eqref{bound-A}, there exist two constants $c,C >0$ such that
\begin{align}
c\|U_n^*\|_V^2 \leq |(U_n^*,AU_n^*)_H| \leq C\|U_n^*\|_V^2, \quad \text{ and } \quad (U_n^*,\Gamma^n U_n^*)_H = \int_\mathcal{S} f(k\times v_n^*) \bcdot v_n^* = 0,
\end{align}
since $k \times v_n^*$ and $v_n^*$ are orthogonal. Moreover, by anti-symmetry of $B$,
\begin{align}
    (U_n^*,B^n U_n^*)_H = (U_n^*,B U_n^*)_H = 0.
\end{align}
Also 
\begin{align}
    \Big|\Big( &U_n^*, \frac{1}{\rho_0} \mathbf{P} \nabla_H p \Big)_H\Big| \leq \Big|\Big(g \:  \nabla_H \int_z^0 {(\beta_T T_n + \beta_S S_n) dz'} , v_n^* \Big)_{L^2}\Big| \nonumber\\
    &+ \Big|\Big( \nabla_H \int_z^0 \Big[K*(u_s \bcdot \nabla_3 w_s) + \frac{1}{2} \nabla_3 \bcdot (a_3^K \nabla_3 w_s) \Big] dz',v_n^* \Big)_{L^2}\Big|\nonumber\\
    &+ \Big|\Big(\nabla_H \int_z^0 K*(u_s \bcdot \nabla_3 w(v_n^*)), v_n^* \Big)_{L^2}\Big| + \frac{1}{2} \Big(\nabla_3 \bcdot (a_3^K \nabla_3 w(v_n^*)), w(v_n^*) \Big)_{L^2}\nonumber\\
    &\leq \Big|\Big(g (\beta_T T_n + \beta_S S_n) , w(v_n^*) \Big)_{L^2}\Big| + \Big|\Big(K*(u_s \bcdot \nabla_3 w(v_n^*)), w(v_n^*) \Big)_{L^2}\Big| \nonumber\\
    &+ \frac{1}{2} \Big(a_3^K \nabla_3 w(v_n^*), \nabla_3 w(v_n^*) \Big)_{L^2} + \Big|\Big(K*(u_s \bcdot \nabla_H w_s) + \frac{1}{2} \nabla_3 \bcdot (a_3^K \nabla_3 w_s), w(v_n^*) \Big)_{L^2}\Big| \nonumber\\
    &\qquad \qquad \qquad \qquad \qquad \qquad \qquad \qquad \qquad \qquad \qquad \qquad \qquad + C\|U_n^*\|_{V}, \label{eq-inner-product-pressure}
\end{align}
where we used the triangle inequality and the boundary conditions on $v_n^*$ and $w_n^*$ \eqref{boundary-conditions}. Notice that the second term on the RHS of \eqref{eq-inner-product-pressure} can be written as
\begin{align*}
    \Big|\Big(K*(u_s \bcdot \nabla_3 w(v_n^*)), w(v_n^*) \Big)_{L^2}\Big| &= \Big|\Big( u_s \bcdot \nabla_3 w(v_n^*), \check{K}*w(v_n^*) \Big)_{L^2}\Big| \\
    &\leq C \|u_s\|_{W^{2,\infty}(\mathcal{S})} \|U_n^*\|_H \|\check{K}*w(v_n^*)\|_{H^2}\\
    &\leq C \|U_n^*\|_H^2 \leq C \|U_n^*\|_H \|U_n^*\|_V.    
\end{align*}
Therefore, since the left factor in the fourth inner product on the RHS of equation \eqref{eq-inner-product-pressure} depends only on stochastic parameters -- $u_s$ and $a_3^K$ -- and since $\|w(v_n^*)\|_{L^2(\mathcal{S},\R)} \leq \|U_n^*\|_V$ from \eqref{w-definition},
\begin{equation}
     \Big|\Big( U_n^*,\frac{1}{\rho_0} \mathbf{P} \nabla_H p \Big)_H\Big| \leq \frac{1}{2} \Big(a_3^K \nabla_3 w(v_n^*), \nabla_3 w(v_n^*) \Big)_{L^2}  + C (1+\|U_n^*\|_H) \|U_n^*\|_V,
\end{equation}
thanks to the noise regularity hypothesis \eqref{smoothness-noise}. Remind that $a_3^K$ was defined in \eqref{def-aK} by
\begin{equation*}
    a^K f = \sum_{k=0}^\infty \phi_k \mathcal{C}_K^* \mathcal{C}_K(\phi_k^\tr f).
\end{equation*}
Additionally, the following relation holds,
\begin{align*}
    \Big(a_3^K \nabla_3 w(v_n^*), \nabla_3 w(v_n^*) \Big)_{L^2} &= \sum_{k=0}^\infty \Big(\mathcal{C}_K (\phi_k \bcdot \nabla_3 w(v_n^*)), \mathcal{C}_K (\phi_k \bcdot \nabla_3 w(v_n^*)) \Big)_{L^2} \\
    &= \sum_{k=0}^\infty \Big\|K* (\phi_k \bcdot \nabla_3 w(v_n^*)) \Big\|_{L^2}^2 \geq 0.
\end{align*}
For $F_\sigma$, we have
\begin{multline}
    (U_n^*, F_\sigma(U_n^*))_H =(d_t v_s + \Big[ B(v_n^*,v_s) - \frac{1}{2} \nabla_3 \bcdot (a_3 \nabla_3 v_s) \Big] dt + \Big[Av_s + \Gamma v_s \Big] dt - \frac{1}{2} \nabla_3 \bcdot (a_3 \nabla_3 v_n^*), v_n^*)_H\\
    = (d_t v_s, v_n^*)_H + (B(v_n^*,v_s),v_n^*)_H - \frac{1}{2} (\nabla_3 \bcdot (a_3 \nabla_3 v_s),v_n^*)_H + (Av_s + \Gamma v_s, v_n^*)_H + \frac{1}{2} (a_3 \nabla_3 v_n^*, \nabla_3 v_n^*)_H,
\end{multline}
using the boundary conditions \eqref{boundary-conditions} and \eqref{eq-noiseBC-weak}. Here the first, third and fourth terms are inner products between parameters and $v_n^*$, so are bounded by $C \|v_n^*\|_H$ by \eqref{smoothness-noise} and Cauchy-Schwarz inequality. Similarly, by anti-symmetry,
$$|(B(v_n^*,v_s),v_n^*)_H| = |(B(v_n^*,v_n^*),v_s)_H| \leq \|v_s\|_{L^\infty(\mathcal{S}, \R^2)} \|v_n^*\|_H \|v_n^*\|_V.$$
Consequently, 
\begin{equation}
    (U_n^*, F_\sigma(U_n^*))_H \leq C \|v_n^*\|_H (1 + \|v_n^*\|_V) + \frac{1}{2} (a_3 \nabla_3 v_n^*, \nabla_3 v_n^*)_H.
\end{equation}

\begin{remark}\label{energy-estimate-boundedVar-P-hat-K}
    Considering $(\hat{\mathcal{P}}_K)$, and using notations of Remark \ref{F-hat-sigma-Lipschitz}, we have by similar arguments,
    \begin{multline}
        (v_n^*,\hat{F}_\sigma(U_n^*))_H \leq |(v_n^*,F_\sigma^{BT}(U_n^*))_H| + \frac{1}{2}(\nabla_3 v_n^*, \overline{a_{BC}} \nabla_3 \bar{v}_n^*)_{L^2} + \frac{1}{2}|(\nabla_3 v_n^*, \overline{a_{BC}} \nabla_3 \bar{v}_s)_{L^2}|\\
        \leq C \|v_n^*\|_H (1+\|v_n^*\|_V) + \frac{1}{2} (a_{BT} \nabla_3 v_n^*, \nabla_3 v_n^*)_{L^2} + \frac{1}{2}(\overline{a_{BC}} \nabla_3 \bar{v}_n^*, \nabla_3 \bar{v}_n^*)_{L^2}.
    \end{multline}
    Hence, the bound on $(v_n^*,\hat{F}_\sigma(U_n^*))_H$ has a similar form.
\end{remark}

\noindent
\emph{Step 2:} Now we estimate the contribution of the martingale terms. First,
\begin{align}
    \Big\|\Big[G_\sigma^n&(U_n^*)- \mathcal{M}^{\sigma,U_n^*}\Big]^* U_n^*\Big\|_H^2 = \sum_{k=0}^\infty \Big[\Big(w(v_n^*), K*[\phi_k \bcdot \nabla_3 w(v_n^*)]\Big)_{L^2} + \underbrace{(w(v_n^*), A \phi_k^z )_{L^2}}_{\leq C \|\phi_k\|_{H^3(\mathcal{S},\R^3)} \|v_n^*\|_H} \nonumber\\
    &- (v_n^*,\phi_k \bcdot \nabla_3 v_n^*)_{L^2} - \underbrace{\Big(v_n^*, \phi_k \bcdot \nabla_3 v_s + \nabla_H \int_z^0 \phi_k \bcdot \nabla_3 w_s + A \phi_k^H + \Gamma \phi_k^H\Big)_{L^2}}_{\leq C \|\phi_k\|_{H^3(\mathcal{S},\R^3)} \|v_n^*\|_H} \Big]^2, \label{basic-energ-cross-term}
\end{align}
noticing that the second and fourth terms under the sum are directly bounded by $C\|v_n^*\|_H$ since they are inner products between $v_n^*$ and functions of stochastic parameters -- $u_s$ and $\phi_k = (\phi_k^H \quad \phi_k^z)^\tr$. Also, the third term under the sum cancels out using that $\phi_k$ is divergence-free and the boundary condition $\phi_k = 0$ on $\Gamma$. Moreover, denoting $\check{K} : x \in \mathcal{S} \mapsto K(-x) \in \R$, remark that
\begin{align*}
    \Big(&w(v_n^*), K*[\phi_k \bcdot \nabla_3 w(v_n^*)]\Big)_{L^2(\mathcal{S},\R)} = -\Big( \nabla_3 \bcdot [\phi_k (\check{K} * w(v_n^*))], w(v_n^*) \Big)_{L^2(\mathcal{S},\R)}\\
    &= -\Big( \phi_k \bcdot [ \nabla_3 (\check{K} * w(v_n^*))], w(v_n^*) \Big)_{L^2(\mathcal{S},\R)} - \Big( (\nabla_3 \bcdot \phi_k) (\check{K} * w(v_n^*)), w(v_n^*) \Big)_{L^2(\mathcal{S},\R)} \\
    &\leq \|\nabla_H \phi_k\|_{W^{1,\infty}(\mathcal{S},\R^{2\times3})} \|\nabla_H \check{K} * w(v_n^*) \|_{H^1(\mathcal{S},\R^2)} \| v_n^* \|_H.
\end{align*}
Also, by Young's convolution inequality,
\begin{align*}
    \|\nabla_H \check{K} * w(v_n^*)\|_{L^2(\mathcal{S},\R^2)} &\leq \|(\Delta_H \check{K}) * \int_z^0 v_n^* dz'\|_{L^2(\mathcal{S},\R^2)} \leq C \|\Delta_H \check{K} \|_{L^1(\mathcal{S},\R)} \|v_n^*\|_H,
\end{align*}
and
\begin{align*}
    \|\nabla_3 \nabla_H \check{K} * w(v_n^*)\|_{L^2(\mathcal{S},\R^{2 \times 3})} &\leq \|\nabla_H \bcdot (\nabla_3 \nabla_H \check{K}) * \int_z^0 v_n^* dz'\|_{L^2(\mathcal{S},\R^{2 \times 3})}\\
    &\leq C \|\Delta_H \nabla_3 \check{K} \|_{L^1(\mathcal{S},\R^3)} \|v_n^*\|_H.
\end{align*}
Thus, gathering the previous estimates,
\begin{align*}
    \Big\|\Big[G_\sigma^n(U_n^*)- \mathcal{M}^{\sigma,U_n^*}\Big]^*U_n^*\Big\|_H^2 \leq \sum_{k=0}^\infty C \|\phi_k\|_{H^3}^2 (1+\|v_n^*\|_H^4) \leq C (1+\|U_n^*\|_H^4).
\end{align*}
Therefore, 
\begin{align}
    \|U_n^*\|_H^{p-4} \Big\|\Big[G_\sigma^n(U_n^*)- \mathcal{M}^{\sigma,U_n^*}\Big]^*U_n^*\Big\|_H^2 \leq C(\|U_n^*\|_H^{p}+1).
\end{align}

\begin{remark} \label{remark-energ-noise-dep-v}
    Notice that, upon considering a non divergence-free noise $\phi_k$, the third term under the sum in \eqref{basic-energ-cross-term}, $(v_n^*,\phi_k \bcdot \nabla_3 v_n^*)_H$, is not zero anymore. However, if we only consider $p=2$, then the term is cancelled out in the It\={o} lemma \eqref{Itô-énergie}. Moreover, following \cite{AHHS_2022,AHHS_2022_preprint}, if we assume that the noise satisfies the following ``parabolicity condition''
    \begin{equation}
        \forall \xi \in \R^3, \forall (t,x) \in \R_+ \times \mathcal{S},  \quad \sum_k (\phi_k(t,x) \bcdot \xi)^2 \leq \chi |\xi|^2, \label{eq-parabolicity-condition}
    \end{equation}
    for some $\chi \in (0, \frac{2}{p-2} \min\{\mu_v,\nu_v\})$ -- see Assumption 3.1 in \cite{AHHS_2022} and \cite{AHHS_2022_preprint} -- then
    \begin{equation}
        \sum_k (v_n^*,\phi_k \bcdot \nabla_3 v_n^*)_H^2 \leq \chi \|v_n^*\|_H^2 \|v_n^*\|_V^2.
    \end{equation}
    Hence,
    \begin{equation}
        \|U_n^*\|_H^{p-4} \Big\|\Big[G_\sigma^n(U_n^*)- \mathcal{M}^{\sigma,U_n^*}\Big]^*U_n^*\Big\|_H^2 \leq C(\|U_n^*\|_H^{p}+1) + \chi \|U_n^*\|_H^{p-2} \|U_n^*\|_V^2.
    \end{equation}
    Therefore, for non divergence-free noises, the previous bound holds when $p=2$, and for all $p \geq 3$ such that the parabolicity condition holds. However, taking for example $p=3$ leads to $\chi \in (0, 2 \min\{\mu_v,\nu_v\})$, which can be interpreted as the noise being controlled by the viscosity term -- that is to say the noise is (very) small.
\end{remark}

\noindent
Additionally,
\begin{align}
    \frac{1}{2} \Big\|G_\sigma^n(U_n^*)- &\mathcal{M}^{\sigma,U_n^*}\Big\|_{\mathcal{L}_2(\mathcal{W},H)}^2 = \frac{1}{2} \sum_{k=0}^\infty \Big\| \mathbf{P} [\phi_k \bcdot \nabla_3 v_n^* + \phi_k \bcdot \nabla_3 v_s + A \phi_k^H + \Gamma \phi_k^H] \Big\|_H^2 \nonumber\\
    &+ \frac{1}{2} \sum_{k=0}^\infty \Big\|\mathbf{P} \nabla_H \pi_k \Big\|_H^2 + \sum_{k=0}^\infty \Big(\phi_k \bcdot \nabla_3 (v_n^*+v_s) + A \phi_k^H + \Gamma \phi_k^H, \mathbf{P} \nabla_H \pi_k\Big)_{L^2}.
\end{align}
Remark that
\begin{multline*}
    \Big\|\phi_k \bcdot \nabla_3 v_n^* + \phi_k \bcdot \nabla_3 v_s + A \phi_k^H + \Gamma \phi_k^H\Big\|_H^2 = \Big(\phi_k \phi_k^\tr \nabla_3 (v_n^* + v_s), \nabla_3 (v_n^* + v_s)\Big)_{L^2} \\
    +  \|A \phi_k^H + \Gamma \phi_k^H\|_H^2 + 2\Big(\phi_k \bcdot \nabla_3 (v_n^*+v_s), A \phi_k^H + \Gamma \phi_k^H\Big)_{L^2}\\
    \leq (\phi_k \phi_k^\tr \nabla_3 v_n^*, \nabla_3 v_n^*)_{L^2} + C\|\phi_k\|_{H^3}^2(1+ \|v_n^*\|_H),
\end{multline*}
and thus
\begin{align*}
    \frac{1}{2}\sum_{k=0}^\infty \Big\|\mathbf{P}[\phi_k \bcdot \nabla_3 v_n^* + \phi_k \bcdot \nabla_3 v_s + A \phi_k^H + \Gamma \phi_k^H] \Big\|_H^2 &\leq \frac{1}{2} (a_3 \nabla_3 v_n^*, \nabla_3 v_n^*)_H + C (1+ \|v_n^*\|_H).
\end{align*}
In addition, by Cauchy-Schwarz and Young inequalities, for all $\xi >0$, there exists $C_\xi >0$ such that
\begin{align*}
    \sum_{k=0}^\infty (\phi_k \bcdot &\nabla_3 (v_n^*+v_s) + A \phi_k^H + \Gamma \phi_k^H, \mathbf{P} \nabla_H \pi_k)_{L^2} \\
    &\leq \xi \sum_{k=0}^\infty \|\mathbf{P} [\phi_k \bcdot \nabla_3 (v_n^*+v_s)] + A \phi_k^H + \Gamma \phi_k^H\|_H^2 + C_\xi \sum_{k=0}^\infty\|\mathbf{P} \nabla_H \pi_k\|_H^2 \\
    &\leq \xi (\|v_n^*\|_V^2 +1) + C_\xi \sum_{k=0}^\infty\|\mathbf{P} \nabla_H \pi_k\|_H^2.
\end{align*}
The same notation $C_\xi$ is used to refer to different constants -- in the same fashion as the notation $C$ being kept from one line to the other -- and it emphasises the $\xi$ dependence of $C_\xi$. Furthermore, the regularising properties of $K$ lead to
\begin{equation}
     \frac{1}{2} \sum_{k=0}^\infty \|\mathbf{P} \nabla_H \pi_k\|_H^2 \leq C(1+ \|U_n^*\|_H^2), \quad  \frac{1}{2} \sum_{k=0}^\infty \|\mathbf{P} \nabla_H \pi_k\|_V^2 \leq C(1+ \|U_n^*\|_V^2).
\end{equation}
Consequently,
\begin{equation}
    \frac{1}{2} \Big\|G_\sigma^n(U_n^*)- \mathcal{M}^{\sigma,U_n^*}\Big\|_{\mathcal{L}_2(\mathcal{W},H)}^2 \leq \frac{1}{2} (a_3 \nabla_3 v_n^*, \nabla_3 v_n^*)_{L^2} + C_\xi (1+ \|v_n^*\|_H^2) + \xi \|v_n^*\|_V^2.
\end{equation}
Therefore,
\begin{align}
    \frac{1}{2} \|G_\sigma^n(U_n^*) - \mathcal{M}^{\sigma,U_n^*}\|_{\mathcal{L}_2(\mathcal{W},H)}^2 \|U_n^*\|_H^{p-2} \leq \frac{1}{2} \bigg(a_3 &\nabla_3 \begin{pmatrix}
        v_n^*\\w(v_n^*)
    \end{pmatrix}, \nabla_3 \begin{pmatrix}
        v_n^*\\w(v_n^*)
    \end{pmatrix} \bigg)_{L^2}\|U_n^*\|_H^{p-2} \nonumber\\
    &+ C_\xi(1+\|U_n^*\|_H^p) + \xi \|U_n^*\|_V^2\|U_n^*\|_H^{p-2}.
\end{align}
Hence, gathering all the previous estimates into equation \eqref{Itô-énergie}, we reach
\begin{align}
    d_t \|U_n^*\|_H^p + &c \|U_n^*\|_H^{p-2} \|U_n^*\|_V^2 dt + \alpha \|U_n^*\|_H^{p-2} \|w(v_n^*)\|_V^2 \label{Itô-énergie-diff} \\
    &\leq C \|U_n^*\|_H^{p-2} (U_n^*, G_\sigma^n(U_n^*)dW_t)_H  + C_\xi (\|U_n^*\|_H^p+1) dt + \xi  \|U_n^*\|_H^{p-2} \|U_n^*\|_V^2 dt.  \nonumber
\end{align}

\begin{remark}\label{energy-estimate-martingale-P-hat-K}
    Consider now the problem $(\hat{\mathcal{P}}_K)$, and use the notations of Remarks \ref{F-hat-sigma-Lipschitz} and \ref{G-hat-sigma-Lipschitz}. The martingale contribution to the energy reads
    \begin{equation}
        \frac{p(p-2)}{2} \Big\|\Big[\hat{G}_\sigma^n(U_n^*)- \mathcal{M}^{\sigma,U_n^*}\Big]^* U_n^*\Big\|_H^2 \|U_n^*\|_H^{p-4} + \frac{p}{2} \Big\|\hat{G}_\sigma^n(U_n^*) - \mathcal{M}^{\sigma,U_n^*} \Big\|_{\mathcal{L}_2(\mathcal{W},H)}^2 \|U_n^*\|_H^{p-2}.
    \end{equation}
    The first term is estimated using the same arguments as for the problem $(\mathcal{P}_K)$. Also, it is straightforward to show that
    \begin{multline*}
        \frac{1}{2} \|\hat{G}_\sigma^n(U_n^*) - \mathcal{M}^{\sigma,U_n^*}\|_{\mathcal{L}_2(\mathcal{W},H)}^2 \|U_n^*\|_H^{p-2} \leq \frac{1}{2} \bigg(a_{BT} \nabla_3 \begin{pmatrix}
            v_n^*\\w(v_n^*)
        \end{pmatrix}, \nabla_3 \begin{pmatrix}
            v_n^*\\w(v_n^*)
        \end{pmatrix} \bigg)_{L^2}\|U_n^*\|_H^{p-2} \nonumber\\
        + \frac{1}{2} \bigg(\overline{a_{BC}} \nabla_3 \bar{v}_n^*, \nabla_3 \bar{v}_n^* \bigg)_{L^2}\|U_n^*\|_H^{p-2}+ C_\xi(1+\|U_n^*\|_H^p) + \xi \|U_n^*\|_V^2\|U_n^*\|_H^{p-2}.
    \end{multline*}
    Consequently, the covariation term of the transport noise $(\tilde{\phi}_k \bcdot \nabla_3) \bar{v}_n^*$ is cancelled out by the stochastic diffusion term $\frac{1}{2} (\overline{a_{BC}} \nabla_3 \bar{v}_n^*, \nabla_3 \bar{v}_n^* )_H\|U_n^*\|_H^{p-2}$. We find the following relation
    \begin{align}
        d_t \|U_n^*\|_H^p + &c \|U_n^*\|_H^{p-2} \|U_n^*\|_V^2 dt + \alpha \|U_n^*\|_H^{p-2} \|w(v_n^*)\|_V^2 \label{Itô-énergie-diff-P-hat-K} \\
        &\leq C \|U_n^*\|_H^{p-2} (U_n^*, \hat{G}_\sigma^n(U_n^*)dW_t)_H  + C_\xi (\|U_n^*\|_H^p+1) dt + \xi  \|U_n^*\|_H^{p-2} \|U_n^*\|_V^2 dt.  \nonumber
    \end{align}
\end{remark}
By time integration and taking the expectation, we have
\begin{align}
    \E[\|U_n^*(t)\|_H^p] + \E[ \int_0^t (c- \xi) \|U_n^*\|_H^{p-2} \|U_n^*\|_V^2 dr] &\leq \E[\|U_n^*(0)\|_H^p] + C_\xi \E[1+\int_0^t \|U_n^*\|_H^p dr]. \label{Itô-énergie-espérance}
\end{align}
Now use the previous inequality \eqref{Itô-énergie-espérance} with $\xi = \frac{c}{2}$, and apply the (classical) Grönwall lemma, so that for $t \in [0,T]$,
\begin{align}
    \E[\|U_n^*(t)\|_H^p] \leq C. \label{majoration-Gronwall}
\end{align}
\emph{Step 3:} We are now in position to deduce the energy estimates we claimed earlier. Using the inequality \eqref{Itô-énergie-espérance},
\begin{align}
    \E\Bigg[\int_0^t \|U_n^*(r)\|_H^{p-2} \|U_n^*(r)\|_V^2 dr\Bigg] \leq C \Bigg(1 + \int_0^{T} \E\Big[\|U_n^*\|_H^p\Big] dr \Bigg),
\end{align}
and applying \eqref{majoration-Gronwall} allows to conclude that
\begin{align}
    \E[ \int_0^{T} \|U_n^*(r)\|_H^{p-2} \|U_n^*(r)\|_V^2 dr] \leq C.
\end{align}
Moreover, by the inequality \eqref{Itô-énergie-diff} with $\xi = \frac{c}{2}$, upon integrating and taking the supremum and the expectation, one also has
\begin{multline}
    \E[\sup_{0\leq t\leq T} \|U_n^*\|_H^{p}] \leq \E[\|U_n^*(0)\|_H^p] + C \int_0^{T} \E[\|U_n^*\|_H^p ] dr + C (T+1) + \frac{c}{2} \E[ \int_0^t \|U_n^*(t)\|_H^{p-2} \|U_n^*\|_V^2 dr]\\
    + \E\Bigg[\sup_{0\leq t\leq T} \bigg|\int_0^t \|U_n^*\|_H^{p-2}(U_n^*, G_\sigma(U_n^*)dW_r - \frac{1}{\rho_0} \mathbf{P} \nabla_H dp_r^\sigma)_H \bigg|\Bigg].
\end{multline}
We estimate the martingale term using the Burkholder-Davis-Gundy inequality (see Theorem 4.36 in \cite{book_DPZ_2014}) and Young's inequality. This yield,
\begin{multline} \label{eq-sup-noise}
    \E\Bigg[\sup_{0\leq t\leq T} \bigg|\int_0^t \|U_n^*\|_H^{p-2}(U_n^*, G_\sigma^n(U_n^*) - \frac{1}{\rho_0} \mathbf{P} \nabla_H p_t^\sigma)_H \bigg|\Bigg] \\
    \leq C \E\Big[\Big(\int_0^{T} \|U_n^*\|_H^{2(p-2)} \underbrace{\|[G_\sigma^n(U_n^*)- \mathcal{M}^{\sigma,U_n^*}]^*U_n^*\|_{L^2}^2}_{C (1+ \|U_n^*\|_H^4)} dr\Big)^{1/2}\Big] \leq C \E\Big[\Big(1+\int_0^{T} \|U_n^*\|_H^{2p} dr \Big)^{1/2}\Big] \\
    \leq C \E\Big[1+\sup_{0\leq t\leq T} \|U_n^*\|_H^{p/2} (\int_0^{T} \|U_n^*\|_H^{p})^{1/2} \Big]\\
    \leq \frac{1}{2} \E\Big[ \sup_{0\leq t\leq T} \|U_n^*\|_H^{p} \Big] + C\E\Big[1 + (\int_0^{T} \|U_n^*\|_H^{p}) \Big]
\end{multline}
Therefore,
\begin{align}
    \E[\sup_{0\leq t\leq T} \|U_n^*\|_H^{p}] \leq C \int_0^{T} \E[1+\|U_n^*\|_H^p ] dr,
\end{align}
and by equation \eqref{majoration-Gronwall}, we reach eventually
\begin{align}
    \E[\sup_{0\leq t\leq T} \|U_n^*\|_H^{p}] \leq C.
\end{align}
\CQFD

\begin{remark} \hfill \label{remark-sup-noise}
\begin{enumerate}
    \item If the noise is not divergence-free, then equation \eqref{eq-sup-noise} rather reads,
    \begin{multline}
        \E\Bigg[\sup_{0\leq t\leq T} \bigg|\int_0^t \|U_n^*\|_H^{p-2} (U_n^*, G_\sigma^n(U_n^*) - \frac{1}{\rho_0} \mathbf{P} \nabla_H p_t^\sigma)_H \bigg|\Bigg]\\
        \leq C \E\Big[\Big(1+\int_0^{T} \|U_n^*\|_H^{2p} + \|U_n^*\|_H^{2p-2} \|U_n^*\|_V^2 dr \Big)^{1/2}\Big] \\
        \leq C \E\Big[1+\sup_{0\leq t\leq T} \|U_n^*\|_H^{p/2} \Big(1+\int_0^{T} \|U_n^*\|_H^{p-2} \|U_n^*\|_V^2 dr \Big)^{1/2} \Big] \\
        \leq \frac{1}{2} \E\Big[\sup_{0\leq t\leq T} \|U_n^*\|_H^p \Big] + C\Big(1+ \E\Big[\int_0^{T} \|U_n^*\|_H^{p-2} \|U_n^*\|_V^2 dr \Big] \Big).
    \end{multline}
    Therefore the argument still holds for $p=2$, and so does it for $p \geq 3$ if the additional parabolicity condition \eqref{eq-parabolicity-condition} holds.
    \item Adapting the argument to the problem $(\hat{\mathcal{P}}_K)$ is straightforward, since equation \eqref{basic-energ-cross-term} gives the following bound,
    \begin{equation}
        \|U_n^*\|_H^{p-4} \Big\|\Big[G_\sigma^n(U_n^*)- \mathcal{M}^{\sigma,U_n^*}\Big]^*U_n^*\Big\|_H^2 \leq C(\|U_n^*\|_H^{p}+1).
    \end{equation}
\end{enumerate}

\noindent
Hence, Lemma \ref{lemma-basic-estimate} holds for both $(\mathcal{P}_K)$ and $(\hat{\mathcal{P}}_K)$, with divergence-free noises. Assuming that the ($p$-dependent) parabolicity condition holds, it also holds for divergent noises.

\end{remark}

\subsection{Tightness of the laws} \label{subsec-tension}

Now we show that the solutions $U_n^*$ have tight laws in the space $L^2([0,T],H) \cap C([0,T],\mathcal{D}(A^{-3}))$. We only consider the ``full'' problem $(\mathcal{P}_K)$ in this section, since adapting the arguments to $(\hat{\mathcal{P}}_K)$ is immediate.

\subsubsection{Tightness in $L^2([0,T],H)$ } \label{subsec-tension1}

Let $\alpha \in (1/3,1/2)$. In the following, if $x_1, x_2 \in \R$, we may use the notation $x_1 \vee x_2 := \max(x_1, x_2)$. By Theorem 2.1 of \cite{FG_1995} (see also \cite{Lions_1969} and \cite{Temam_1995}), the embedding
\begin{align}
    L^2\big([0,T],V\big) \cap W^{\alpha,2}\big([0,T],\mathcal{D}(A^{-2})\big) \hookrightarrow L^2\big([0,T],H\big) \text{ is compact.} \label{compact-embedding}
\end{align}
So it is enough to show that $(U_n^*)$ is bounded in probability in both spaces $L^2([0,T],V)$ and $W^{\alpha,2}([0,T],\mathcal{D}(A^{-2}))$.

It can be readily shown that the processes $U_n^*$ are bounded in probability in $L^2([0,T],V)$, using Markov inequality and the previous energy estimates. Hence we focus on showing that the processes $(U_n^*)$ are bounded in probability in $W^{\alpha,2}\Big([0,T],\mathcal{D}(A^{-2})\Big)$, with the inclusion relation
\begin{equation}
    W^{1,2}\Big([0,T],\mathcal{D}(A^{-2})\Big) \subset W^{1/2,2}\Big([0,T],\mathcal{D}(A^{-2})\Big) \subset W^{\alpha,2}\Big([0,T],\mathcal{D}(A^{-2})\Big).
\end{equation}
Remind the Galerkin approximation problem $(\mathcal{P}_n)$, 
\begin{multline}
    U_n^*  = P^n(U_0^*) - \int_0^t [A U_n^* + B^n (U_n^*) + \Gamma^n U_n^* + \frac{1}{\rho_0} (\nabla_H p)^n dt + F_\sigma^n (U_n^*)] dr \\
    + \int_0^t [G_\sigma^n(U_n^*)dW_r - \frac{1}{\rho_0} (\mathbf{P} \nabla_H dp_r^{\sigma})^n],
\end{multline}
which we rewrite as $U_n^* := \sum_{k=1}^8 J_n^k$, with
\begin{align*}
    J_n^1 := P^n(U_0^*), && J_n^2:= - \int_0^t A U_n^* dr, && J_n^3 := - \int_0^t B^n (U_n^*) dr,\\
    J_n^4 := - \int_0^t \Gamma^n U_n^* dr, && J_n^5 := -\int_0^t \frac{1}{\rho_0} (\mathbf{P} \nabla_H p)^n dr, && J_n^6 := -\int_0^t F_\sigma^n (U_n^*) dr,\\
    J_n^7 := \int_0^t G_\sigma^n(U_n^*) dW_r && J_n^8 := \int_0^t \frac{1}{\rho_0} (\mathbf{P} \nabla_H dp_r^\sigma)^n.
\end{align*}
Now we show that each term is bounded in $L^1\Big(\Omega, W^{\alpha,2}\Big([0,T],\mathcal{D}(A^{-2})\Big)\Big)$, which allows to conclude. Thanks to the energy estimates \eqref{estimation-énerg} and the properties proven in subsection \ref{subsec-prelim}, the following holds,
\begin{equation}
    \E\Big[\sum_{i=1}^5 \|J_n^i\|_{W^{1,2}([0,T],\mathcal{D}(A^{-\gamma}))}\Big] \leq C.
\end{equation}
For $J_n^6$, by continuity of $F_\sigma : H \rightarrow \mathcal{D}(A^{-\gamma})$,
\begin{multline}
    \E\Big[\|J_n^6\|_{W^{1,2}([0,T],\mathcal{D}(A^{-\gamma}))}^2\Big] = \E\Big[\|\int_0^t F_\sigma^n(U_n^*) dr \|_{W^{1,2}([0,T],\mathcal{D}(A^{-\gamma}))}^2\Big] \\ \leq C \E\Big[\|F_\sigma^n(U_n^*)\|_{L^2([0,T],\mathcal{D}(A^{-\gamma})}^2\Big]  \leq C \E\Big[\int_0^{T} \|F_\sigma(U_n^*)\|_{\mathcal{D}(A^{-\gamma})}^2 dt\Big] \\ \leq  C \E\Big[\int_0^{T} \|U_n^*\|_H^2 dt\Big] \leq C.
\end{multline}
For $J_n^7$, by the Burkholder-Davis-Gundy inequality (see Lemma 2.1 in \cite{FG_1995}), since $\alpha <1/2$,
\begin{multline}
    \E\Big[\|J_n^7\|_{W^{\alpha,2}([0,T],\mathcal{D}(A^{-\gamma}))}^2\Big] \leq C_\alpha \E\Big[\int_0^{T} \|G_\sigma U_n^*\|_{\mathcal{L}_2(L^2(\mathcal{S},\mathbbm{R^3}),\mathcal{D}(A^{-\gamma}))}^2 dt\Big] \\ \leq C_\alpha \E\Big[1+\int_0^{T} \|U_n^*\|_H^2 dr\Big] \leq C,
\end{multline}
by continuity of $G_\sigma \Psi : H \rightarrow \mathcal{D}(A^{-\gamma})$ and using \eqref{estimation-énerg}. Furthermore, for $J_n^8$, since the operator $\frac{1}{\rho_0} \nabla_H p_t^\sigma : H \rightarrow \mathcal{L}_2(\mathcal{W}, H)$ is continuous,
\begin{multline}
    \E\Big[\|J_n^8\|_{W^{\alpha,2}([0,T],\mathcal{D}(A^{-2}))}^2\Big] \leq C_\alpha \E\Big[\int_0^{T} \|\frac{1}{\rho_0} \mathbf{P} \nabla_H p^\sigma\|_{\mathcal{L}_2(L^2(\mathcal{S},\mathbbm{R^3}),\mathcal{D}(A^{-2}))}^2 dt\Big] \\ \leq C_\alpha \E\Big[1+\int_0^{T} \|U_n^*\|_H^2 dr\Big] \leq C.
\end{multline}
Gathering all the previous estimates, we get
\begin{align}
    \E\Big[\|U_n^*\|_{W^{\alpha,2}([0,T],\mathcal{D}(A^{-2}))}\Big] \leq C,
\end{align}
so that $U_n^*$ is bounded in probability in $W^{\alpha,2}([0,T],\mathcal{D}(A^{-2}))$ by Markov inequality. Therefore we deduce that $U_n^*$ is tight in $L^2([0,T],H)$ using that the embedding \eqref{compact-embedding} is compact.

\CQFD

\subsubsection{Tightness in $C([0,T],\mathcal{D}(A^{-3}))$}

Since $\alpha > 1/3$, using Theorem 2.2 of \cite{FG_1995} (see also \cite{Lions_1969} and \cite{Temam_1995}), the embeddings
\begin{align}
    W^{1,2}([0,T],\mathcal{D}(A^{-2})) &\hookrightarrow C([0,T],\mathcal{D}(A^{-3})) \label{compact-embedding-1}
\end{align}
and
\begin{align}
    W^{\alpha,3}([0,T],\mathcal{D}(A^{-2})) &\hookrightarrow C([0,T],\mathcal{D}(A^{-3}))\label{compact-embedding-2}
\end{align}
are compact. Using the previous subsection \ref{subsec-tension1}, we directly have
\begin{align}
    \Big\|U_n^* - \int_0^t \Big[ G(U_n^*)dW_r + \frac{1}{\rho_0} \mathbf{P} \nabla_H dp_r^\sigma  \Big] \Big\|_{ W^{1,2}([0,T],\mathcal{D}(A^{-2}))} \leq C.
\end{align}
Moreover,
\begin{multline}
    \E\Bigg[\Big\|\int_0^t \Big[ G_\sigma(U_n^*)dW_r + \frac{1}{\rho_0} \mathbf{P} \nabla_H dp_r^\sigma \Big]\Big\|_{W^{\alpha,3}([0,T],\mathcal{D}(A^{-2}))}^3 \Bigg]\\
    \leq C \E\Bigg[ \int_0^{T}  \Big\|\Big[ G_\sigma(U_n^*) + \frac{1}{\rho_0} \mathbf{P} \nabla_H dp^\sigma \Big]\Big\|_{\mathcal{L}_2(\mathcal{W},\mathcal{D}(A^{-2}))}^3 dr \Bigg] \leq C \Bigg(\E\Bigg[\int_0^{T} \|U_n^*\|_H^3\Bigg] + 1 \Bigg) \leq C,
\end{multline}
using the estimate $\eqref{estimation-énerg}$.

\bigskip

Combining these two energy estimates, we conclude that $U_n^*$ is tight in $C([0,T],\mathcal{D}(A^{-3}))$, using that the embeddings \eqref{compact-embedding-1} and \eqref{compact-embedding-2} are compact.

\CQFD

\subsection{Taking the limit of the approximate solutions} \label{subsec-limite}

Now we take the limit of the Galerkin equation and show that limit of $v_n^*$ is a solution of $(\mathcal{P}_K)$. Again, adapting the arguments to $(\hat{\mathcal{P}}_K)$ is immediate.

By Prohorov theorem (see Theorem 2.3 of \cite{book_DPZ_2014}), the probability measures of $(U_n^*)_n$ are relatively compact since they are tight. In addition, defining $W_n =W$ -- where $W$ is the Wiener process of the stochastic basis $(\Omega, \mathcal{F}, (\mathcal{F}_t)_t, \mathbbm{P}, W)$ -- we deduce that $(U_n^*,W)$ is tight in the product space $[L^2([0,T],H) \cap C([0,T],\mathcal{D}(A^{-3}))] \times C^{1/4}([0,T], H^{-2}(\mathcal{S}, \R^3))$. This means there exists a subsequence $(U_{\phi(n)}^*, W_{\phi(n)})_n$ of $(U_n^*,W_n)_n$ such that the associated probability measures converge in the space of probability measures. Consequently, the ordered pair $(U_{\phi(n)}^*, W_{\phi(n)})_n$ converges in law in the product space $[L^\infty([0,T], H) \cap L^2([0,T],V)] \times C^{1/4}([0,T], H^{-2}(\mathcal{S}, \R^3))$. By Skorokhod theorem (see Theorem 2.4 of \cite{book_DPZ_2014}), there exists a probability space $(\overline{\Omega}, \overline{\mathcal{F}}, \overline{\mathbbm{P}} )$ equipped with cylindrical Wiener processes $(\overline{W}^n)_n$ and $\overline{W}$, a sequence $(\overline{U}_n)_n$ and a process $\overline{U}$ on this space such that
\begin{align}
    \overline{W}^n \overset{a.s}{\longrightarrow} \overline{W}, \quad \overline{W}^n &\overset{\mathcal{L}}{=} W_n = W, \quad \text{in $C^{1/4}([0,T], H^{-2}(\mathcal{S}, \R^3))$},\\
    \overline{U}_n \overset{a.s}{\longrightarrow} \overline{U}, \quad \overline{U}_n &\overset{\mathcal{L}}{=} U_{\phi(n)}^*, \quad \text{in $L^\infty([0,T], H) \cap L^2([0,T],V)$}. \nonumber
\end{align}
Then, denote by $\overline{\mathcal{F}}^n_t = \sigma((\overline{U}_n(s))_{0 \leq s \leq t}, (\overline{W}^n(s))_{0 \leq s \leq t})$ the $\sigma$-algebra generated by $\overline{U}_n$ and $\overline{W}^n$ on $[0,T]$. Similarly, define $\overline{\mathcal{F}}_t = \sigma((\overline{U}(s))_{0 \leq s \leq t}, (\overline{W}(s))_{0 \leq s \leq t})$. Thus, $(\overline{\Omega}, \overline{\mathcal{F}}, \overline{\mathbbm{P}}, (\overline{\mathcal{F}}^n_t)_t, \overline{W}^n )$ is a stochastic basis for all $n$, so that $\overline{U}_n$ and $\overline{W}^n$ are adapted to $(\overline{\mathcal{F}}^n_t)_t$. Likewise, $(\overline{\Omega}, \overline{\mathcal{F}}, \overline{\mathbbm{P}}, (\overline{\mathcal{F}}_t)_t, \overline{W} )$ is a stochastic basis such that $\overline{U}$ and $\overline{W}$ are adapted to $(\overline{\mathcal{F}}_t)_t$.

The equality $\overline{U}_n \overset{\mathcal{L}}{=} U_{\phi(n)}^*$ implies, in particular, that $\overline{U}_n$ is a solution of some projected Galerkin problem $(\mathcal{P}_{\phi(n)})$ -- see \cite{Bensoussan1995} (Section 4.3.4). Therefore the previously established energy estimates also hold for $\overline{U}_n$. Now, we aim to show that the limit $\overline{U}$ is a solution of $(\mathcal{P}_K)$ -- equation \eqref{problem-PK}. The solution $\overline{U}_n$ of $(\mathcal{P}_K)$ satisfies the following equation,
\begin{multline}
    \overline{U}_n = P^n(\overline{U}_0) - \int_0^t [A \overline{U}_n + B^n (\overline{U}_n) + \Gamma^n \overline{U}_n + \Big(\frac{1}{\rho_0} \mathbf{P} \nabla_H p^{\overline{U}_n}\Big)^n + F_\sigma^n (\overline{U}_n)] dr\\ + \int_0^t \Big[G_\sigma^n(\overline{U}_n) d\overline{W}^n_r - (\frac{1}{\rho_0} \mathbf{P} \nabla_H dp_r^{\sigma,\overline{U}_n})^n\Big].
\end{multline}
We deduce, for $z \in \mathcal{D}(A^3)$,
\begin{multline}
    (\overline{U}_n - P^n(\overline{U}_n(0)),z)_H + \int_0^t (A \overline{U}_n,z)_H dr + \int_0^t (B^n (\overline{U}_n),z)_H dr + \int_0^t (\Gamma^n \overline{U}_n,z)_H dr \\
    + \int_0^t \Big(\Big(\frac{1}{\rho_0} \mathbf{P} \nabla_H p^{\overline{U}_n}\Big)^n, z\Big)_H dr + \int_0^t (F_\sigma^n (\overline{U}_n),z)_H dr = \int_0^t \Big(G_\sigma^n (\overline{U}_n) d\overline{W}^n_r, z\Big)_H \\
    - \int_0^t \Big(\frac{1}{\rho_0} \mathbf{P} \nabla_H dp_r^{\sigma,\overline{U}_n}\Big)^n, z\Big)_H,
\end{multline}
which we rewrite
\begin{equation}
    J_n^1 + J_n^2 + J_n^3 + J_n^4 + J_n^5 + J_n^6 = J_n^7 + J_n^8,
\end{equation}
to show that each term converges. Regarding $J_1$, $J_2$, $J_4$, $J_5$ and $J_6$, it is immediate that, $\overline{\mathbbm{P}}$-almost surely,
\begin{gather}
    J_n^1 \underset{n\longrightarrow \infty}{\longrightarrow} (\overline{U} - \overline{U}(0),z)_H, \quad
    J_n^2 \underset{n\longrightarrow \infty}{\longrightarrow} \int_0^t (A \overline{U},z)_H dr, \quad
    J_n^4 \underset{n\longrightarrow \infty}{\longrightarrow} \int_0^t (\Gamma \overline{U},z)_H dr,  \nonumber\\
    J_n^5 \underset{n\longrightarrow \infty}{\longrightarrow} \int_0^t (\frac{1}{\rho_0} \mathbf{P} \nabla_H p^{\overline{U}},z)_H dr, \quad J_n^6 \underset{n\longrightarrow \infty}{\longrightarrow} \int_0^t (F_\sigma (\overline{U}) ,z)_H dr. 
\end{gather}
For $J_n^3$,
\begin{align}
    |J_n^3 - \int_0^t (B(\overline{U}),z)_H dr| &\leq \int_0^t |(B^n(\overline{U}_n)-B(\overline{U}),z)_H| dr\\
    &\leq \underbrace{\int_0^t |(B^n(\overline{U}_n)-B(\overline{U}_n),z)_H| dr}_{J_n^{3,1}} + \underbrace{\int_0^t |(B(\overline{U}_n)-B(\overline{U}),z)_H| dr}_{J_n^{3,2}}. \nonumber
\end{align}
Remark that 
\begin{align}
    J_n^{3,1} \leq \int_0^{T} |(B(\overline{U}_n),P^n(z)-z)_H| dr &\leq C \int_0^{T}\|B(\overline{U}_n)\|_{V'} \|P^n(z)-z\|_{\mathcal{D}(A^3)}dr \nonumber \\
    &\leq C \int_0^{T}\|\overline{U}_n\|_V^2  \|P^n(z)-z\|_{\mathcal{D}(A^{3})} dr \nonumber\\
    &\leq C \|P^n(z)-z\|_{\mathcal{D}(A^3)}, 
\end{align}
using the energy estimate \eqref{estimation-énerg}. Therefore $\E[J_n^{3,1}] \underset{n\longrightarrow \infty}{\longrightarrow} 0.$ Moreover, since $B(\overline{U}_n)-B(\overline{U}) = B(\overline{U}_n-\overline{U},\overline{U}_n) - B(\overline{U},\overline{U}_n-\overline{U})$, use Cauchy-Schwarz inequality to obtain
\begin{align}
    J_n^{3,2} &\leq C \int_0^t \|\overline{U}_n\|_H \|\overline{U}_n-\overline{U}\|_H \|z\|_{\mathcal{D}(A^3)} dr + C \int_0^t \|\overline{U}\|_H \|\overline{U}_n-\overline{U}\|_H \|z\|_{\mathcal{D}(A^3)} dr \\
    &\leq C [ \|\overline{U}\|_{L^2([0,T],H)}+\sup_{n \geq 0} \|\overline{U}_n\|_{L^2([0,T],H)}] \|\overline{U}_n-\overline{U}\|_{L^2([0,T],H)} \|z\|_{\mathcal{D}(A^3)}, \nonumber
\end{align}
and thus $\E[J_n^{3,2}] \underset{n\longrightarrow \infty}{\longrightarrow} 0$. Gathering the estimates on $J_n^{3,1}$ and $J_n^{3,2}$, $\overline{\mathbbm{P}}$-a.s yields,
$$J_n^3 \underset{n\longrightarrow \infty}{\longrightarrow} \int_0^t (B(\overline{U}),z)_H dr.$$
We show now that $J_n^7$ and $J_n^8$ converge in probability. For $J_n^7$, we write, by the triangle inequality,
\begin{multline*}
    \|G_\sigma^n(\overline{U}_n) - G_\sigma(\overline{U})\|_{\mathcal{L}_2(\mathcal{W},\mathcal{D}(A^{-3}))}^2 \leq 2 \underbrace{\|G_\sigma^n(\overline{U}_n) - G_\sigma^n(\overline{U})\|_{\mathcal{L}_2(\mathcal{W},\mathcal{D}(A^{-3}))}^2}_{\leq C \|\overline{U}-\overline{U}_n\|_{H}^2} \\
    + 2\|G_\sigma^n(\overline{U}) - G_\sigma(\overline{U})\|_{\mathcal{L}_2(\mathcal{W},\mathcal{D}(A^{-3}))}^2,
\end{multline*}
and remark that
\begin{align*}
    \|G_\sigma^n(\overline{U}) - G_\sigma(\overline{U})\|_{\mathcal{L}_2(\mathcal{W},\mathcal{D}(A^{-3}))}^2 = \|(P^n-I) G_\sigma(\overline{U})\|_{\mathcal{L}_2(\mathcal{W},H)}^2 &\leq C \lambda_{n+1}^{-3} \|G_\sigma(\overline{U})\|_{\mathcal{L}_2(\mathcal{W},\mathcal{D}(A^{-3}))}^2 \\
    &\leq C \lambda_{n+1}^{-3} (1+\|\overline{U}\|_H^2).
\end{align*}
We deduce
\begin{align*}
    \E\Big[\int_0^{T} \|G_\sigma^n(\overline{U}) - G_\sigma\overline{U})\|_{\mathcal{L}_2(\mathcal{W},\mathcal{D}(A^{-3})}^2 dr\Big] \leq C \lambda_{n+1}^{-3} (1+\E[\|\overline{U}\|_H^2]) &\leq C \lambda_{n+1}^{-3} (1+\E[\sup_{[0,T]}\|\overline{U}\|_H^2])\\
    &\leq C \lambda_{n+1}^{-3} \underset{n\longrightarrow \infty}{\longrightarrow}0.
\end{align*}
Therefore, $G_\sigma^n(\overline{U}_n) \underset{n\longrightarrow \infty}{\longrightarrow} G_\sigma(\overline{U})$, in probability. By Lemma 2.1 of \cite{DGHT_2011}, we reach
\begin{align}
    \int_0^t G_\sigma^n(\overline{U}_n) d\overline{W}_r^n \underset{n\longrightarrow \infty}{\longrightarrow} \int_0^t G_\sigma(\overline{U}) d\overline{W}_r \text{ in probability in $L^2([0,T],\mathcal{L}_2(L^1(\mathcal{S},\R^3),\mathcal{D}(A^{-3}))$}.
\end{align}
Using similar arguments, we also have
\begin{align}
    \int_0^t \Big(\frac{1}{\rho_0} \mathbf{P} \nabla_H dp_r^{\sigma,{\overline{U}_n}}\Big)^n \underset{n\longrightarrow \infty}{\longrightarrow} \int_0^t \frac{1}{\rho_0} \mathbf{P} \nabla_H dp_r^{\sigma,{\overline{U}}} \text{ in probability in $L^2([0,T],\mathcal{L}_2(L^1(\mathcal{S},\R^3),\mathcal{D}(A^{-3}))$}.
\end{align}
Thus,
\begin{align}
    J_n^7 \overset{\mathbbm{P}}{\underset{n\longrightarrow \infty}{\longrightarrow}} \int_0^t (G_\sigma(\overline{U}) d\overline{W}_r,z)_H, \text{ and } \quad J_n^8 \overset{\mathbbm{P}}{\underset{n\longrightarrow \infty}{\longrightarrow}} \int_0^t (\frac{1}{\rho_0} \mathbf{P} \nabla_H dp_r^{\sigma,{\overline{U}}},z)_H.
\end{align}
Consequently, there exists a subsequence of the family $\overline{U}_n$ such that $J_n^7$ and $J_n^8$ converge almost surely to $\int_0^t (G_\sigma(\overline{U}) d\overline{W}_r,z)_H$ and $\int_0^t (\frac{1}{\rho_0} \mathbf{P} \nabla_H dp_r^{\sigma,{\overline{U}}},z)_H$, respectively.

\bigskip

Taking the limit in the Galerkin equation for such a subsequence, for $z \in \mathcal{D}(A^3)$, the following holds,
\begin{multline}
    \Big(\overline{U} - \overline{U}(0) + \int_0^t [A \overline{U} + B (\overline{U}) + \Gamma\overline{U} + \frac{1}{\rho_0} \mathbf{P} \nabla_H p^{\overline{U}} + F_\sigma(\overline{U})] dr, z\Big)_H \\
    = \int_0^t \Big(G_\sigma (\overline{U}) d\overline{W}_r - \frac{1}{\rho_0} \mathbf{P} \nabla_H dp_r^{\sigma,\overline{U}}, z\Big)_H.
\end{multline}
Hence, the following equality holds true almost surely in $\mathcal{D}(A^{-3})$,
\begin{align}
    d_t\overline{U} + [A \overline{U} + B (\overline{U}) +  \Gamma\overline{U} + \frac{1}{\rho_0} \mathbf{P} \nabla_H p^{\overline{U}} + F_\sigma(\overline{U})]dt= G_\sigma (\overline{U}) d\overline{W}_t - \frac{1}{\rho_0} \mathbf{P} \nabla_H dp_t^{\sigma,\overline{U}}.
\end{align}
Furthermore, as proven above, $\overline{U}_n$ fulfils the energy estimates \eqref{estimation-énerg}. Therefore, by taking the limit as $n \rightarrow 0$, we deduce
\begin{align}
    \overline{U} \in L^2(\overline{\Omega},L^2([0,T],V)) \cap L^2(\overline{\Omega},L^\infty([0,T],H)).
\end{align}
\CQFD
Consequently, Theorem \ref{theorem-filtered-weak} holds for $(\mathcal{P}_K)$ and $(\hat{\mathcal{P}}_K)$.

\section{Existence and uniqueness of a local-in-time pathwise solution for the low-pass filtered problem} \label{section-loc-wpn}

In this section, we show the second point of Theorem \ref{theorem-filtered-weak}, that is there exists a unique local-in-time pathwise solution to the (low-pass) filtered problem $(\mathcal{P}_K)$ -- equation \eqref{problem-PK}. The arguments can be adapted to the approximated problem $(\hat{\mathcal{P}}_K)$ -- equation \eqref{problem-approx-PK}, see Remark \ref{remark-PK-approx-strong} in particular. In the following, using the method presented in \cite{BS_2021}, we consider the following cut-off problem $(\mathcal{P}^\kappa)$,
\begin{multline}
    d_t X + [A X + \theta_\kappa(\|X - X^{ref}\|_V) B (X) + \Gamma X + \Big(\frac{1}{\rho_0} \mathbf{P} \nabla_H p^{X}\Big) + F_\sigma (X)]dt \\
    = \theta_\kappa(\|X - X^{ref}\|_V) G_\sigma (X) dW_t - \Big(\frac{1}{\rho_0} \mathbf{P} \nabla_H dp_t^{\sigma,X}\Big) \quad (\mathcal{P}^\kappa),
\end{multline}
with the initial condition $U_0 \in V$.

\subsection{Energy estimate of a cut-off system} \label{subsec-strong-energ}

In this subsection, we consider the projected cut-off problems $(\mathcal{P}_n^\kappa)$,
\begin{multline}
    d_t X_n + [A X_n + \theta_\kappa(\|X_n - X_n^{ref}\|_V) B^n (X_n) + \Gamma^n X_n + \Big(\frac{1}{\rho_0} \mathbf{P} \nabla_H p^{X_n}\Big)^n + F_\sigma^n (X_n)]dt \\
    = G_\sigma^n (X_n) dW_t - \Big(\frac{1}{\rho_0} \mathbf{P} \nabla_H dp_t^{\sigma,X_n}\Big)^n \quad (\mathcal{P}_n^\kappa),
\end{multline}
with the initial condition $U_0 \in V$. The aim of this subsection is to prove that, for all $T>0$, there exists a constant $C_{\kappa, T}$, such that for all $n$,
\begin{align}
    \E[\sup_{0 \leq t\leq T} \|X_n\|_V^2 ] +\E[\int_0^{T} \|X_n\|_{\mathcal{D}(A)}^2 dt] \leq C_{\kappa, T}.
\end{align}
Define $X^{ref}$ the unique \emph{deterministic} solution to the heat equation,
\begin{equation}
    d_t X^{ref} + A X^{ref} = 0, \quad (X^{ref})_{t=0} = U_0,
\end{equation}
which belongs to $L^\infty\big([0,T], V\big) \cap L^2\big([0,T], \mathcal{D}(A)\big)$ by a standard argument. We denote by $X_n^{ref} := P^n X^{ref}$ the projection of $X^{ref}$ onto $H_n = P^n H$. Then, define the problems $(\mathcal{P}_n^\kappa)$ on $H_n$ as follows, with the initial condition $(X_n)_{t=0} = P^n U_0$,
\begin{multline}
    d_t X_n + [A X_n + \theta_\kappa(\|X_n - X_n^{ref}\|_V) B^n (X_n) + \Gamma^n X_n + \Big(\frac{1}{\rho_0} \mathbf{P} \nabla_H p^{X_n}\Big)^n + F_\sigma^n (X_n)]dt \\
    = \theta_\kappa(\|X - X^{ref}\|_V) G_\sigma^n (X_n) dW_t - \Big(\frac{1}{\rho_0} \mathbf{P} \nabla_H dp_t^{\sigma,X_n}\Big)^n \quad (\mathcal{P}_n^\kappa),
\end{multline}
where $\theta_\kappa \in C^\infty(\R)$ is such that
\begin{align}
    \mathbbm{1}_{[-\kappa/2,\kappa/2]} \leq \theta_\kappa \leq \mathbbm{1}_{[-\kappa,\kappa]},
\end{align}
and is referred to as the cut-off coefficient. Also, $p^{X_n}$ and $p_t^{\sigma,X_n}$ are respectively the pressure bounded variation and martingale terms associated to $X_n$. Now we infer that $\bar{X}_n = X_n - X_n^{ref}$ fulfils
\begin{multline}
    d_t \bar{X}_n + \Big[A \bar{X}_n + \theta_\kappa(\|\bar{X}_n\|_V) B^n (\bar{X}_n + X_n^{ref}) + \Gamma^n (\bar{X}_n + X_n^{ref}) + \Big(\frac{1}{\rho_0} \mathbf{P} \nabla_H p^{\bar{X}_n + X_n^{ref}}\Big)^n \\
    + F_\sigma^n (\bar{X}_n + X_n^{ref})\Big]dt = \theta_\kappa(\|\bar{X}_n\|_V) G_\sigma^n (\bar{X}_n + X_n^{ref}) dW_t - \Big(\frac{1}{\rho_0} \mathbf{P} \nabla_H dp_t^{\sigma,\bar{X}_n + X_n^{ref}}\Big)^n, \label{eq-X-bar}
\end{multline}
with $(\bar{X}_n)_{t=0} = 0$. Remind that, for all $v^\sharp \in \mathcal{D}(A)$,  the following holds
\begin{gather}
    A v^\sharp = A^v v^\sharp, \quad \mathcal{A} A^v v^\sharp = A^v \bar{v}^\sharp, \quad \mathcal{R} A^v v^\sharp = A^v \tilde{v}^\sharp
\end{gather}
Let $n$ a positive integer, and $X_n$ a solution of $(\mathcal{P}_n^\kappa)$. Applying It\={o}'s lemma to equation \eqref{eq-X-bar} yields
\begin{align}
    d_t &\|\bar{X}_n\|_V^2 \nonumber \\
    = &- 2 \Big( A \bar{X}_n, A \bar{X}_n + \theta_\kappa(\|\bar{X}_n\|_V) B^n (\bar{X}_n + X_n^{ref}) + \Gamma^n (\bar{X}_n + X_n^{ref}) + \Big(\frac{1}{\rho_0} \mathbf{P} \nabla_H p^{\bar{X}_n + X_n^{ref}}\Big)^n \Big)_H dt \nonumber \\ 
    &- 2 \Big( A \bar{X}_n, F_\sigma^n (\bar{X}_n + X_n^{ref})\Big)_H dt + 2 \Big( A \bar{X}_n, \theta_\kappa(\|\bar{X}_n\|_V) G_\sigma^n (\bar{X}_n + X_n^{ref}) dW_t \nonumber\\
    &- \Big(\frac{1}{\rho_0} \mathbf{P} \nabla_H dp_t^{\sigma, \bar{X}_n + X_n^{ref}}\Big)^n \Big)_H+ \|\theta_\kappa(\|\bar{X}_n\|_V) G_\sigma^n (\bar{X}_n + X_n^{ref}) - \mathcal{M}^{\sigma, \bar{X}_n + X_n^{ref}} \|_{\mathcal{L}_2(\mathcal{W},V)}^2 dt. \label{Itô-énergie-faible}
\end{align}
where $\mathcal{M}^{\sigma, \bar{X}_n + X_n^{ref}}$ is defined similarly as in subsection \ref{subsec-energy-1}. Now we estimate each term in the previous relation. Begin with the advection term: by Young's inequality, for all $\xi>0$ there exists $C_\xi >0$ such that
\begin{multline}
    \theta_\kappa(\|\bar{X}_n\|_V) (A \bar{X}_n, B^n (\bar{X}_n + X_n^{ref}))_H \\
    \leq C \theta_\kappa(\|\bar{X}_n\|_V) \|\bar{X}_n\|_{\mathcal{D}(A)}^{2} \|\bar{X}_n\|_V + C \|\bar{X}_n\|_{\mathcal{D}(A)}(\|\bar{X}_n\|_{\mathcal{D}(A)}^{1/2} \|\bar{X}_n\|_V^{1/2} +1) \\
    \leq C \kappa \|\bar{X}_n\|_{\mathcal{D}(A)}^2 + C_\xi (1+\|\bar{X}_n\|_V^2) + \xi \|\bar{X}_n\|_{\mathcal{D}(A)}^2, \label{eq-bound-cut-off}
\end{multline}
where we used equation \eqref{majoration PTZ}. In addition,
\begin{align}
    (\Gamma^n (\bar{X}_n + X_n^{ref}), A \bar{X}_n)_H \leq C_\xi (1+\|\bar{X}_n\|_H^2) + \xi \|\bar{X}_n\|_{\mathcal{D}(A)}^2,
\end{align}
where the notation $C_\xi$ is kept from line to line, as in subsection \ref{subsec-energy-1}. For the bounded variation pressure term, we denote by $\bar{X}_n = (v_n^{\bar{X}}, T_n^{\bar{X}}, S_n^{\bar{X}})^\tr$, so that
\begin{align}
    \Big|\Big( A \bar{X}_n ,\frac{1}{\rho_0} \mathbf{P} &\nabla_H p^{\bar{X}_n + X_n^{ref}} \Big)_H\Big| \leq \Big|\Big(g \:  \nabla_H \int_z^0 {(\beta_T T_n^{\bar{X}} + \beta_S S_n^{\bar{X}}) dz'} , A v_n^{\bar{X}} \Big)_{L^2}\Big| \nonumber\\
    &+ \Big|\Big(\nabla_H \int_z^0 K* \Big[ v_s \bcdot \nabla_H w(v_n^{\bar{X}}) + w_s \partial_z w(v_n^{\bar{X}}) - \frac{1}{2} \nabla_3 \bcdot (a_3^K \nabla_3 w(v_n^{*,\bar{X}}))\Big] , A v_n^{\bar{X}} \Big)_{L^2}\Big| \nonumber\\
    &+ \Big|\Big( \nabla_H \int_z^0 K* \Big[v_s \bcdot \nabla_H w_s + w_s \partial_z w_s\Big] + \frac{1}{2} \nabla_3 \bcdot (a_3^K \nabla_3 w_s) \Big] dz', A v_n^{\bar{X}} \Big)_{L^2}\Big| + C_n^{ref} \nonumber\\
    &\leq C(1 + \|w(v_n^{\bar{X}})\|_H + \|\bar{X}_n\|_V)\|\bar{X}_n\|_{\mathcal{D}(A)} \leq C_\xi(1 + \|\bar{X}_n\|_V^2) + \xi \|\bar{X}_n\|_{\mathcal{D}(A)}^2,
\end{align}
where $C_n^{ref}$ is a constant depending only on $\|X_n^{ref}\|_{\mathcal{D}(A)}$ and the parameters of the problem. Notice that $C_n^{ref}$ is bounded independently of $n$ since $\|X_n^{ref}\|_{\mathcal{D}(A)}$ is. To evaluate the noise associated drift term, we denote $D_{\mu,\nu} = Diag(\mu,\mu,\nu)$, so that, for $\xi >0$,
\begin{align}
    (F_\sigma^n (\bar{X}_n + X_n^{ref}), A \bar{X}_n )_H &= \Big(d_t v_s + B(v_n^{\bar{X}},v_s) - \frac{1}{2} \nabla_3 \bcdot (a_3 \nabla_3 v_s) + \big[Av_s + \Gamma v_s \big], A v_n^{\bar{X}}\Big)_{L^2} \nonumber\\
    &\qquad \qquad \qquad \qquad \qquad \qquad + C_n^{ref2} + \frac{1}{2} \Big(\nabla_3 \bcdot (a_3 \nabla_3 v_n^{\bar{X}}), A v_n^{\bar{X}} \Big)_{L^2} \nonumber\\
    &\leq C(1+ \| \bar{X}\|_V) \| \bar{X}\|_{\mathcal{D}(A)} + \frac{1}{2}\Big(D_{\mu,\nu} a_3 \nabla_3 (\nabla_3 v_n^{\bar{X}}), \nabla_3 (\nabla_3 v_n^{\bar{X}})\Big)_{L^2}\nonumber\\
    &\leq C_\xi (1+ \| \bar{X}\|_V^2) + \xi \| \bar{X}\|_{\mathcal{D}(A)}^2 + \frac{1}{2}\Big(D_{\mu,\nu} a_3 \nabla_3 (\nabla_3 v_n^{\bar{X}}), \nabla_3 (\nabla_3 v_n^{\bar{X}})\Big)_{L^2}, \label{eq-F-eq-h1h2}
\end{align}
where $C_n^{ref2}$ is a constant that is bounded independently of $n$, and which depends only on $\|X_n^{ref}\|_{\mathcal{D}(A)}$ and the parameters of the problem. Here we used the fact that $A v_n^{\bar{X}} = A^v v_n^{\bar{X}}$, and the boundary conditions \eqref{eq-noiseBC-strong} and \eqref{boundary-conditions}, remarking the $\Gamma_u$ and $\Gamma_b$ are flat surfaces.

To evaluate $\frac{1}{2} \|\theta_\kappa(\|\bar{X}_n\|_V) G_\sigma^n (\bar{X}_n + X_n^{ref}) - \mathcal{M}^{\sigma, \bar{X}_n + X_n^{ref}}\|_{\mathcal{L}_2(\mathcal{W},V)}^2$, we define
\begin{align}
    \psi_k = &-A \phi_k^H - B(\phi_k,v_s) - \Gamma \phi_k^H \in \mathcal{D}(A),
\end{align}
so that,
\begin{align}
    \frac{1}{2} &\|\theta_\kappa(\|\bar{X}_n\|_V) G_\sigma^n (\bar{X}_n + X_n^{ref}) - \mathcal{M}^{\sigma, \bar{X}_n + X_n^{ref}}\|_{\mathcal{L}_2(\mathcal{W},V)}^2 \nonumber\\
    &= \frac{1}{2} \sum_{k=0}^\infty \|\theta_\kappa(\|\bar{X}_n\|_V)(\psi_k - \mathbf{P}[\phi_k \bcdot \nabla_3 v_n^{\bar{X}}]) - \mathbf{P} \nabla_H \pi_k^{\bar{X}}\|_V^2 + C_n^{ref3} \nonumber \\
    &\leq  (\frac{1}{2} + \xi) \sum_{k=0}^\infty \theta_\kappa(\|\bar{X}_n\|_V) \|\psi_k - \mathbf{P}[\phi_k \bcdot \nabla_3 v_n^{\bar{X}}]\|_V^2 + C_\xi \Big(1+\sum_{k=0}^\infty \|\mathbf{P} \nabla_H \pi_k^{\bar{X}}\|_V^2\Big),
\end{align}
with $C_n^{ref3}$ enjoying, again, the same properties as $C_n^{ref1}$ and $C_n^{ref2}$. Consequently,
\begin{equation}
    \frac{1}{2} \sum_{k=0}^\infty \|\psi_k - \mathbf{P}[\phi_k \bcdot \nabla_3 v_n^{\bar{X}}]\|_V^2 \leq \frac{1}{2}\Big(D_{\mu,\nu} a_3 \nabla_3 (\nabla_3 v_n^{\bar{X}}), \nabla_3 (\nabla_3 v_n^{\bar{X}})\Big)_{L^2} + C(1+\|\bar{X}_n\|_V^2),
\end{equation}
and, using the regularisation property of $\mathcal{C}_K$ -- equation \eqref{regularisation-CK} -- we deduce the following,
\begin{align}
    \sum_{k=0}^\infty \|\mathbf{P} \nabla_H \pi_k^{\bar{X}}\|_V^2 \leq C (1 + \|w(v_n^{\bar{X}})\|_H^2) &\leq C(1+ \|v_n^{\bar{X}}\|_{V}^2).
\end{align}
Therefore,
\begin{multline}
    \frac{1}{2} \|\theta_\kappa(\|\bar{X}_n\|_V) G_\sigma^n (\bar{X}_n + X_n^{ref}) - \mathcal{M}^{\sigma,\bar{X}_n + X_n^{ref}}\|_{\mathcal{L}_2(\mathcal{W},V)}^2
    \leq \frac{1}{2}\Big(D_{\mu,\nu} a_3 \nabla_3 (\nabla_3 v_n^{\bar{X}}), \nabla_3 (\nabla_3 v_n^{\bar{X}})\Big)_{L^2}\\
    + C_\xi(1+\|\bar{X}_n\|_V^2) + \xi \|\bar{X}_n\|_{\mathcal{D}(A)}^2.
\end{multline}
In addition, thanks to the Burkholder-Davis-Gundy inequality, for all $T>0$,
\begin{multline}
    \E\Bigg[\sup_{0\leq t\leq T} \bigg|\int_0^t ( A \bar{X}_n, \theta_\kappa(\|\bar{X}_n\|_V) G_\sigma^n (\bar{X}_n + X_n^{ref}) dW_t - \frac{1}{\rho_0} \mathbf{P} \nabla_H dp_t^{\sigma, \bar{X}_n + X_n^{ref}})_H \bigg|\Bigg] \\
    \leq  C \E\Big[\Big(\int_0^{T} \|\bar{X}_n\|_{\mathcal{D}(A)}^2 \underbrace{\|[\theta_\kappa(\|\bar{X}_n\|_V) G_\sigma^n(\bar{X}_n + X_n^{ref})- \mathcal{M}^{\sigma,\bar{X}_n + X_n^{ref}}]^*(\bar{X}_n + X_n^{ref})\|_{L^2}^2}_{\leq C (1+ \theta_\kappa^2(\|\bar{X}_n\|_V) \|\bar{X}_n + X_n^{ref}\|_V^2)} dr\Big)^{1/2}\Big] \\
    \leq C \E\Big[ \Big(\int_0^{T} \|\bar{X}_n\|_{\mathcal{D}(A)}^2 dr\Big)^{1/2} \sup_{0\leq t\leq T} \sqrt{1+ \theta_\kappa^2(\|\bar{X}_n\|_V) \|\bar{X}_n + X_n^{ref}\|_V^2}\Big] \\
    \leq \xi \E\Big[ \int_0^{T} \|\bar{X}_n\|_{\mathcal{D}(A)}^2 dr \Big] +  C_\xi \E\Big[\sup_{0\leq t\leq T} (1+ \theta_\kappa(\|\bar{X}_n\|_V) \|\bar{X}_n\|_V^2)\Big].
\end{multline}
Gathering these estimates, we obtain on the RHS $2(C\kappa + 4 \xi)\|\bar{X}_n\|_{\mathcal{D}(A)}^2$ (see equation \ref{eq-bound-cut-off}). Thus, we choose $\xi = \frac{1}{16}$ and $\kappa = \frac{1}{4C}$, so that, using analogous arguments as in section \eqref{subsec-energy-1}, the following estimates hold for all $T>0$,
\begin{align}
    \E[\sup_{0\leq t \leq T} \|\bar{X}_n\|_V^2 ] &\leq C_{\kappa,T}, \quad \E[\int_0^{T} \|\bar{X}_n\|_{\mathcal{D}(A)}^2 dt] \leq C_{\kappa,T},
\end{align}
where $C_{\kappa,T}$ is a constant independent of $n$.

\begin{remark} \label{remark-PK-approx-strong}
    Regarding the problem $(\hat{\mathcal{P}}_K)$, all the arguments are similar. In particular, the terms $(\hat{F}_\sigma (\bar{X}_n + X_n^{ref}), A v_n^{\bar{X}})_H$ and $\frac{1}{2} \|\hat{G}_\sigma (\bar{X}_n + X_n^{ref}) - \mathcal{M}^{\sigma,\bar{X}_n + X_n^{ref}}\|_{\mathcal{L}_2(\mathcal{W}, V)}^2$ can be estimated as follows,
    \begin{multline}
        \Big(\hat{F}_\sigma (\bar{X}_n + X_n^{ref}), A v_n^{\bar{X}}\Big)_H \leq C_\xi(1+ \| \bar{X}_n\|_V^2) + \xi \| \bar{X}_n\|_{\mathcal{D}(A)}^2 \\
        + \frac{1}{2}\Big(D_{\mu,\nu} a_{BT} \nabla_3 (\nabla_3 v_n^{\bar{X}}), \nabla_3 (\nabla_3 v_n^{\bar{X}})\Big)_{L^2} + \frac{1}{2}\Big(D_{\mu,\nu} \overline{a_{BC}} \nabla_3 (\nabla_3 \bar{v}_n^{\bar{X}}), \nabla_3 (\nabla_3 \bar{v}_n^{\bar{X}})\Big)_{L^2},
    \end{multline}
    and
    \begin{multline}
        \frac{1}{2} \|\hat{G}_\sigma (\bar{X}_n + X_n^{ref}) - \mathcal{M}^{\sigma,\bar{X}_n + X_n^{ref}}\|_{\mathcal{L}_2(\mathcal{W},V)}^2 \leq C_\xi(1+\|\bar{X}_n\|_V^2) + \xi \|\bar{X}_n\|_{\mathcal{D}(A)}^2 \\
        + \frac{1}{2}\Big(a_{BT} \nabla_3 (\nabla_3 v_n^{\bar{X}}), \nabla_3 (\nabla_3 v_n^{\bar{X}})\Big)_{L^2} + \frac{1}{2} \Big(\overline{a_{BC}} \nabla_3 (\nabla_3 \bar{v}_n^{\bar{X}}), \nabla_3 (\nabla_3 \bar{v}_n^{\bar{X}})\Big)_{L^2}.
    \end{multline}
    Therefore, we conclude that the existence of a pathwise solution of $(\hat{\mathcal{P}}_K)$ holds in the same space, adapting the rest of the proof being immediate.
\end{remark}

\subsection{Uniqueness of the solution} \label{subsec-uniqueness}

In this subsection, we show that pathwise solutions to the problem $(\mathcal{P}^\kappa)$ are unique. Our proof is analogous to the one of Proposition 3.5 in \cite{BS_2021}.

Let $(\mathcal{S},X^1)$ and $(\mathcal{S},X^2)$ two pathwise solutions of $(\mathcal{P}^\kappa)$ on the same stochastic basis. Let $R:=X^1-X^2$. Then remark that, for all $T>0$,
\begin{align}
    R \in L^2([0,T],\mathcal{D}(A)) \cap L^\infty([0,T],V), \quad \text{a.s.}
\end{align}
Substract the equations satisfied by $X^1$ and $X^2$, to get
\begin{multline}
    d_t R + A R dt + \Big[ B (X^1,X^1) - B (X^2,X^2)\Big] dt + C (R) dt + \frac{1}{\rho_0} \mathbf{P} \nabla_H (p^{X^1} - p^{X^2}) dt \\
    + [F_\sigma (R) - F_\sigma(0)] dt = [G_\sigma (R) - G_\sigma(0)] dW_t - \frac{1}{\rho_0} \mathbf{P} \nabla_H (dp_t^{\sigma,X^1} - dp_t^{\sigma,X^2}),
\end{multline}
with $R(0) = X^1(0) - X^2(0).$ Using It\={o}'s lemma, we reach
\begin{align}
    d_t \|R\|_H^2 = &- 2 (R, AR + B(X^1) - B(X^2) + C (R) +  \frac{1}{\rho_0} \mathbf{P} \nabla_H (p^{X^1} - p^{X^2}))_H dt \\
    &+ 2(R, [F_\sigma (R) - F_\sigma(0)])_H dt + \|[G_\sigma (R) - G_\sigma(0)] - [\mathcal{M}^{\sigma, X^1} - \mathcal{M}^{\sigma, X^2}]\|_{\mathcal{L}_2(\mathcal{W},H)}^2 dt \nonumber\\
    &+ 2 (R, [G_\sigma (R) - G_\sigma(0)] dW_t - \frac{1}{\rho_0} \mathbf{P} \nabla_H (dp_t^{\sigma,X^1} - dp_t^{\sigma,X^2}))_H  \nonumber.
\end{align}
Let $\eta, \zeta$ two stopping times such that with $0 \leq \eta < \zeta < \tau$. Integrate between $\eta$ and $\zeta$, then take the supremum and the expectation, so that
\begin{align}
    \E\Big[ \sup_{[\eta,\zeta]} &\|R\|_H^2 + c \int_{\eta}^{\zeta} \|R\|_{V}^2 ds \Big] \leq \E\Big[ \|R(\eta)\|_H^2 \Big] + 2\E\Big[ \int_{\eta}^{\zeta} |(R, B(X^1) - B(X^2))_H| ds \Big] \nonumber\\
    &+ 2\E\Big[ \int_{\eta}^{\zeta} |(R, CR + \frac{1}{\rho_0} \mathbf{P} \nabla_H (p^{X^1} - p^{X^2}) + [F_\sigma(R) - F_\sigma(0)] )_H| ds \Big] \nonumber\\
    &+ \E\Big[ \int_{\eta}^{\zeta} \|[G_\sigma (R) - G_\sigma(0)] - [\mathcal{M}^{\sigma, X^1} - \mathcal{M}^{\sigma, X^2}]\|_{\mathcal{L}_2(\mathcal{W},H)}^2 ds \Big] \nonumber \\
    &+ 2\E\Big[ \sup_{[\eta,\zeta]} \Big|\int_{\eta}^{t} (R, [G_\sigma (R) - G_\sigma(0)] dW_t - \frac{1}{\rho_0} \mathbf{P} \nabla_H (dp_t^{\sigma,X^1} - dp_t^{\sigma,X^2}))_H\Big| \Big] \nonumber\\
    &:= J_1 + J_2 + J_3 + J_4.
\end{align}
For $J_1$, by Hölder's inequality,
\begin{align}
    J_1 &= 2\E\Big[ \int_{\eta}^{\zeta} |(R, B(R,X^1))_H + (R, B(X^2,R))_H| ds \Big] \nonumber\\
    &\leq C \E\Big[ \int_{\eta}^{\zeta} \|R\|_H \|R\|_{V} \Big( \|X^1\|_{\mathcal{D}(A)} +  \|X^2\|_{\mathcal{D}(A)} \Big)  ds\Big] \nonumber \\
    &\leq C_\xi \E\Big[ \int_{\eta}^{\zeta} \|R\|_H^2 \Big( \|X^1\|_{\mathcal{D}(A)}^2 + \|X^2\|_{\mathcal{D}(A)}^2 \Big) ds\Big] + \xi \E[\int_{\eta}^{\zeta} \|R\|_{V}^2 ds] .
\end{align}
For $J_2$ and $J_3$, repeating the arguments of subsection \ref{subsec-strong-energ},
\begin{align}
    J_2 + J_3 &\leq C_\xi \E\Big[\int_{\eta}^{\zeta}\|R\|_H^2 \Big] + \xi \E\Big[ \int_{\eta}^{\zeta}\|R\|_V^2 \Big].
\end{align}
For $J_4$,
\begin{multline}
    J_4 \leq C \sum_{k=0}^\infty \E\Bigg[\Bigg(\int_{\eta}^{\zeta} \Big(R, \big[G_{\sigma} R - G_{\sigma}(0)\big]\phi_k- \mathbf{P} \nabla_H [\pi_k^{U_1} - \pi_k^{U_2} ] \Big)_H^2 \: dr\Bigg)^{1/2}\Bigg] \\
    \leq C_\xi \E\Big[\int_{\eta}^{\zeta} \|R\|_H^2\Big] + \xi \E\Big[ \int_{\eta}^{\zeta} \|R\|_{V}^2\Big].
\end{multline}
Then, using the stochastic Grönwall lemma (see Lemma 5.3 in \cite{GHZ_2009}),
\begin{align}
    \E \Big[ \sup_{[0,T]}\|R\|_H^2 \Big] = 0.
\end{align}
Therefore, $R=0$ a.s., so the solution to $(\mathcal{P}^\kappa)$ is pathwise unique.
\CQFD

\subsection{Maximal solution with improved regularity} \label{subsec-max-sol-and-regularity}

Using similar tightness arguments as in subsection \ref{subsec-tension}, together with an argument due to Gyöngy and Krylov as in \cite{BS_2021,DGHT_2011} (see also \cite{GK_2022}), we deduce that there exists a unique solution $X$ to $(\mathcal{P}^\kappa)$, such that $X \in L^\infty([0,T], V) \cap L^2([0,T], \mathcal{D}(A))$ almost surely. In addition, the convergence $X_n \rightarrow X$ happens in probability in $L^2([0,T], V)$. Hence, by setting
\begin{equation}
    \tau^\kappa := \Big\{t>0 \Big| \|X - X_{ref}\|_V \geq \kappa/2 \Big\},
\end{equation}
we infer that $U = X_{\tau^\kappa \wedge \bcdot}$ is a (pathwise unique) local-in-time solution to the problem $(\mathcal{P}_K)$. By repeating the arguments of \cite{GHZ_2009}, we conclude that there exists a maximal solution $U$ to $(\mathcal{P}_K)$, associated to the stopping time $\tau$. Moreover, this solution is global-in-time whenever the following holds,
\begin{equation}
    \forall T>0, \quad  \sup_{[0,\tau \wedge T]} \|U\|_V^2 + \int_0^{\tau \wedge T} \|U\|_{\mathcal{D}(A)}^2 < \infty, \text{ a.s.} \label{eq-global-in-time-condition}
\end{equation}
In the rest of the subsection, we show the time-continuity of $U$, that is, for all stopping time $0< \tau' < \tau$ and $T>0$,
\begin{equation}
    U_{\tau' \wedge \bcdot}^* \in C([0,T], V).
\end{equation}
For this purpose, we use an argument proposed in \cite{DGHT_2011,DHM_2023}. Consider the following equation, 
\begin{equation}
    d_t z(t) + Az(t) dt + \mathbf{P} \nabla_H p^{U^*}dt = G_\sigma (U^*)d\overline{W}_t - \frac{1}{\rho_0} \mathbf{P} \nabla_H dp_t^{\sigma,U^*},
\end{equation}
with the initial condition $z(0) = U^*(0)$, and keeping the same notation $U^*$ as before. Again, $p^{U^*}$ and $p_t^{\sigma,U^*}$ are respectively the pressure bounded variation and martingale terms associated to $U^*$. This equation has a unique solution satisfying $z(\tau' \wedge \bcdot) \in L^2(\Omega, C([0,T],V))$.

Using It\={o}'s lemma with $x \mapsto \frac{1}{2} \|x\|_V^2 $, the following estimate is straightforward,
\begin{align}
    \E[\sup_{0\leq t \leq T} \|z(\tau' \wedge t)\|_V^2] + \E[\int_0^{T} \|z(\tau' \wedge t)\|_{\mathcal{D}(A)}^2 dt] < \infty.
\end{align}
Thus, $z \in L^2([0,T],\mathcal{D}(A))$ almost surely. Let $\tilde{U}:= U^*-z \in L^2([0,T],\mathcal{D}(A))$, so that
\begin{equation}
    d_t \tilde{U} + A\tilde{U}dt + B(\tilde{U} + z)dt + C(\tilde{U} + z) dt + F_\sigma (\tilde{U}+z)dt =0.
\end{equation}
Since $\tilde{U},z \in L^2([0,T],\mathcal{D}(A))$, we have 
\begin{align}
    \int_0^T \|B(\tilde{U}+z)(\tau' \wedge s)\|_H^2  ds &\leq C \int_0^T \|(\tilde{U}+z)(\tau' \wedge s)\|_{V}^2 \|(\tilde{U}+z)(\tau' \wedge s)\|_{\mathcal{D}(A)}^2\\
    &\leq C \Big( \sup_{[0,T]}\|(\tilde{U}+z)(\tau' \wedge s)\|_{V}^2 \Big) \int_0^T \|(\tilde{U}+z)(\tau' \wedge s)\|_{\mathcal{D}(A)}^2 \leq C, \nonumber
\end{align}
and
\begin{align}
    \int_0^T \|A\tilde{U}(\tau' \wedge s)\|_{H}^2 ds &\leq C \int_0^T \|\tilde{U}(\tau' \wedge s)\|_{\mathcal{D}(A)}^2 \leq C, \\
    \int_0^T \|F_\sigma(\tilde{U}+z)(\tau' \wedge s)\|_{H}^2 ds &\leq C \int_0^T \|(\tilde{U}+z)(\tau' \wedge s)\|_{\mathcal{D}(A)}^2 \leq C.
\end{align}
Therefore $A\tilde{U}_{(\tau' \wedge \bcdot)}, B(\tilde{U} + z)_{(\tau' \wedge \bcdot)}, F_\sigma(\tilde{U}+z)_{(\tau' \wedge \bcdot)} \in L^2([0,T],H)$ almost surely, consequently $(\frac{d}{dt} \tilde{U})_{(\tau' \wedge \bcdot)} \in L^2([0,T],H)$ almost surely as well. Since $\tilde{U} \in L^2([0,T],\mathcal{D}(A))$ $\overline{\mathbbm{P}}$-a.s, we can apply a result from \cite{Temam_2001} (Chapter 3, Lemma 1.2), to reach 
\begin{align}
    \tilde{U}_{(\tau' \wedge \bcdot)} \in L^2(\Omega, C([0,T],V)).
\end{align}
As a consequence 
\begin{align}
    z_{(\tau' \wedge \bcdot)},\tilde{U}_{(\tau' \wedge \bcdot)},U^*_{(\tau' \wedge \bcdot)} \in L^2(\Omega,L^2([0,T],\mathcal{D}(A))) \cap L^2(\Omega,C([0,T],V)).
\end{align}

\CQFD

Consequently, Theorem \ref{theorem-filtered-weak} holds for $(\mathcal{P}_K)$. As shown in the previous remarks, the same holds for $(\hat{\mathcal{P}}_K)$, with minor changes in the proof.

\section{Existence of a global-in-time pathwise solution under \modif{stronger} assumptions} \label{section-globality}

In this section, we show that the approximated low-pass filtered problem $(\hat{\mathcal{P}}_K)$ -- equation \eqref{problem-approx-PK} -- admits a global pathwise solution. This implies that $(\mathcal{P}_K)$ is well-posed as well for a purely barotropic noise $\sigma dW_t$, since it corresponds to $(\hat{\mathcal{P}}_K)$ with $\sigma_{BC} dW_t = 0$. We remind that, as we assume $K \in H^3$, the functional $\mathcal{C}_K [\bcdot] = K * [\bcdot]$ fulfils
\begin{equation}
    \|\mathcal{C}_K \bcdot \|_{H^k} \leq C (1+\|\bcdot\|_{L^2}), \quad \forall k \in \{0,1,2,3\}.
\end{equation}
This remark will be used further to establish energy estimates. In addition, we remind that the noise is assumed to be regular enough in the sense of equation \eqref{smoothness-noise}, that is
\begin{align}
    \sup_{t\in [0,T]} \sum_{k=0}^\infty \|\phi_k\|_{H^4(\mathcal{S}, \R^3)}^2 < \infty, \quad u_s \in L^\infty\Big([0,T], H^4(\mathcal{S}, \R^3)\Big)\\
   d_t u_s \in L^\infty\Big([0,T],H^3(\mathcal{S}, \R^3)\Big), \quad a \nabla u_s \in L^\infty\Big([0,T],H^2(\mathcal{S},\R^3)\Big). \nonumber
\end{align}
Let $U$ be a maximal pathwise solution of $(\hat{\mathcal{P}}_K)$, and denote by $\tau$ its associated stopping time. As mentioned in subsection \ref{subsec-max-sol-and-regularity}, it suffices to show that equation \eqref{eq-global-in-time-condition} holds. The following proof relies on a decomposition of the velocity between the barotropic and the baroclinic modes, which originates from \cite{CT_2007}. Remind first that the barotropic and baroclinic operators $\mathcal{A}_2$ and $\mathcal{R}$ are defined as
\begin{align}
    (\mathcal{A}_2 v)(x,y) = \frac{1}{h} \int_{-h}^0 v(x,y,z') \: dz', \quad (\mathcal{A} v)(x,y,z) =  (\mathcal{A}_2 v)(x,y), \quad \mathcal{R} v = (I - \mathcal{A})v.
\end{align}
Then, decompose the velocity $v$ as follows
\begin{align}
    v= \bar{v} + \tilde{v}, \quad \bar{v} = \mathcal{A}_2 v, \quad \tilde{v} = \mathcal{R} v.
\end{align}
Using the same computation as in \cite{CT_2007}, one shows that the following holds on the (2D) domain $\mathcal{S}_H$ for $\bar{v}$, using that $\overline{\mathbf{P} \nabla_H p} = \mathbf{P}_{2D} \nabla_H \bar{p} =0$, and similarly that $\overline{\mathbf{P} \nabla_H dp_t^\sigma} =0$,
\begin{multline}
    d_t\bar{v} + \mathbf{P}_{2D} \Big[ - \mu_v \Delta_H \bar{v} + (\bar{v} \bcdot \nabla_H)\bar{v} + \mathcal{A}_2 [(\tilde{v}\bcdot \nabla_H) \tilde{v} + (\nabla_H \bcdot \tilde{v}) \tilde{v}] + f \bold{k} \times \bar{v}\Big] dt\\
    + \mathcal{A}_2 \hat{F}_\sigma^v(U) dt = \mathcal{A}_2 \hat{G}_\sigma^v(U) dW_t,
\end{multline}
with $\bar{v}_0 = \mathcal{A}_2 v_0$, $\nabla_H \bcdot \bar{v} =0$ on $\mathcal{S}_H$ and $\bar{v}\bcdot\bold{n} = \frac{\partial\bar{v}}{\bold{n}} \times \bold{n}=0$ on $\partial \mathcal{S}_H$. Moreover, $\tilde{v}$ fulfils the following on $\mathcal{S}$
\begin{multline}
    d_t\tilde{v} + \Big[ A \tilde{v} + B(\tilde{v}) + (\bar{v}\bcdot \nabla_H) \tilde{v} + (\tilde{v}\bcdot \nabla_H) \bar{v} - \mathcal{A} [(\tilde{v}\bcdot \nabla_H) \tilde{v} + (\nabla_H \bcdot \tilde{v}) \tilde{v}] + f \bold{k} \times \tilde{v} - \frac{1}{\rho_0 } \widetilde{\mathbf{P} \nabla_H p} \Big]dt \\
    + \mathcal{R} \hat{F}_\sigma^v(U) dt = \mathcal{R} \hat{G}_\sigma^v(U) dW_t - \frac{1}{\rho_0 } \widetilde{\mathbf{P} \nabla_H dp_t^\sigma},
\end{multline}
with $\tilde{v}_0 = \mathcal{R} v_0$, and $\tilde{v} \bcdot\bold{n} = \frac{\partial\bar{v}}{\bold{n}} \times \bold{n} = 0$ on $\Gamma_l$, $\partial_z \tilde{v} =0$ on $\Gamma_h \cup \Gamma_b$. By the energy estimate \eqref{estimation-énerg} of section \ref{section-martingale}, the following estimate holds for all $p \geq 2$ and all $T>0$,
\begin{align}
    \E\Bigg[\sup_{s \in [0,T \wedge \tau]} \|U_s\|_H^p + \int_0^{T \wedge \tau} \|U_s\|_H^{p-2} \|U_s\|_V^2 dr \Bigg] \leq C.
\end{align}

\subsection{$H^1$ estimate for the barotropic velocity}

Apply It\={o}'s lemma with the function $\|\mathcal{A} A^{1/2} \bcdot \|_H^2$ to the problem $(\hat{\mathcal{P}}_K)$, so that
\begin{align}
    d_t\|\bar{v}\|_V^2 + \|A \bar{v}\|_H^2 dt \leq &-2 \Big(\mathbf{P}_{2D} \mathcal{A}_2[ (\tilde{v} \bcdot \nabla_H)\tilde{v} + \tilde{v} (\nabla_H \bcdot \tilde{v}) ], A \bar{v}\Big)_H dt  \label{ito-barotropic} \\
    &- 2 \Big(\mathbf{P}_{2D}(\bar{v} \bcdot \nabla_H) \bar{v}, A \bar{v}\Big)_H dt - 2 \Big(f \bold{k} \times \bar{v}, A \bar{v}\Big)_H dt - 2 \Big(\mathcal{A}_2 \hat{F}_\sigma^v(U), A \bar{v}\Big)_H dt \nonumber\\
    &+ \Big\|\mathcal{A}_2 \hat{G}_\sigma^v(U)\Big\|_{\mathcal{L}_2(\mathcal{W}, V)}^2 dt 
    + 2 \sum_{k=1}^\infty \Big( \mathcal{A}_2 \hat{G}_\sigma^v(U) \phi_k, A \bar{v}\Big)_H d\beta^k \nonumber\\
    := &\sum_{j=1}^5 I_j \: dt + \sum_{k=0}^\infty I_6^k d\beta^k. \nonumber
\end{align}
Using the same arguments as in \cite{CT_2007}, we have
\begin{align}
    I_1 \leq C \Big(\int_\mathcal{S} |\tilde{v}|^2 |\nabla_3 \tilde{v}|^2 d\mathcal{S}\Big)^{1/2} \|A \bar{v}\|_H, \quad \text{ and } \quad I_2 \leq C \|\bar{v}\|_H^{1/2} \|\bar{v}\|_V \|A \bar{v}\|_H^{3/2}.
\end{align}
Hence, by Young's inequality, for all $\xi >0$, there exists $C_\xi>0$ such that
\begin{align}
    I_1 \leq \xi \|A \bar{v}\|_H^2 + C_\xi \int_\mathcal{S} |\tilde{v}|^2 |\nabla_3 \tilde{v}|^2 d\mathcal{S}, \quad \text{ and } \quad I_2 &\leq \xi \|A \bar{v}\|_H^2 + C_\xi \|U\|_H^2 \|U\|_V^2 \|\bar{v}\|_V^2,
\end{align}
using that $\|\bar{v}\|_V \leq \|U\|_V$. Once again, the notation $C_\xi$ is carried from line to line, even if it may refer to different $\xi$-dependent constants. Furthermore,
\begin{align}
    I_3 &\leq \xi \|A \bar{v}\|_H^2 + C_\xi\|U\|_H^2, \\
    I_4 & \leq \xi \|A \bar{v}\|_H^2 + C_\xi (1+ \|U\|_V^2) - \mu (\nabla_3 \bcdot (a_{BT} \nabla_3 \bar{v}), A \bar{v})_{L^2} - (\nabla_3 \bcdot (\overline{a_{BC}} \nabla_3 \bar{v}), A \bar{v})_{L^2} \nonumber\\
    &\leq \xi \|A \bar{v}\|_H^2 + C_\xi (1+ \|U\|_V^2) - \mu (\nabla_3 \bcdot (a_{BT} \nabla_3 \bar{v}), A \bar{v})_{L^2} - \underbrace{(\overline{a_{BC}} \nabla_3 (\nabla_3 \bar{v}), \nabla_3 (\nabla_3 \bar{v}))_{L^2}}_{\geq 0}. \label{barotropic-estimate-I4}
\end{align}
For $I_5$, define $\psi_k = -\mathcal{A}_2 (A \phi_k^H + B(\phi_k,v_s) + \Gamma \phi_k^H)$, so we have
\begin{align}
    \|\mathcal{A}_2 &G_\sigma^v(U)\|_{\mathcal{L}_2(\mathcal{W},V)}^2 = \sum_{k=0}^\infty \|\psi_k - \mathbf{P}_{2D} [\bar{\phi}_k^H \bcdot \nabla_H \bar{v}] - \mathbf{P}_{2D} \mathcal{A}_2(\phi_k^z \partial_z v)\|_V^2 \nonumber\\
    &\leq \sum_{k=0}^\infty \|\psi_k\|_V^2 + \|\mathbf{P}_{2D}[\bar{\phi}_k^H \bcdot \nabla_H \bar{v} + \mathcal{A}_2(\phi_k^z \partial_z v)]\|_V^2 + 2 \underbrace{|(\mathbf{P}_{2D}[\bar{\phi}_k^H \bcdot \nabla_H \bar{v} + \mathcal{A}_2(\phi_k^z \partial_z v)], \psi_k)_V|}_{\leq |(b(\bar{\phi}_k^H, \bar{v}, \Delta_3 \psi_k)| \leq \|\phi_k\|_{L^\infty} \|\psi_k\|_{H^2} \|U\|_V} \nonumber\\
    &\leq \sum_{k=0}^\infty \|\mathbf{P}_{2D}[\bar{\phi}_k \bcdot \nabla_3 \bar{v} + \mathcal{A}_2(\phi_k^z \partial_z v)]\|_V^2 + C (1+\|U\|_V^2).
\end{align}
Moreover, remark that
\begin{align*}
    \|\mathbf{P}_{2D}[\bar{\phi}_k^H \bcdot \nabla_H \bar{v} + \mathcal{A}_2(\phi_k^z \partial_z v)]\|_V^2 &\leq \|\mathbf{P}_{2D}[\bar{\phi}_k \bcdot \nabla_3 \bar{v}]\|_V^2 + \xi \|\phi_k\|_{H^4}^2 \|\bar{v}\|_{\mathcal{D}(A)}^2 \\
    &+ C_\xi \|\phi_k\|_{H^4}^2 (\|\partial_z v\|_V^2 + \|v\|_{V}^2+1).
\end{align*}
Since
$$\|\mathbf{P}_{2D}[\bar{\phi}_k \bcdot \nabla_3 \bar{v}]\|_V^2 \leq \mu (\nabla_3 \bcdot (\bar{\phi}_k \bar{\phi}_k^\tr \nabla_3 \bar{v}), A \bar{v})_{L^2} + C_\xi\|\phi_k\|_{H^4}^2 (\|U\|_V^2 + \|\partial_z v\|_V^2) + \xi \|\phi_k\|_{H^4}^2 \|\bar{v}\|_{\mathcal{D}(A)}^2,$$
we derive
\begin{align}
    I_5 &\leq \mu (\nabla_3 \bcdot (a_3 \nabla_3 \bar{v}), A \bar{v})_{L^2} + C_\xi (\|\partial_z v\|_V^2 + \|U\|_V^2 + 1) + \xi \|\bar{v}\|_{\mathcal{D}(A)}^2.
\end{align}
Gathering the previous estimates into \eqref{ito-barotropic}, then time integrating between two stopping times $0<\eta<\zeta<T \wedge \tau$, and taking the expectation, we reach
\begin{multline}
    \E\Big[\|\bar{v}\|_{V}^2(\zeta)\Big] + (c - \alpha_1) \E\Big[ \int_\eta^\zeta \|A \bar{v}\|_H^2 \: ds\Big]
    \leq C \E\Big[\|\bar{v}(\eta)\|_{V}^2 + 1\Big] \label{ito-barotrope-simple}\\
    + C\E\Big[\int_\eta^\zeta \|U\|_H^2 \|U\|_V^2 \|\bar{v}\|_{V}^2 + \|U\|_V^2 + \int_\mathcal{S} |\tilde{v}|^2 |\nabla_3 \tilde{v}|^2 + \|\partial_z v\|_V^2 ds\Big],
\end{multline}
for some $\alpha_1 \in (0,c)$. Similarly, time integrating, then taking the supremum over $[\eta,\zeta]$ and the expectation, we have also
\begin{multline}
    \E\Big[\sup_{[\eta,\zeta]} \|\bar{v}\|_{V}^2\Big] + (c - \alpha_1) \E\Big[ \int_\eta^\zeta \|A \bar{v}\|_H^2 \: ds\Big]
    \leq C \E\Big[\|\bar{v}(\eta)\|_{V}^2 + 1\Big]\label{ito-barotrope-sup-2} \\
    + C\E\Big[\int_\eta^\zeta \|U\|_H^2 \|U\|_V^2 \|\bar{v}\|_{V}^2 + \|U\|_V^2 + \int_\mathcal{S} |\tilde{v}|^2 |\nabla_3 \tilde{v}|^2 + \|\partial_z v\|_V^2 ds\Big] + \E\Bigg[ \sup_{[\eta,\zeta]} \Big|\int_\eta^s \sum_{k=0}^\infty I_6^k d\beta_r^k\Big|\Bigg].
\end{multline}
Additionally, by the Burkholder-Davis-Gundy inequality (Theorem 4.36 in \cite{book_DPZ_2014}),
\begin{align}
    \E\Bigg[ \sup_{[\eta,\zeta]} &\Big|\int_\eta^s \sum_{k=0}^\infty I_6^k d\beta_r^k\Big|\Bigg] \leq C \E\Bigg[\Bigg( \int_\eta^\zeta \sum_{k=0}^\infty \Big|\Big(\psi_k - \bar{\phi}_k \bcdot \nabla \bar{v} - \mathcal{A}_2(\phi_k^z \partial_z v), A \bar{v} \Big)_{L^2}\Big|^2 ds \Bigg)^{1/2}\Bigg] \nonumber\\
    &\leq C \E\Bigg[\Bigg( \int_\eta^\zeta \sum_{k=0}^\infty |(\psi_k,\bar{v})_V|^2 + |(\bar{\phi}_k \bcdot \nabla \bar{v}, (-\Delta_H)\bar{v})_{L^2}|^2 + |(\mathcal{A}_2(\phi_k^z \partial_z v), A \bar{v} )_{L^2}|^2 ds \Bigg)^{1/2}\Bigg] \nonumber\\
    &\leq C  \E\Bigg[\Bigg( \int_\eta^\zeta \sum_{k=0}^\infty \|\psi_k\|_V^2 \|\bar{v}\|_V^2 + \|\bar{\phi}_k\|_{L^\infty}^2 \|A\bar{v}\|_H^2 \|\bar{v}\|_V^2 + \|\phi_k^z\|_{L^\infty}^2\|\partial_z v\|_V^2 \|\bar{v}\|_V^2 ds \Bigg)^{1/2}\Bigg] \nonumber\\
    &\leq \xi \E\Big[\sup_{[\eta,\zeta]} \|\bar{v}\|_{V}^2\Big] + C_\xi \E\Bigg[1+ \int_\eta^\zeta \|A\bar{v}\|_H^2 + \|\partial_z v\|_V^2 ds \Bigg]. 
\end{align}
Therefore, by \eqref{ito-barotrope-simple}, we reach
\begin{multline}
    \E\Big[\sup_{[\eta,\zeta]} \|\bar{v}\|_{V}^2\Big] + \E\Big[ \int_\eta^\zeta \|A \bar{v}\|_H^2 \: ds\Big]
    \leq C_1 \E\Big[\|\bar{v}(\eta)\|_{V}^2 + 1\Big] \label{estimation-barotrope}\\
    + C_1 \E\Big[\int_\eta^\zeta \|U\|_H^2 \|U\|_V^2 \|\bar{v}\|_{V}^2 + \|U\|_V^2 + \int_\mathcal{S} |\tilde{v}|^2 |\nabla_3 \tilde{v}|^2 + \|\partial_z v\|_V^2 ds\Big].
\end{multline}

\subsection{$L^2$ estimate for the vertical gradient of velocity} \label{subsec-vertical-gradient}

Now remark that
\begin{align}
    \partial_z [(v \bcdot \nabla_H) v + w(v) \partial_z v ] &= (\partial_z v \bcdot \nabla_H) v + (v \bcdot \nabla_H) \partial_z v - (\nabla_H \bcdot v) \partial_z v + w(v) \partial_{zz} v,\\
    \text{ and } \quad &((v \bcdot \nabla_H) \partial_z v + w(v) \partial_{zz} v, \partial_z v)_{L^2} = 0. \nonumber
\end{align}
So
\begin{align}
   (\partial_z [(v \bcdot \nabla_H) v + w(v) \partial_z v ], \partial_z v)_H = ((\partial_z v \bcdot \nabla_H) v - (\nabla_H \bcdot v) \partial_z v, \partial_z v)_H.
\end{align}
Applying It\={o}'s lemma with $\|\partial_z \bcdot \|_H^2$ to the problem $(\hat{\mathcal{P}}_K)$, we reach
\begin{align}
    d_t\|\partial_z v\|_{L^2}^2 + 2 &c \|\nabla_3 \partial_z v\|_{L^2} ^2 dt \leq - 2 \Big((\partial_z v \bcdot \nabla_H) v - (\nabla_H \bcdot v) \partial_z v + g [\beta_T \nabla_H T + \beta_S \nabla_H S],\partial_z v\Big)_{L^2} dt \nonumber\\
    &+ \Big(\mathbf{P} \nabla_H \Big[\nabla_3 \bcdot (a_3^K \nabla_3 (w(v)+w_s))+ K*[(u_s \bcdot \nabla_3) (w(v)+w_s)]\Big], \partial_z v\Big)_{L^2} dt \nonumber\\
    &- 2 \Big(\partial_z \hat{F}_\sigma^v(U), \partial_z v\Big)_{L^2} dt + \Big\|\partial_z \hat{G}_\sigma^v(U)- \partial_z\mathcal{M}^{U,\sigma}\Big\|_{\mathcal{L}_2(\mathcal{W},{L^2})}^2 dt \nonumber\\
    &+ 2 \sum_{k=0}^\infty \int_\mathcal{S} \Big[\partial_z \hat{G}_\sigma^v (U) (\phi_k) - \partial_z{\nabla_H \pi_k}\Big]\bcdot \partial_z v \: d\beta_t^k \nonumber\\
    &:= \sum_{j=1}^4 I_j dt + \sum_{k=0}^\infty I_5^k \: d\beta_t^k,
\end{align}
with $c=\min_i\{\mu_i, \nu_i\}$ for example, since $(f \bold{k} \times \partial_z  v, \partial_z v )_H =0$. Repeating the argument of \cite{HH_2020} (Step 2, p.24), for all $\xi > 0$,
\begin{align}
    I_1 \leq C_\xi (1+\|v\|_V^2) \|\partial_z v\|_{L^2}^2 + C_\xi \int_\mathcal{S} |\nabla \tilde{v}|^2 |\tilde{v}|^2 + \xi \|\nabla_3 \partial_z v\|_{L^2}^2.
\end{align}
\modif{
For $I_2$, by the regularising properties of $K$ and Sobolev embeddings, since $u_s$ and $\phi_k$ are divergence-free,
\begin{align}
    I_2 &= \Big( \mathbf{P} \nabla_H \Big[ \nabla_3 \bcdot (a_3^K \nabla_3 (w(v)+w_s)) + K*[u_s \bcdot \nabla_3 (w(v)+w_s)] \Big], \partial_{z} v\Big)_{L^2}\\
    &\leq \sum_k \|\nabla_H [ \nabla_3 \bcdot (\phi_k \mathcal{C}_K^* \mathcal{C}_K [\nabla_3 \bcdot (w(v)+w_s)\phi_k]) \|_{L^2} \|\partial_z v\|_{L^2} + \|K\|_{H^2} \|u_s\|_{H^{2}} \|w(v)\|_{L^2} \|\partial_z v\|_{L^2} \nonumber\\
    &\leq \Big(\|K\|_{H^2}^2(\|w_s\|_{L^2} + \|w(v)\|_{L^2}) \sum_k \|\phi_k\|_{H^4}^2 + \|K\|_{H^2} \|u_s\|_{H^{2}} \|w(v)\|_{L^2}\Big) \|\partial_z v\|_{L^2}. \nonumber\\
    &\qquad \qquad \qquad \qquad \qquad \qquad \qquad \qquad \qquad \qquad \qquad \qquad \qquad \qquad \qquad \leq C (1+\|U\|_V^2). \nonumber
\end{align}
}
In addition,
\begin{align}
    I_3 &\leq C - (a_{BT} \nabla_3 \partial_z v, \nabla_3 \partial_z v)_{L^2} \\
    I_4 &\leq \Big[ C + \underbrace{\Big(\sum_{k=0}^\infty \partial_z \phi_k^{BC} (\partial_z \phi_k^{BC})^\tr \nabla_3 v, \nabla_3 v\Big)_{L^2}}_{\leq C \|U\|_V^2} + (a_{BT} \nabla_3 \partial_z v, \nabla_3 \partial_z v)_{L^2} \Big] \\
    &+ \underbrace{\sum_{k=0}^\infty \|\nabla_H \Big[K * [\phi_k \bcdot \nabla_3 (w(v)+ w_s)] + A \phi_k^z\Big]\|_{L^2}^2}_{\leq C (1+ \sum_{k=0}^\infty \|(\nabla_3 \bcdot \nabla_H K) * (\phi_k w(v))\|_{L^2}^2) \leq C (1+ \|K\|_{H^2}^2\|v\|_V^2)} \nonumber\\
    &+ \Big|\underbrace{2\sum_{k=0}^\infty (\partial_z[\phi_k \bcdot \nabla_3 v + \phi_k \bcdot \nabla_3 v_s + A \phi_k^H + \Gamma \phi_k^H], \mathbf{P} \nabla_H[ K* [\phi_k \bcdot \nabla_3 (w(v)+w_s)] + A \phi_k^z])_{L^2}}_{\leq C_\xi(1+ \|v\|_V^2) + \xi \|\nabla_3 \partial_z v\|_{L^2}^2}\Big|.  \nonumber
\end{align}
Gathering the previous estimates, integrating and taking expectation, we have
\begin{align}
    \E\Big[ \|\partial_z v\|_{L^2}^2(\zeta)\Big] + (c - \alpha_2) \E\Big[ \int_\eta^\zeta &\|\nabla_3 \partial_z v\|_{L^2}^{2} \: ds\Big]
    \leq C \E\Big[\|\partial_z v(\eta)\|_{L^2}^2 +1\Big] \label{ito-gradient-z-simple} \\
    &+ C\E\Big[\int_\eta^\zeta (1+\|U\|_V^2) (1+\|\partial_z v\|_{L^2}^2) + \int_\mathcal{S} |\nabla \tilde{v}|^2 |\tilde{v}|^2 ds \Big]. \nonumber
\end{align}
Similarly, taking the supremum, we reach
\begin{multline}
    \E\Big[\sup_{[\eta,\zeta]} \|\partial_z v\|_{L^2}^2\Big] + (c - \alpha_2) \E\Big[ \int_\eta^\zeta \|\nabla_3 \partial_z v\|_{L^2}^{2} \: ds\Big]
    \leq C \E\Big[\|\partial_z v(\eta)\|_{H}^2 +1\Big] \\
    + C\E\Big[\int_\eta^\zeta (1+\|U\|_V^2) (1+\|\partial_z v\|_{L^2}^2) + \int_\mathcal{S} |\nabla \tilde{v}|^2 |\tilde{v}|^2 ds \Big] + \E\Bigg[ \sup_{[\eta,\zeta]} \Big|\int_\eta^s \sum_{k=0}^\infty I_6^k d\beta_r^k\Big|\Bigg].
\end{multline}
Denoting by $\psi_k = \mathbf{P}[\phi_k \bcdot \nabla_3 v_s] + A \phi_k^H + C \phi_k^H$ and $\psi_k^w = \phi_k \bcdot \nabla_3 w_s + A \phi_k^z$, and using the Burkholder-Davis-Gundy inequality (Theorem 4.36 in \cite{book_DPZ_2014}),
\begin{align}
    \E\Bigg[ &\sup_{[\eta,\zeta]} \Big|\int_\eta^s \sum_{k=0}^\infty I_6^k d\beta_r^k\Big|\Bigg] \leq C \E\Bigg[\Bigg( \int_\eta^\zeta \sum_{k=0}^\infty \Big|\Big(\partial_z \psi_k  - \partial_z(\phi_k \bcdot \nabla_3 v) - \partial_z{ \mathbf{P} \nabla_H \pi_k}, \partial_z v \Big)_{L^2}\Big|^2 ds \Bigg)^{1/2}\Bigg] \nonumber\\
    &\leq C \E\Bigg[\Bigg( \int_\eta^\zeta \sum_{k=0}^\infty \Big|\Big(\partial_z \psi_k + \partial_z \psi_k^w - \partial_z(\phi_k \bcdot \nabla_3 v) +  \mathbf{P} \nabla_H K* (\phi_k \bcdot \nabla_3 w(v) ), \partial_z v \Big)_{L^2}\Big|^2 ds \Bigg)^{1/2}\Bigg]. \label{ito-gradient-z-sup}
\end{align}
Remark that
\begin{multline*}
    \Big|\Big(\partial_z \psi_k + \partial_z \psi_k^w - \partial_z(\phi_k \bcdot \nabla_3 v) +  \mathbf{P} \nabla_H K*(\phi_k \bcdot \nabla_3 w(v) ), \partial_z v \Big)_{L^2}\Big|^2 \\
    \leq C\Big[\|\partial_z \psi_k + \partial_z \psi_k^w\|_{L^2}^2 \|\partial_z v\|_{L^2}^2+ \|\phi_k\|_{W^{1,\infty}}^2 \|\nabla_3 \partial_z v\|_{L^2}^2 \|\partial_z v\|_{L^2}^2
    + \Big|\Big( \nabla_H K*[\phi_k \bcdot \nabla_3 w(v)], \partial_z v \Big)_{L^2}\Big|^2 \Big],
\end{multline*}
and that
\begin{align*}
    \Big( \nabla_H K*[\phi_k \bcdot \nabla_3 w(v)], \partial_z v \Big)_{L^2} \leq C \|w(v)\|_{L^2} \|\partial_z v\|_{L^2} &\leq  C (\|v\|_V^2 +  \|\partial_z v\|_{L^2}^2) \leq  C (1 + \|\partial_z v\|_{L^2}^2).
\end{align*}
Consequently,
\begin{multline}
   C \E\Bigg[\Bigg( \int_\eta^\zeta \sum_{k=0}^\infty \Big|\Big(\partial_z \psi_k + \partial_z \psi_k^w - \partial_z(\phi_k \bcdot \nabla_3 v) +  \mathbf{P} \nabla_H K* (\phi_k \bcdot \nabla_3 w(v) ), \partial_z v \Big)_{L^2}\Big|^2 ds \Bigg)^{1/2}\Bigg] \\
   \leq C \E\Bigg[\Bigg(1+ \sup_{[\eta,\zeta]} \|\partial_z v\|_{L^2}\Bigg) \Bigg( \int_\eta^\zeta (1+\|\nabla_3 \partial_z v\|_{L^2}^2) ds \Bigg)^{1/2}\Bigg] \\ \leq \xi \E\Bigg[1+ \sup_{[\eta,\zeta]} \|\partial_z v\|_{L^2}^2 \Bigg] + C_\xi \E\Bigg[ \int_\eta^\zeta (1+\|\nabla_3\partial_z v\|_{L^2}^2) ds \Bigg].
\end{multline}
Thus, by equation \eqref{ito-gradient-z-simple},
\begin{align}
    \E\Big[\sup_{[\eta,\zeta]} \|\partial_z v\|_{L^2}^2\Big] + \E\Big[ \int_\eta^\zeta \|\nabla_3 \partial_z v\|_{L^2}^{2} \: ds\Big]
    \leq &C_2 \E\Big[\|\partial_z v(\eta)\|_{L^2}^2 +1\Big] \label{estimation-gradient-z} \\
    + &C_2 \E\Big[\int_\eta^\zeta (1+\|U\|_V^2) (1+\|\partial_z v\|_{L^2}^2) + \int_\mathcal{S} |\nabla \tilde{v}|^2 |\tilde{v}|^2 ds \Big]. \nonumber
\end{align}

\subsection{$L^4$ estimate for the baroclinic velocity} \label{subsec-baroclinic-estimate}

Remind first that the Leray type projector $\mathbf{P}$ corresponds to the identity on baroclinic vector fields. Also, using the regularising properties of the kernel $K$, we have
\begin{equation}
    \sum_{k=0}^\infty \|\nabla_H \pi_k\|_{L^4}^2 \leq C(1+\|w(v)\|_H^2) \leq C(1+\|v\|_V^2).
\end{equation}
Apply It\={o}'s lemma to the problem $(\hat{\mathcal{P}}_K)$ with $\|\mathcal{R} \bcdot\|_{L^4}^4$, so that
\begin{align}
    d_t\|\tilde{v}\|_{L^4}^4 &+ 4 \int_{\mathcal{S}} (\mu_v |\nabla_H\tilde{v}|^2 |\tilde{v}|^2 + \mu_v |\nabla_H|\tilde{v}|^2|^2 + \nu_v |\partial_z \tilde{v}|^2 |\tilde{v}|^2 + \nu_v |\partial_z |\tilde{v}|^2 |^2 ) d\mathcal{S} \: dt \\
    = &- 4 \int_{\mathcal{S}} |\tilde{v}|^2 \tilde{v} \bcdot \Big( (\tilde{v} \bcdot \nabla_H) \bar{v} - \mathcal{A}((\tilde{v} \bcdot \nabla_H) \tilde{v} +  \tilde{v} (\nabla_H\bcdot \tilde{v} )) + g\mathcal{R}\int_z^0 (\beta_T \nabla_H T + \beta_S \nabla_H S) \: dz' \Big) d\mathcal{S} \: dt \nonumber\\
    &+ 4 \int_{\mathcal{S}} |\tilde{v}|^2 \tilde{v} \bcdot \mathcal{R} \nabla_H \int_z^0 \Big[ \frac{1}{2} \nabla_3 \bcdot (a_3^K \nabla_3 (w(v)+w_s)) + K*[u_s \bcdot \nabla_3 (w(v)+w_s)]\Big] dz' \nonumber\\
    &- 4 \int_{\mathcal{S}} |\tilde{v}|^2 \tilde{v} \bcdot \mathcal{R} \hat{F}_\sigma^v(U) d\mathcal{S} \: dt  + 2 \sum_{k=0}^\infty \int_{\mathcal{S}} |\tilde{v}|^2 |\mathcal{R} \hat{G}_\sigma^v(U) \phi_k - \nabla_H \tilde{\pi}_k|^2 d\mathcal{S} \: dt \nonumber\\
    &+ 4 \int_{\mathcal{S}} |[\mathcal{R} \hat{G}_\sigma^v(U) - \mathcal{R}\mathcal{M}^{U,\sigma}]^* \tilde{v}|^2 d\mathcal{S} \: dt + 4 \sum_{k=0}^\infty \int_{\mathcal{S}} |\tilde{v}|^2 \tilde{v} \bcdot (\mathcal{R} G_\sigma^v(U) \phi_k - \nabla_H \tilde{\pi}_k) d\mathcal{S} \: d\beta^k \nonumber\\
    := &\sum_{j=1}^5 I_j dt + \sum_{k=0}^\infty I_6^k d\beta^k, \nonumber
\end{align}
remarking that $\int_\mathcal{S} |\tilde{v}|^2 \tilde{v} \bcdot (f\bold{k} \times \tilde{v}) = \int_\mathcal{S} |\tilde{v}|^2 \tilde{v} \bcdot (\bar{v} \bcdot \nabla_H\tilde{v}) =0$, since $\nabla_H \bcdot \bar{v} = 0$. For $I_1$, we reason similarly as in \cite{CT_2007} (Section 3.2), with $L^4$ rather than $L^6$, so that
\begin{align}
    |I_1| \leq \xi \int_\mathcal{S} |\nabla_3 \tilde{v}|^2 |\tilde{v}|^2 d\mathcal{S} + C_\xi \|\tilde{v}\|_{L^4}^4 (\|v\|_H^2\|v\|_V^2+1) + C\|T\|_H^2\|T\|_V^2.
\end{align}
\modif{
For $I_2$, by Sobolev embedding,
\begin{align}
    I_2 &\leq C \||\tilde{v}|^2 \nabla_3 \tilde{v}\|_{L^1} \|\nabla_3 \bcdot (a_3^K \nabla_3 (w(v)+w_s)) + K*[u_s \bcdot \nabla_3 (w(v)+w_s)]\|_{L^\infty} \nonumber\\
    &\leq C \|v\|_V \|\tilde{v}\|_{L^4}^2 \Big(\|\nabla_3 \bcdot (a_3^K \nabla_3 (w(v)+w_s))\|_{H^{2}} + \|K*[u_s \bcdot \nabla_3 (w(v)+w_s)]\|_{L^\infty}\Big) \nonumber\\
    &\leq C \|v\|_V \|\tilde{v}\|_{L^4}^2 (1+\|v\|_{H} \|K\|_{H^3}^2 + \|K\|_{H^2} \|u_s\|_{L^\infty} \|v\|_H + \|K\|_{H^1} \|u_s\|_{W^{1,\infty}} \|v\|_H) \nonumber\\
    &\leq C \|v\|_V \|\tilde{v}\|_{L^4}^2 (1+\|U\|_{H}) \leq C\|U\|_V^2(\|U\|_H^2 + 1) + C\|\tilde{v}\|_{L^4}^4.
\end{align}
}
For $I_3$,
\begin{align}
    I_3 &= - 4 \int_{\mathcal{S}} |\tilde{v}|^2 \tilde{v} \bcdot \Big[\mathcal{R} B(v,v_s) - \frac{1}{2} \nabla_3 \bcdot( a_{BT} \nabla_3 \tilde{v} ) \Big]  d\mathcal{S} \nonumber\\
    &- 4 \int_{\mathcal{S}} |\tilde{v}|^2 \tilde{v} \bcdot \mathcal{R}[ d_t v_s + \frac{1}{2} \nabla_3 \bcdot (a_{BT} \nabla_3 v_s) + Av_s + \Gamma v_s] \nonumber\\
    &\leq C_\xi (\|\tilde{v}\|_{L^4}^4+1) + 2 \int_{\mathcal{S}} |\tilde{v}|^2 \tilde{v} \bcdot (\nabla_3 \bcdot a_{BT} \nabla_3 \tilde{v}) \: d\mathcal{S} + \xi \int_{\mathcal{S}} |\tilde{v}|^2 |\nabla_3 \tilde{v}|^2 \: d\mathcal{S} \nonumber\\
    &\leq C_\xi (\|\tilde{v}\|_{L^4}^4+1) - 2 \int_\mathcal{S} ( a_{BT} \nabla_3 \tilde{v} : \nabla_3 \tilde{v}) |\tilde{v}|^2 - \underbrace{\frac{1}{2}( a_{BT} \nabla_3 |\tilde{v}|^2 \bcdot \nabla_3 |\tilde{v}|^2) d\mathcal{S}}_{\geq 0} +  \xi \int_{\mathcal{S}} |\tilde{v}|^2 |\nabla_3 \tilde{v}|^2 \: d\mathcal{S}.
\end{align}
For $I_4$, define $\psi_k = \mathcal{R}[\phi_k \bcdot \nabla_3 v + \phi_k \bcdot \nabla_3 v_s + A \phi_k^H + \Gamma \phi_k^H]$ and $\psi_k^w = \mathcal{R}[\phi_k \bcdot \nabla_3 w_s + A \phi_k^z]$, so that
\begin{align}
    |I_4| &\leq 2 \sum_{k=0}^\infty \int_\mathcal{S} |\tilde{v}|^2 |\psi_k + \psi_k^w - \phi_k^{BT} \bcdot \nabla_3 \tilde{v} - \phi_k^{BC} \bcdot \nabla_3 \bar{v} + \mathcal{A}(\phi_k^z \partial_z v) - \nabla_H \tilde{\pi}_k|^2 d\mathcal{S} \nonumber \\
    &\leq 2 \int_\mathcal{S} (a_{BT} \nabla_3 \tilde{v} : \nabla_3 \tilde{v}) \: |\tilde{v}|^2 d\mathcal{S} + 2 \int_\mathcal{S} |a_{BC} \nabla_3 \bar{v}| \: |\nabla_3 \bar{v}| \: |\tilde{v}|^2 d\mathcal{S} + \xi \int_\mathcal{S} |\nabla_3 \tilde{v}|^2 \: |\tilde{v}|^2 d\mathcal{S} \nonumber \\
    &+ \underbrace{C \sum_{k=0}^\infty \|\tilde{v}\|_{L^4}^2 \|\mathcal{A} \phi_k^z \partial_z v\|_{L^4}^2}_{\leq C \|\tilde{v}\|_{L^4}^2 \|\partial_z v\|_{L^2} \|\nabla_3 \partial_z v\|_{L^2}} + C_\xi (\|\tilde{v}\|_{L^4}^4 + \|U\|_V^2 + \|\partial_z \tilde{v}\|_{L^2}^2 +1) + C \sum_{k=0}^\infty \underbrace{\int_\mathcal{S} |\tilde{v}|^2 |\nabla_H \tilde{\pi}_k|^2 d\mathcal{S}}_{\leq C \|\tilde{v}\|_{L^4}^2 \|\nabla_H \tilde{\pi}_k\|_{L^4}^2} \nonumber\\
    &\leq 2 \int_\mathcal{S} (a_{BT} \nabla_3 \tilde{v} : \nabla_3 \tilde{v}) \: |\tilde{v}|^2 d\mathcal{S} + C_\xi (1+\|U\|_V^2) (\|\tilde{v}\|_{L^4}^4 + \|U\|_V^2 + \|\partial_z v\|_{L^2}^2 +1) \nonumber\\
    &+ \xi \int_\mathcal{S} |\nabla_3 \tilde{v}|^2 \: |\tilde{v}|^2 d\mathcal{S} +  \xi \|\bar{v}\|_{\mathcal{D}(A)}^2,
\end{align}
remarking that, by Hölder's inequality, the 2D Sobolev embedding, and interpolation inequality,
\begin{align}
    2 \int_\mathcal{S} |a_{BC} \nabla_3 \bar{v}| \: |\nabla_3 \bar{v}| \: |\tilde{v}|^2 d\mathcal{S} &\leq C \|\nabla_3 \bar{v}\|_{L^4}^2 \|\tilde{v}\|_{L^4}^2 \leq C \|\bar{v}\|_{\mathcal{D}(A^{3/4})}^2 \|\tilde{v}\|_{L^4}^2 \nonumber\\
    &\leq C \|\bar{v}\|_{V} \|\bar{v}\|_{\mathcal{D}(A)} \|\tilde{v}\|_{L^4}^2 \leq C \|U\|_{V}^2 \|\tilde{v}\|_{L^4}^4 + \xi \|\bar{v}\|_{\mathcal{D}(A)}.
\end{align}
For $I_5$, we have for any $\xi_1 >0$,
\begin{align}
    |I_5| &\leq 4 \sum_{k=0}^\infty \int_{\mathcal{S}} (\tilde{v} \bcdot (\psi_k + \psi_k^w- \phi_k \bcdot \nabla_3 \tilde{v} + \mathcal{A}(\phi_k^z \partial_z v) - \nabla_H \tilde{\pi}_k))^2 \label{baroclinic-I5}\\
    &\leq 4 \sum_{k=0}^\infty \int_{\mathcal{S}} \underbrace{ [\tilde{v} \bcdot (\phi_k \bcdot \nabla_3 \tilde{v} - \mathcal{A}(\phi_k^z \partial_z v))]^2}_{\leq \frac{1}{4} |\phi_k \bcdot \nabla_3 |\tilde{v}|^2|^2 + |\tilde{v}|^2 |\mathcal{A}(\phi_k^z \partial_z v)|^2 } d\mathcal{S} + C (\|\tilde{v}\|_{L^4}^4 + \|\nabla_H \tilde{\pi}_k\|_{L^4}^2 + 1) \nonumber \\
    &\leq \int_\mathcal{S} |a_3 \nabla_3 |\tilde{v}|^2| \: |\nabla_3 |\tilde{v}|^2| d\mathcal{S} + C_{\xi_1} (\|\partial_z v\|_{L^2}^2 +1)\|\tilde{v}\|_{L^4}^4 + \xi_1 \|\nabla_3\partial_z v\|_{L^2}^2 + C \sum_{k=0}^\infty \|\nabla_H \tilde{\pi}_k\|_{L^4}^2. \nonumber
\end{align}
Then, integrate and take the expectation, to reach
\begin{align}
    \E\Big[ |\tilde{v}\|_{L^4}^4\Big] + &(c-\alpha_3)\E\Big[ \int_\eta^\zeta \int_{\mathcal{S}} (|\nabla_H \tilde{v}|^2 |\tilde{v}|^2 + |\nabla_H |\tilde{v}|^2|^2 + |\partial_z \tilde{v}|^2 |\tilde{v}|^2 + |\partial_z |\tilde{v}|^2 |^2) d\mathcal{S} \: ds\Big] \\
    \leq C_{\xi_1} &\E\Big[1+ \|\tilde{v}(\eta)\|_{L^4}^4 + \int_\eta^\zeta (\|U\|_V^2 +1)\|\tilde{v}\|_{L^4}^4  \Big] + \xi_1 \E\Big[\int_\eta^\zeta \|\nabla_3 \partial_z v\|_{L^2}^2 \: ds \Big] \nonumber
\end{align}
Similarly by taking the supremum,
\begin{multline}
    \E\Big[\sup_{[\eta,\zeta]} |\tilde{v}\|_{L^4}^4\Big] + (c-\alpha_3)\E\Big[ \int_\eta^\zeta \int_{\mathcal{S}} (|\nabla_H \tilde{v}|^2 |\tilde{v}|^2 + |\nabla_H |\tilde{v}|^2|^2 + |\partial_z \tilde{v}|^2 |\tilde{v}|^2 + |\partial_z |\tilde{v}|^2 |^2) d\mathcal{S} \: ds\Big] \label{ito-barocline} \\
    \leq C_{\xi_1} \E\Big[1+ \|\tilde{v}(\eta)\|_{L^4}^4 + \int_\eta^\zeta (\|U\|_V^2 +1)\|\tilde{v}\|_{L^4}^4  \Big] + \xi_1 \E\Big[\int_\eta^\zeta \|\nabla_3 \partial_z v\|_{L^2}^2 \: ds \Big] \\
    + \E\Bigg[ \sup_{[\eta,\zeta]} \Big|\int_\eta^s \sum_{k=0}^\infty I_6^k d\beta_r^k\Big|\Bigg].
\end{multline}
Also, by the Burkholder-Davis-Gundy inequality (Theorem 4.36 in \cite{book_DPZ_2014}), one obtains
\begin{align}
    \E&\Bigg[\sup_{[\eta,\zeta]} \Big|\int_\eta^s \sum_{k=0}^\infty I_6^k d\beta_r^k \Big|\Bigg] \leq \E\Bigg[ \int_\eta^\zeta \Big|\int_{\mathcal{S}} |\tilde{v}|^2 \tilde{v} \bcdot \sum_{k=0}^\infty \Big(\mathcal{R} G_\sigma^v(U) \phi_k - \nabla_H \tilde{\pi}_k\Big) d\mathcal{S}\Big|^2 ds\Bigg] \nonumber\\
    &\leq C\E\Bigg[ \int_\eta^\zeta \underbrace{\Big|\int_{\mathcal{S}} |\tilde{v}|^2 \tilde{v} \bcdot \sum_{k=0}^\infty \Big(\mathcal{R} G_\sigma^v(U) \phi_k\Big) d\mathcal{S}\Big|^2}_{\leq \|\tilde{v}\|_{L^4}^2 \||\tilde{v}||\nabla_3 \tilde{v}|\|_{L^2} \leq C_\xi\|\tilde{v}\|_{L^4}^4 + \xi \||\tilde{v}||\nabla_3 \tilde{v}|\|_{L^2}^2} + \underbrace{\Big|\int_{\mathcal{S}} |\tilde{v}|^2 \tilde{v} \bcdot \sum_{k=0}^\infty \Big(\nabla_H \tilde{\pi}_k\Big) d\mathcal{S}\Big|^2}_{\leq \|\tilde{v}\|_{L^4}^3 \sum_{k=0}^\infty \|\nabla_H \tilde{\pi}_k\|_{L^4} } ds\Bigg] \nonumber\\
    &\leq C_\xi \E\Bigg[\int_\eta^\zeta(\|\tilde{v}\|_{L^4}^4+1)\Bigg] + \xi \E\Bigg[\int_\eta^\zeta \||\tilde{v}||\nabla_3 \tilde{v}|\|_{L^2}^2\Bigg].
\end{align}
Therefore, using \eqref{ito-barocline}, the following result holds upon choosing $\xi_1 = \frac{1}{16(C_2 \vee 1) }$ in equation \eqref{baroclinic-I5},
\begin{multline}
    \E\Big[\sup_{[\eta,\zeta]} |\tilde{v}\|_{L^4}^4\Big] + \E\Big[ \int_\eta^\zeta \int_{\mathcal{S}} (|\nabla_H \tilde{v}|^2 |\tilde{v}|^2 + |\nabla_H |\tilde{v}|^2|^2 + |\partial_z \tilde{v}|^2 |\tilde{v}|^2 + |\partial_z |\tilde{v}|^2 |^2) d\mathcal{S} \: ds\Big] \label{estimation-barocline}\\
    \leq C_3 \E\Big[1+ \|\tilde{v}(\eta)\|_{L^4}^4 + \int_\eta^\zeta (\|U\|_V^2 +1)\|\tilde{v}\|_{L^4}^4  \Big] + \frac{1}{16(C_2 \vee 1) } \E\Big[\int_\eta^\zeta \|\nabla_3 \partial_z v\|_{L^2}^2 \: ds \Big].
\end{multline}

\subsection{Intermediate estimates}

In this subsection, we gather the previous estimates on the vertical gradient and the barotropic and baroclinic modes of velocity, to obtain two intermediate estimates. This is used in the following to show the existence of a global pathwise solution. Define first, for $T>0$, and two stopping times $\eta < \zeta < T \wedge \tau$,
\begin{gather}
    X_s := \|\bar{v}\|_V^2 + \|\partial_z v\|_{L^2}^2 + \|\tilde{v}\|_{L^4}^4, \quad Y_s := \|\bar{v}\|_{\mathcal{D}(A)}^2 + \|\nabla_3 \partial_z v\|_{L^2}^2 + \int_{\mathcal{S}} |\tilde{v}|^2 |\nabla \tilde{v}|^2, \\
    \mathcal{E}(\eta, \zeta) := \E\Big[\sup_{s \in [\eta,\zeta]} X_s \Big] + \E\Big[\int_\eta^\zeta Y_s ds\Big].
\end{gather}
Then, gathering the previous estimates as follows: $\frac{\eqref{estimation-barotrope}}{8C_1 C_2} + \frac{\eqref{estimation-gradient-z}}{4 C_2} + \eqref{estimation-barocline}$, we reach
\begin{align}
    \mathcal{E}(\eta, \zeta) \leq C \E[1 + X_\eta] + C\E\Big[\int_\eta^\zeta (1+\|U\|_H^2)(1+\|U\|_V^2 )X_s \: ds\Big].
\end{align}
Define, for all $k>1$,
\begin{equation}
    l_k := \inf \Big\{t > 0\: \Big| \: \|U(t)\|_H + (\int_0^t \|U\|_V^2)^{1/2} \geq k \Big\} \wedge \tau,
\end{equation}
so that $\lim_{k \rightarrow \infty} \mathbbm{P}(l_k = \tau) =1$. Hence, using the stochastic Grönwall lemma of \cite{GHZ_2009} with $\eta < \zeta < l_k \wedge T$, for all $k$ there exists a $k$-dependent constant $C_k$, such that
\begin{align}
    \mathcal{E}(0, l_k \wedge T) \leq C_k (1+\|v_0\|_V^4).\label{intermediate-estimate}
\end{align}
Furthermore, by It\={o}'s lemma,
\begin{align}
    d_t\|\partial_z v\|_{L^2}^4 &+ 4 c \|\partial_z v\|_{L^2}^2 \|\nabla_3 \partial_z v\|_{L^2} ^2 dt \leq - 4 \|\partial_z v\|_{L^2}^2 ((\partial_z v \bcdot \nabla_H) v - (\nabla_H \bcdot v) \partial_z v,\partial_z v)_{L^2} dt \nonumber\\
    &+ 4 \|\partial_z v\|_{L^2}^2 ( g [\beta_T \nabla_H T + \beta_S \nabla_H S],\partial_z v)_{L^2} dt - 4 \|\partial_z v\|_{L^2}^2(\partial_z \hat{F}_\sigma^v(U), \partial_z v)_{L^2} dt \nonumber\\
    &+ 4 \|\partial_z v\|_{L^2}^2 \Big(\mathbf{P} \nabla_H \Big[\nabla_3 \bcdot (a_3^K \nabla_3 (w(v)+w_s))+ K*[(u_s \bcdot \nabla_3) (w(v)+w_s)]\Big], \partial_z v\Big)_{L^2} dt \nonumber\\
    & + 2 \|\partial_z v\|_{L^2}^2 \|\partial_z \hat{G}_\sigma^v(U)- \partial_z\mathcal{M}^{U,\sigma}\|_{\mathcal{L}_2(\mathcal{W},{L^2})}^2 dt + 4 \|[\partial_z \hat{G}_\sigma^v(U)- \partial_z\mathcal{M}^{U,\sigma}]^* \partial_z v\|_{\mathcal{L}_2(\mathcal{W},{L^2})}^2 dt \nonumber \\
    &+ 4 \|\partial_z v\|_{L^2}^2 \sum_{k=0}^\infty \int_\mathcal{S} \Big[\partial_z \hat{G}_\sigma^v (U) (\phi_k) - \partial_z{\nabla_H \pi_k}\Big]\bcdot \partial_z v \: d\beta_t^k.
\end{align}
Thus, reasoning as in subsection \ref{subsec-vertical-gradient},
\begin{multline}
    \E\Big[\sup_{[\eta,\zeta]} \|\partial_z v\|_{L^2}^4\Big] + \E\Big[ \int_\eta^\zeta \|\partial_z v\|_{L^2}^2 \|\nabla_3 \partial_z v\|_{L^2}^{2} \: ds\Big]
    \leq C \E\Big[ 1 + \|\partial_z v(\eta)\|_{L^2}^4 \Big] \\
    + C \E\Big[\int_\eta^\zeta (1+\|U\|_V^2  + \int_\mathcal{S} |\nabla \tilde{v}|^2 |\tilde{v}|^2 ) (1+\|\partial_z v\|_{L^2}^4) ds \Big] \\
    + C \E\Big[\int_\eta^\zeta\|[\partial_z \hat{G}_\sigma^v(U)- \partial_z\mathcal{M}^{U,\sigma}]^* \partial_z v\|_{\mathcal{L}_2(\mathcal{W},{L^2})}^2 ds\Big].
\end{multline}
It remains to estimate the last term on the RHS. Using the notations of subsection \ref{subsec-vertical-gradient}, notice that,
\begin{multline}
    \E\Big[\int_\eta^\zeta\|[\partial_z \hat{G}_\sigma^v(U)- \partial_z\mathcal{M}^{U,\sigma}]^* \partial_z v\|_{\mathcal{L}_2(\mathcal{W},{L^2})}^2 ds\Big] \\
    \leq C \E\Bigg[\int_\eta^\zeta \sum_{k=0}^\infty \Big|\Big(\partial_z \psi_k  - \partial_z(\phi_k \bcdot \nabla_3 v) - \mathbf{P} \nabla_H (\phi_k \bcdot \nabla_3 w(v) ), \partial_z v \Big)_{L^2}\Big|^2 ds \Bigg] \\
   \leq C \E\Bigg[(1+ \sup_{[\eta,\zeta]} \|\partial_z v\|_{L^2}^4 ) \int_\eta^\zeta (1+\|\nabla_3 \partial_z v\|_{L^2}^2) ds\Bigg],
\end{multline}
so that
\begin{multline}
    \E\Big[\sup_{[\eta,\zeta]} \|\partial_z v\|_{L^2}^4\Big] + \E\Big[ \int_\eta^\zeta \|\partial_z v\|_{L^2}^2 \|\nabla_3 \partial_z v\|_{L^2}^{2} \: ds\Big]
    \leq C \E\Big[ 1 + \|\partial_z v(\eta)\|_{L^2}^4 \Big] \\
    + C \E\Big[\int_\eta^\zeta (1+\|U\|_V^2  + \int_\mathcal{S} |\nabla \tilde{v}|^2 |\tilde{v}|^2 + \|\nabla_3 \partial_z v\|_{L^2}^2) (1+\|\partial_z v\|_{L^2}^4) ds \Big].
\end{multline}
Hence, the stochastic Grönwall lemma yields, with $\eta < \zeta < l_k \wedge T$ again,
\begin{equation}
    \E\Big[\sup_{[0,l_k \wedge T]} \|\partial_z v\|_{L^2}^4\Big] + \E\Big[ \int_0^{l_k \wedge T} \|\partial_z v\|_{L^2}^2 \|\nabla_3 \partial_z v\|_{L^2}^{2} \: ds\Big] \leq C_k \E\Big[ 1 + \|\partial_z v(0)\|_{L^2}^4 \Big].
\end{equation}
    
\CQFD

\subsection{Globality-in-time of the solution} \label{subsec-globality-final}

Our final argument for showing the existence of a global pathwise solution is inspired by what was proposed in \cite{AHHS_2022} (see the arguments p.38-39), yet different. In the following, we use a result that is similar to Proposition 4.1 of \cite{AHHS_2022} [Part 4], which we prove in the appendix (Appendix \ref{sec-appendix-1}): let $T>0$, and define two stopping times $\eta, \zeta < \tau \wedge T$, such that $\eta < \zeta$. Then the following holds,
\begin{equation}
    \E\Big[\sup_{s \in [\eta,\zeta]} \|v\|_{V}^2 \Big] + \E\Big[\int_\eta^\zeta \|v\|_{\mathcal{D}(A)}^2 ds\Big]
    \leq C \E\Big[1 + \|v(\eta)\|_{V}^2 
    + \int_\eta^\zeta \|v \bcdot \nabla v\|_{L^2}^2 ds + \int_\eta^\zeta \|w(v) \partial_z v\|_{L^2}^2 ds \Big]. \label{majoration-Agresti}
\end{equation}
%------------------------------------------------------
Moreover, remark that by Hölder inequality and Sobolev embeddings,
\begin{equation}
    \E\Big[ \int_\eta^\zeta \|v \bcdot \nabla v\|_{L^2}^2 \Big] \leq C_\xi \E\Big[ \int_\eta^\zeta (1+ \|v\|_V^2) \Big] +\xi \E \Big[\int_\eta^\zeta \|v\|_{\mathcal{D}(A)}^2 \Big].
\end{equation}
Furthermore,
\begin{align}
    \|w(v) \partial_z v\|_{L^2(\mathcal{S})}^2 &\leq \int_{-h}^0 \|w(v)\|_{L^4(\mathcal{S}_H)}^2 \|\partial_z v\|_{L^4(\mathcal{S}_H)}^2 dz \nonumber\\
    &\leq \Bigg(\sup_{[-h,0]} \|w(v)\|_{L^4(\mathcal{S}_H)}^2 \Bigg) \Bigg(\int_{-h}^0 \|\partial_z v\|_{L^4(\mathcal{S}_H)}^2 dz\Bigg) \nonumber\\
    &\leq C \Bigg( \int_{-h}^0 \|\nabla \bcdot v\|_{L^4(\mathcal{S}_H)}^2 dz\Bigg) \Bigg(\int_{-h}^0 \|\partial_z v\|_{L^4(\mathcal{S}_H)}^2 dz \Bigg) \nonumber\\
    &\leq C \Bigg( \int_{-h}^0 \|v\|_{H^1(\mathcal{S}_H)} \|v\|_{H^2(\mathcal{S}_H)} dz\Bigg) \Bigg(\int_{-h}^0 \|\partial_z v\|_{L^2(\mathcal{S}_H)} \|\partial_z v\|_{H^1(\mathcal{S}_H)} dz \Bigg) \nonumber\\
    &\leq C \|v\|_{H^1} \|v\|_{H^2} \|\partial_z v\|_{L^2} \|\partial_z v\|_{H^1} \leq C_\xi \|v\|_{V}^2  \|\partial_z v\|_{L^2}^2 \|\partial_z v\|_{H^1}^2 + \xi  \|v\|_{\mathcal{D}(A)}^2.
\end{align}
Thus,
\begin{align}
    \E\Big[ \int_\eta^\zeta \|w(v) \partial_z v\|_{L^2(\mathcal{S})}^2 \Big] \leq C_\xi \E \Big[ \int_\eta^\zeta \|v\|_{V}^2 \|\partial_z v\|_{L^2}^2 \|\partial_z v\|_{H^1}^2 \Big] + \xi \E \Big[ \int_\eta^\zeta \|v\|_{\mathcal{D}(A)}^2 \Big].
\end{align}
Using those estimates in equation \eqref{majoration-Agresti}, we reach
\begin{multline}
    \E\Big[  \sup_{[\eta,\zeta]} \|v\|_{V}^2 \Big] + \E\Big[ \int_\eta^\zeta \|v\|_{\mathcal{D}(A)}^2 ds \Big] \leq C \E\Big[1+ \|v(\eta)\|_{V}^2\Big] \\
    + C \E\Big[  \sup_{[\eta,\zeta]} \|v\|_{V}^2 \int_\eta^\zeta (1+  \|\partial_z v\|_{L^2}^2 \|\partial_z v\|_{H^1}^2) ds \Big].
\end{multline}
Therefore, by the stochastic Grönwall lemma of \cite{AV2024} (rather than the one of \cite{GHZ_2009}), the following holds a.s.,
\begin{equation}
    \sup_{[0,l_k \wedge T]} \|v\|_{V}^2 + \int_0^{l_k \wedge T} \|v\|_{\mathcal{D}(A)}^2 ds < \infty.
\end{equation}
Since $\lim_{k \rightarrow \infty} \mathbbm{P}(l_k = \tau) =1$, we deduce that, almost surely,
\begin{equation}
    \sup_{[0, \tau \wedge T]} \|v\|_{V}^2 + \int_0^{ \tau \wedge T} \|v\|_{\mathcal{D}(A)}^2 ds < \infty.
\end{equation}
Thus, $\tau = \infty$ a.s.
The aforementioned solution to the problem $(\hat{\mathcal{P}}_K)$ is hence global-in-time, with values almost surely in the space $$L_{loc}^2\Big([0,+\infty),\mathcal{D}(A)\Big) \cap C\Big([0,+\infty),V\Big).$$
\CQFD

\subsection{Continuity in probability with respect to the initial data and the noise operator} \label{subsec-continuity}

This subsection concerns the continuity of the solutions of $(\hat{\mathcal{P}}_K)$ with respect to initial data and to the noise. We prove the second point of Theorem \ref{theorem-filtered-strong} below.

\bigskip

Let $U_0^n$, $U_0$, $\sigma^n$, $\sigma$, $U^n$, $U$ and $R^n$ as in the second point of Theorem \ref{theorem-filtered-strong}, and let $T>0$. Remark that,
\begin{align}
    R^n \in L^2\big([0,T],\mathcal{D}(A)\big) \cap C\big([0,T],V\big), \quad \forall T >0 \text { a.s.}
\end{align}
Let $\eta < \zeta$ two positive stopping times. Reasoning as in the proof for uniqueness in subsection \ref{subsec-uniqueness}, we infer that
\begin{multline}
    \E\Big[ \sup_{[\eta,\zeta]} \|R^n\|_V^2 + \int_{\eta}^{\zeta} \|R^n\|_{\mathcal{D}(A)}^2 ds \Big] \leq \E\Big[ \|R^n(0)\|_V^2 \Big] + C \E\Big[\int_{\eta}^{\zeta} \|R\|_V^2 \Big( 1+ \|U\|_{\mathcal{D}(A)}^2 +  \|U^n\|_{\mathcal{D}(A)}^2 \Big) ds \Big]\\
    + C\E\Big[ \int_{\eta}^{\zeta} \delta_\sigma^n (1+\|U\|_{\mathcal{D}(A)}^2 + \|U^n\|_{\mathcal{D}(A)}^2 + (\beta_\sigma^n)^2) + (\delta_\sigma^n)^2 (\beta_\sigma^n + 1) ds \Big], \label{eq-continuity-before-gronwall}
\end{multline}
where
\begin{gather}
    \delta_\sigma^n := \sup_{[0,T]}\|\sigma^n - \sigma\|_{\mathcal{L}_2(\mathcal{W},H^4)}^2  + \|a_3 - a_3^n\|_{H^1([0,T],H^3)}, \\
    \beta_\sigma^n := \sup_{[0,T]}(\|\sigma \|_{\mathcal{L}_2(\mathcal{W},H^4)}^2 + \|\sigma^n\|_{\mathcal{L}_2(\mathcal{W},H^4)}^2) + \|a_3\|_{H^1([0,T],H^3)} +\| a_3^n\|_{H^1([0,T],H^3)}. \nonumber
\end{gather}
In addition, we define
\begin{equation}
    \mathcal{N}(\eta, \zeta ; X) := \sup_{[\eta,\zeta]} \|X\|_V^2 + \int_{\eta}^{\zeta} \|X\|_{\mathcal{D}(A)}^2 ds,
\end{equation}
and denote by
\begin{equation}
    \tau_{k,n}^* := \inf\{t > 0 \: | \: \mathcal{N}(0, t;U) + \mathcal{N}(0, t;U^n) > k\}.
\end{equation}
Using equation \eqref{eq-continuity-before-gronwall} with $(U^1,U^2,\tau') = (U, U^n, \tau_{k,n}^*)$, by the stochastic Grönwall lemma of \cite{AV2024}, we infer that there exists a constant $C_0$ such that, for all $\epsilon, K >0$, for all integer $k>0$,
\begin{multline}
    \mathbbm{P}(\mathcal{N}(0, \tau_{k,n}^* \wedge T; R^n)  \geq \epsilon) \leq \frac{C_0}{\epsilon} e^{C_0K} (\|R_0^n\|_{V}^2 + \delta_\sigma^n) \\
    + \mathbbm{P}\Big(\int_0^{\tau_{k,n}^* \wedge T} ( \|U\|_{\mathcal{D}(A)}^2 + \|U^n\|_{\mathcal{D}(A)}^2) ds + \chi_n \geq K \Big),
\end{multline}
where $\chi_n := T \delta_\sigma^n (1+\beta_\sigma^n) + T (1+ (\beta_\sigma^n)^2)$. Since $\delta_\sigma^n$ and $\beta_\sigma^n$ are bounded, so is $\chi_n$. Then, denote by $\hat{\chi} = \sup_{n} \chi_n$, and set $K = k+1+\hat{\chi}$, so that
\begin{equation}
    \mathbbm{P}(\mathcal{N}(0, \tau_{k,n}^* \wedge T; R^n)  \geq \epsilon) \leq \frac{C_0}{\epsilon} e^{C_0(k+1+\hat{\chi})} (\|R_0^n\|_{V}^2 + \delta_\sigma^n). \label{eq-tail-estimate}
\end{equation}
Hence, by choosing $k_n = \Big\lfloor \frac{- \ln (\|R_0^n\|_{V}^2 + \delta_\sigma^n)}{2C_0} \Big\rfloor \rightarrow \infty$ for a large enough $n$, we find that
\begin{equation}
    \mathbbm{P}(\mathcal{N}(0, \tau_{k_n,n}^* \wedge T; R^n)  \geq \epsilon)
    \leq \frac{C_0}{\epsilon} e^{C_0(1+\hat{\chi})} \sqrt{\|R_0^n\|_{V}^2 + \delta_n} \rightarrow 0.
\end{equation}
As the solutions $U$ and $U^n$ are global-in-time a.s., $\tau_{k_n,n}^* \rightarrow \infty$ in probability. Therefore,
\begin{equation}
    \mathbbm{P}(\mathcal{N}(0, T; R^n)  \geq \epsilon) \rightarrow 0,
\end{equation}
that is to say $U^n \rightarrow U$ in probability, in the space $C([0,T],V) \cap L^2([0,T], \mathcal{D}(A))$.

\CQFD

%\renewtheorem{definition}{Definition}[section]

\appendix

\section{Appendix 1: Proof of equation \eqref{majoration-Agresti} -- An estimate on the velocity and the temperature} \label{sec-appendix-1}

We show equation \eqref{majoration-Agresti}, which is similar to Proposition 4.1 of \cite{AHHS_2022}, and which we use in the final globality argument. The result is the following lemma,
\begin{lemma}
    Let $T>0$ a real number. Let $U = (v,T,S)^\tr$ the solution of either $(\mathcal{P}_K)$ or $(\hat{\mathcal{P}}_K)$ -- equations \eqref{problem-PK} and \eqref{problem-approx-PK} -- given by the second point of Theorem \ref{theorem-filtered-weak}, and $\tau_0$ its associated a stopping time. In addition, let two stopping times $\eta, \zeta < \tau_0 \wedge T$, such that $\eta < \zeta$. Then,
    \begin{align}
        \E\Big[\sup_{s \in [\eta,\zeta]} \|v\|_{V}^2 \Big] + \E\Big[\int_\eta^\zeta \|v\|_{\mathcal{D}(A)}^2 ds\Big] \leq &C \E\Big[ 1+ \|v(\eta)\|_{V}^2 + \int_\eta^\zeta \|v \bcdot \nabla v\|_{L^2}^2 ds + \int_\eta^\zeta \|w(v) \partial_z v\|_{L^2}^2 ds \Big].
    \end{align}
\end{lemma}
\emph{Proof:} We will prove this lemma for the problem $(\mathcal{P}_K)$, the proof for $(\hat{\mathcal{P}}_K)$ being similar. Remind first the definition of $(\mathcal{P}_K)$,
\begin{equation}
    \partial_t U^* + \Big[A U^* + B (U^*) + \Gamma U^* + \frac{1}{\rho_0} \mathbf{P} \nabla_H p + F_\sigma (U^*)\Big] dt =  G_\sigma(U^*) dW_t - \frac{1}{\rho_0} \mathbf{P} \nabla_H dp_t^\sigma.
\end{equation}
Then, let $\theta = (T,S)^\tr$, so that
\begin{align}
    \partial_t v^* + \Big[A v^* + B (v^*) + \Gamma v^* + \frac{1}{\rho_0} \mathbf{P} \nabla_H p + F_\sigma^v (v^*)\Big] dt &=  G_\sigma^v(v^*) dW_t - \frac{1}{\rho_0} \mathbf{P} \nabla_H dp_t^\sigma, \\
    \partial_t \theta + \Big[A \theta + B (v^*,\theta) +  F_\sigma^\theta(\theta) \Big] dt &=  G_\sigma^\theta(\theta) dW_t,
\end{align}
where $F_\sigma =: (F_\sigma^v, F_\sigma^\theta)^\tr$ and $G_\sigma =: (G_\sigma^v, G_\sigma^\theta)^\tr$. Let $t >0$, and apply It\={o}'s lemma with $\|.\|_V^2$ and $\|.\|_H^2$,
\begin{align}
    d_t \|v_n^*\|_V^2 = &- 2 (Av_n^*, A v_n^* + B(v_n^*) + \Gamma v_n^* + \frac{1}{\rho_0} \mathbf{P} \nabla_H p + F_\sigma^v (v_n^*))_H dt \label{ito-v_n^*}  \\
    &+ 2 (Av_n^*, G_\sigma^v(v_n^*)dW_t -  \frac{1}{\rho_0} \mathbf{P} \nabla_H dp_t^\sigma)_H + \|G_\sigma^v(v_n^*) - \mathcal{M}^{U,\sigma}\|_{\mathcal{L}_2(\mathcal{W},V)}^2 dt,\nonumber \\
    d_t \|\theta_n\|_H^2 = &- 2 (\theta_n, A \theta_n + B(v_n^*,\theta_n) + F_\sigma^\theta (\theta_n))_H dt  + 2 (\theta_n, G_\sigma^\theta(\theta_n)dW_t)_H + \|G_\sigma^\theta(\theta_n) \|_{\mathcal{L}_2(\mathcal{W},H)}^2 dt. \label{ito-theta_n}
\end{align}
Also, remind that the regularising properties of $K$ lead to
\begin{equation}
     \frac{1}{2} \sum_{k=0}^\infty \|\mathbf{P} \nabla_H \pi_k\|_H^2 \leq C(1+ \|U\|_H^2), \quad  \frac{1}{2} \sum_{k=0}^\infty \|\mathbf{P} \nabla_H \pi_k\|_V^2 \leq C(1+ \|U\|_V^2).
\end{equation}

\paragraph{Estimate of the bounded variation terms}

For the advection term we have
\begin{align}
    (Av_n^*, B^n(v_n^*))_H &\leq  \|A v_n^*\|_H \|B^n(v_n^*)\|_H \leq \xi \|AU_n^*\|_H^2 + C_\xi \Big( \|v_n^* \bcdot \nabla v_n^*\|_{L^2}^2 + \|w_n^* \partial_z v_n^*\|_{L^2}^2 \Big).
\end{align}
Moreover, for all $\xi > 0$,
\begin{align}
    |(Av_n^*, \Gamma^n (v_n^*))_H| + |(A v_n^*, (\frac{1}{\rho_0} \mathbf{P} \nabla_H p)^n)_H| &\leq C\|v_n^*\|_{\mathcal{D}(A)} \|U_n^*\|_V \leq C_\xi (1+\|U_n^*\|_V^2) + \xi \|v_n^*\|_{\mathcal{D}(A)}^2.
\end{align}
For $F_\sigma^v$, by Young's inequality,
\begin{align}
    (A v_n^*, F_\sigma^v (v_n^*))_H \leq M_\sigma^1 + M_\sigma^2 + C\|Av_n^*\|_H \leq M_\sigma^1 + M_\sigma^2 + C_\xi + \xi\|Av_n^*\|_H,
\end{align}
with
\begin{align}
    |M_\sigma^1| &\leq |b(v_n^*,v_s, Av_n^*)| \leq C \|v_n^*\|_V^{1/2} \|v_n^*\|_{\mathcal{D}(A)}^{3/2} \leq C_\xi \|v_n^*\|_V^2 + \xi \|v_n^*\|_{\mathcal{D}(A)}^2,\\
    M_\sigma^2 &= -\frac{1}{2} (D_{\mu,\nu}\nabla_3 \bcdot a_3 \nabla_3 v_n^*, Av_n^*)_{L^2} = -\frac{1}{2}(D_{\mu,\nu}\partial_i (a_{ij} \partial_j v_{n,k}^*), - \mu_l \partial_{ll} v_{n,k}^* )_{L^2} \nonumber\\
    &= -\frac{1}{2}(D_{\mu,\nu}\partial_l (a_{ij} \partial_j v_{n,k}^*), - \mu_l \partial_{il} v_{n,k}^* )_{L^2} \nonumber \\
    &= -\frac{1}{2}(D_{\mu,\nu}(a_{3,ij} \partial_{jl} v_{n,k}^*), - \mu_l \partial_{il} v_{n,k}^* )_{L^2} + C_\xi \|U_n^*\|_V^2 + \xi \|U_n^*\|_{\mathcal{D}(A)}^2,
\end{align}
where $D_{\mu,\nu} = Diag(\mu,\mu,\nu)$, and using Einstein's notation convention of summation over repeated indices, with $i,j,l \in \{x,y,z\}$ and $k \in \{x,y\}$. We also used the notation $\mu_x= \mu_y = \mu$ and $\mu_z = \nu$.

\paragraph{Estimate of the martingale terms}

Let $\psi_k = \mathbf{P}[\phi_k \bcdot \nabla_3 v + \phi_k \bcdot \nabla_3 v_s] + A \phi_k^H + C \phi_k^H$ and $\psi_k^w = \mathbf{P} K*[\phi_k \bcdot \nabla_3 w_s] + A \phi_k^z$, so that
\begin{multline}
    \frac{1}{2} \|G_\sigma^n(v_n^*) - \mathcal{M}^{U,\sigma}\|_{\mathcal{L}_2(\mathcal{W},V)}^2  = \frac{1}{2} \sum_{k=0}^\infty \|\psi_k + \psi_k^w - \mathbf{P}[\phi_k\bcdot \nabla_3 v_n^*] - \mathbf{P} \nabla_H \pi_k\|_V^2 \\
    \leq \frac{1}{2} \sum_{k=0}^\infty \|\psi_k + \psi_k^w\|_V^2 + \|\mathbf{P}[\phi_k \bcdot \nabla_3 v_n^*]\|_V^2 + 2 C \underbrace{|(\nabla_3 (\phi_k \bcdot \nabla_3 v_n^*), \nabla_3 (\psi_k + \psi_k^w))_{L^2}|}_{\leq C |b(\phi_k, v_n^*, \Delta_3 (\psi_k + \psi_k^w))|} \\
    + \|\mathbf{P} \nabla_H \pi_k\|_V^2 + C \underbrace{|(\nabla_3 \mathbf{P} \nabla_H \pi_k,\nabla_3 (\psi_k + \psi_k^w))_H|}_{\leq C \|\mathbf{P} \nabla_H \pi_k\|_H \|\psi_k + \psi_k^w\|_{H^2} } + C |(\mathbf{P}[\phi_k \bcdot \nabla_3 v_n^*], \mathbf{P} \nabla_H \pi_k)_V|.
\end{multline}
Additionally, remark that
\begin{align*}
    \sum_{k=0}^\infty (\mathbf{P}[\phi_k \bcdot \nabla_3 v_n^*], &\mathbf{P} \nabla_H \pi_k)_V \leq C_\xi \sum_{k=0}^\infty \|\mathbf{P} \nabla_H \pi_k\|_V^2 + \xi \|v_n^*\|_{\mathcal{D}(A)}^2.
\end{align*}
Moreover, for all $\xi >0$,
\begin{align*}
    \sum_{k=0}^\infty \|\mathbf{P}[\phi_k &\bcdot \nabla_3 v_n^*]\|_V^2 \leq \sum_{k=0}^\infty \|D_\mu^{1/2} \nabla_3 (\phi_k \bcdot \nabla_3 v_n^*)\|_{L^2(\mathcal{S},\R^3)}^2\\
    &\leq C_\xi \sum_{k=0}^\infty \|(\nabla_3 \phi_k)^\tr . \nabla_3 v_n^* \|_{L^2(\mathcal{S},\R^{3\times2})}^2 + (1+\xi) \sum_{k=0}^\infty \|D_\mu^{1/2} \phi_k \bcdot \nabla_3 (\nabla_3 v_n^*)\|_{L^2(\mathcal{S},\R^{3\times2})}^2\\
    &\leq C_\xi \|v_n^*\|_V^2 + (1+\xi) (D_{\mu,\nu} a_{3,ij} \partial_{il} v_{n,k}^*, \mu_l \partial_{il} v_{n,k}^*)_H,
\end{align*}
where we used Einstein's notation with $i,j,l \in \{x,y,z\}$ and $k \in \{x,y\}$ again. Hence, for all $\xi >0$,
\begin{align}
    \frac{1}{2} \|G_\sigma^n(v_n^*) - \mathcal{M}^{U_n,\sigma}\|_{\mathcal{L}_2(\mathcal{W},V)}^2 &\leq - M_\sigma^2 + C_\xi (1+\|v_n^*\|_V^2) + \xi \|U_n^*\|_{\mathcal{D}(A)}^2.
\end{align}
Thus, using the equation \eqref{ito-v_n^*}, we reach
\begin{multline}
    d_t \|v_n^*\|_V^2 + \|v_n^*\|_{\mathcal{D}(A)}^2 dt \leq C \Big(\|v_n^* \bcdot \nabla v_n^*\|_H^2 + \|w_n^* \partial_z v_n^*\|_H^2 \Big)\\
    + 2 (v_n^*, G_\sigma^v(v_n^*)dW_t)_H  + C (1+\|U_n^*\|_V^2) dt.
\end{multline}
For all $T>0$ and all stopping times $0 \leq \eta < \zeta \leq T$, by integrating, then taking the supremum and the expectation, we have
\begin{align}
    \E[\sup_{[\eta,\zeta]}\|v_n^*\|_V^2] + \E[ \int_{\eta}^{\zeta} \|v_n^*\|_{\mathcal{D}(A)}^2 dr] &\leq \E[\|v_n^*(\eta)\|_V^2] +  C \E\Big[ \int_{\eta}^{\zeta} \|v_n^* \bcdot \nabla v_n^*\|_{L^2}^2 + \|w_n^* \partial_z v_n^*\|_{L^2}^2 dr\Big] \label{énergie-v_n^*} \\
    &+ C \E[\int_{\eta}^{\zeta} (1+ \|U_n^*\|_V^2) dr] + \E\Bigg[ \sup_{[\eta,\zeta)} \Big|\int_{\eta}^s \sum_{k=0}^\infty I_7^k d\beta_s^k\Big|\Bigg]. \nonumber
\end{align}
Furthermore,
\begin{align}
    \E\Bigg[ \sup_{[\eta,\zeta]} \Big|\int_{\eta}^{\bcdot} \sum_{k=0}^\infty I_7^k d\beta_r^k\Big|\Bigg] &\leq C \E\Bigg[\Bigg( \int_{\eta}^{\zeta} \sum_{k=0}^\infty \underbrace{ \Big|\Big(\psi_k + \psi_k^w - \mathbf{P}[\phi_k \bcdot \nabla_3 v_n^*] - \mathbf{P} \nabla_H \pi_k, Av_n^* \Big)_H\Big|^2 }_{\leq C(\|\psi_k\|_H^2 + \|\phi_k\|_{L^\infty}^2 (1+\|v_n^*\|_V^2) + \|\mathbf{P} \nabla_H \pi_k\|_H^2 )\|Av_n^*\|_H^2} ds \Bigg)^{1/2} \Bigg] \nonumber \\
    &\leq C \E\Bigg[\Bigg(\sup_{[\eta,\zeta]} (1+\| v_n^*\|_V^2) \Bigg)^{1/2} \Bigg( \int_{\eta}^{\zeta} \|Av_n^*\|_H^2 ds \Bigg)^{1/2}\Bigg] \nonumber\\
    &\leq C_\xi \E\Bigg[\sup_{[\eta,\zeta]} (1+\|v_n^*\|_V^2) \Bigg] + \xi \E\Bigg[ \int_{\eta}^{\zeta} \|U_n^*\|_{\mathcal{D}(A)}^2 ds \Bigg]. 
\end{align}
Hence, using Lemma \ref{lemma-basic-estimate},
\begin{multline}
    \E[\sup_{[\eta,\zeta]}\|v_n^*\|_V^2] + \E[ \int_{\eta}^{\zeta} \|v_n^*\|_{\mathcal{D}(A)}^2 dr] \leq C \Big(1+ \E[\|v_n^*(\eta)\|_V^2]\Big) \\
    + C \E\Big[ \int_{\eta}^{\zeta} \|v_n^* \bcdot \nabla v_n^*\|_{L^2}^2 + \|w_n^* \partial_z v_n^*\|_{L^2}^2 dr\Big].
\end{multline}
That is the result we seek.
\CQFD

\section{Appendix 2: Energy estimates for Theorems \ref{theorem-EB} \& \ref{theorem-EV}}

Now we give a sketch of the proof for the energy estimates for the problems $(\mathcal{P}^{\gamma_r}_{EB})$ and $(\mathcal{P}^{\gamma_r}_{EV})$ -- equations \eqref{problem-PEB} and \eqref{problem-PEV}. These are used in the proof of Theorem 2 and 3, which run similarly as for Theorem 1. Therefore, the other steps, including the straightforward existence of Galerkin solutions, are omitted.

\subsection{Energy estimate for the problem $(\mathcal{P}^{\gamma_r}_{EB})$ (Theorem \ref{theorem-EB})}

Assume that $\gamma_r > 1$. For the energy estimates, we only need to check for the pressure term, the noise terms and the additional correction terms $\frac{1}{2} \mathbf{P} \nabla \bcdot \check{a}[\nabla w]$ and $\mathbf{P} \mathcal{C}_\sigma$.

\paragraph{Covariation correction term}
\begin{align*}
    \frac{1}{2} (\mathbf{P} \nabla \bcdot \check{a}[\nabla w(v_n^*)],v_n^*)_H = - \frac{1}{2} (\check{a}[\nabla w(v_n^*)],\nabla v_n^*)_{L^2}
\end{align*}
and $\|\mathbf{P} \mathcal{C}_\sigma\|_H \leq \sum_{k=0}^\infty \|\phi_k\|_{L^\infty} \|\phi_k\|_{H^4} \leq (\sum_{k=0}^\infty \|\phi_k\|_{L^\infty}^2)^{1/2} (\sum_{k=0}^\infty\|\phi_k\|_{H^4}^2)^{1/2} \leq C.$
\paragraph{Bounded variation pressure term}
\begin{align*}
    \Big|\Big( &U_n^*,\frac{1}{\rho_0} \mathbf{P} \nabla_H p \Big)_H\Big| \leq \Big|\Big(g \:  \nabla_H \int_z^0 {(\beta_T T_n + \beta_S S_n) dz'} , v_n^* \Big)_{L^2}\Big| \\
    &+ \Big|\Big(\nabla_H \int_z^0 \Big[v_s \bcdot \nabla_H w(v_n^*) + w_s \partial_z w(v_n^*)\Big] dz', v_n^* \Big)_{L^2}\Big|\\
    &+ \Big(\nabla_H \int_z^0 \Big[ \frac{1}{2} \nabla \bcdot \hat{a}[\nabla v_n^*] + \frac{1}{2} \nabla \bcdot \hat{\hat{a}} [\nabla w(v_n^*)]  - \alpha (-\Delta)^{\gamma_r} w(v_n^*)\Big] dz' ,v_n^*  \Big)_{L^2}\\
    &+ \Big|\Big( \nabla_H \int_z^0 \Big[v_s \bcdot \nabla_H w_s + w_s \partial_z w_s - \Big[\frac{1}{2} \nabla \bcdot \hat{a}[\nabla v_s] + \frac{1}{2} \nabla \bcdot \hat{\hat{a}} [\nabla w_s] \Big] \Big] dz',v_n^* \Big)_{L^2}\Big|\\
    &+\Big|\Big( \nabla_H \int_z^0 \Big[\frac{1}{2}\sum_k \phi_k \bcdot \nabla \partial_z \pi_k + \frac{1}{2} \sum_k (\phi_k \bcdot \nabla) \Delta \phi_k^z dz'\Big],v_n^* \Big)_{L^2}\Big|\\
    &\leq \cancel{\Big|\Big( v_s \bcdot \nabla_H w(v_n^*) + w_s \partial_z w(v_n^*), w(v_n^*)  \Big)_{L^2}\Big|} + \Big( \frac{1}{2} \Big[\nabla \bcdot \hat{a}[\nabla v_n^*] + \nabla \bcdot \hat{\hat{a}} [\nabla w(v_n^*)] \Big], w(v_n^*)  \Big)_{L^2}\\
    &- (\alpha (-\Delta)^{\gamma_r} w(v_n^*) , w(v_n^*)  \Big)_{L^2} + C \|U_n^*\|_H \|U_n^*\|_V + C \|U_n^*\|_V + \frac{1}{2} \Big|\Big(\sum_k \phi_k \bcdot \nabla \partial_z \pi_k, w(v_n^*) \Big)_{L^2}\Big|\\
    &\leq -\frac{1}{2} \Big( \hat{a}[\nabla v_n^*] + \hat{\hat{a}} [\nabla w(v_n^*)], \nabla_3 w(v_n^*)  \Big)_{L^2} - \alpha \|w(v_n^*)\|_{H^{\gamma_r}}^2+ C \|U_n^*\|_H \|U_n^*\|_V + C \|U_n^*\|_V\\
    &+C \Big(\sum_k \|\phi_k\|_{L^\infty}\Big) \underbrace{\|\nabla w(v_n^*)\|_{L^2} \|\partial_z \pi_k\|_{L^2}}_{\leq C(1+\|w(v_n^*)\|_{H^1}^2)},
\end{align*}
where we used that $\phi_k$ is divergence-free, and that $\partial_z \pi_k = -\phi_k \bcdot \nabla w + \Delta \phi_k^z$. Then, we may conclude by Sobolev interpolation and Young's inequality.
\paragraph{First term of the total covariation noise contribution}
\begin{align*}
    \Big\|\Big[G_\sigma^n(U_n^*) - &\mathcal{M}^{\sigma,U_n^*}\Big]^*U_n^*\Big\|_H^2 = \sum_{k=0}^\infty \Big[\cancel{(w(v_n^*), \phi_k \bcdot \nabla_3 w(v_n^*))_{L^2}} + (w(v_n^*), \phi_k \bcdot \nabla_3 w_s)_{L^2}\Big] \\
    &+ \sum_{k=0}^\infty \Big[(w(v_n^*), A \phi_k^z )_{L^2} - \cancel{(v_n^*,[\phi_k \bcdot \nabla_3 v_n^*])_{L^2}} - (v_n^*, \phi_k \bcdot \nabla_3 v_s + A \phi_k^H + \Gamma \phi_k^H)_{L^2} \Big]^2 \\
    &\leq \sum_{k=0}^\infty C \|\phi_k\|_{H^3}^2 (\|v_n^*\|_H^2 + \|w(v_n^*)\|_{L^2}^2) \leq C (\|U_n^*\|_H^2 +\|w(v_n^*)\|_{L^2}^2).
\end{align*}
\paragraph{Second term of the total covariation noise contribution}
\begin{multline*}
    \frac{1}{2} \Big\|G_\sigma^n(U_n^*) - \mathcal{M}^{\sigma,U_n^*}\Big\|_{\mathcal{L}_2(\mathcal{W},H)}^2  = \frac{1}{2} \sum_{k=0}^\infty \Big\|\mathbf{P} [\phi_k \bcdot \nabla_3 v_n^* + \phi_k \bcdot \nabla_3 v_s] + A \phi_k^H + \Gamma \phi_k^H\Big\|_H^2 \\
    + \frac{1}{2} \sum_{k=0}^\infty \Big\|\mathbf{P} \nabla_H \pi_k \Big\|_H^2 + \sum_{k=0}^\infty (\mathbf{P}[\phi_k \bcdot \nabla_3 (v_n^*+v_s)] + A \phi_k^H + \Gamma \phi_k^H, \mathbf{P} \nabla_H \pi_k)_{L^2}\\
    \leq \frac{1}{2} (a_3\nabla_3 v_n^*,\nabla_3 v_n^*)_{L^2} + C_\xi (1+\|U_n^*\|_H^2) + \xi (\|U_n^*\|_V^2 + \|w(v_n^*)\|_{H^{\gamma_r}}^2) \\
    + \Big( \hat{a}[\nabla_3 v_n^*] + \frac{1}{2} \hat{\hat{a}} [\nabla_3 w(v_n^*)], \nabla_3 w(v_n^*) \Big)_{L^2}.
\end{multline*}
Using those inequalities, we conclude by similar arguments as in subsection \ref{subsec-energy-1} (Step 3). Notice that the ``new'' covariation compensation terms lead to the following $w$-dependent energy term $( \hat{a}[\nabla_3 v_n^*] + \frac{1}{2} \hat{\hat{a}} [\nabla_3 w(v_n^*)], \nabla_3 w(v_n^*) )_{L^2},$ arising from the computation of the bounded variation pressure associated energy $( U_n^*,\frac{1}{\rho_0} \mathbf{P} \nabla_H p )_H$. As expected, it balances exactly the quadratic covaration term arising from the martingale component $\frac{1}{2} \|G_\sigma^n(U_n^*)- \mathcal{M}^{\sigma,U_n^*}\|_{\mathcal{L}_2(\mathcal{W},H)}^2$. Furthermore, the term $\xi \|w(v_n^*)\|_{H^{\gamma_r}}^2$ is absorbed by the additional viscosity term $\alpha \|w(v_n^*)\|_{H^{\gamma_r}}^2 \geq \xi \|w(v_n^*)\|_{H^1}^2$. In particular, for all $T>0$, there exists a constant $C$ such that, for all $n \in \mathbbm{N}$,
$$\int_0^T \|w(v_n^*)\|_{H^{\gamma_r}}^2 \leq C.$$
\CQFD

\subsection{Energy estimate for the problem $(\mathcal{P}^{\gamma_r}_{EV})$ (Theorem \ref{theorem-EV})}

Assume that $\gamma_r > 2$. For the energy estimates, we only need to check for the bounded variation pressure and noise terms. Furthermore, notice that
\begin{align*}
    \Big\|\mathbf{P} \nabla_H \pi_k \Big\|_H^2 \leq C (1+ \|w(v_n^*)\|_{H^2}^2) \leq C_\xi (1+ \|U_n^*\|_H^2) + \xi \|w(v_n^*)\|_{H^{\gamma_r}}^2.
\end{align*}

\paragraph{Bounded variation pressure term}
\begin{align*}
    \Big|\Big( U_n^*,\frac{1}{\rho_0} \mathbf{P} \nabla_H p \Big)_H\Big| &\leq \Big|\Big(g \:  \nabla_H \int_z^0 {(\beta_T T_n + \beta_S S_n) dz'} , v_n^* \Big)_{L^2}\Big| \\
    &+ \Big|\Big(\nabla_H \int_z^0 \Big[(v_s \bcdot \nabla_H + w_s \partial_z) w(v_n^*)\Big] dz' ,v_n^*  \Big)_{L^2}\Big|\\
    &+ \Big(\nabla_H \int_z^0 \Big[ \frac{1}{2} \nabla \bcdot (a\nabla v_n^*)  - \alpha (-\Delta)^{\gamma_r} w(v_n^*)\Big] dz' ,v_n^*  \Big)_{L^2} \\
    &+ \Big|\Big( \nabla_H \int_z^0 \Big[(v_s \bcdot \nabla_H + w_s \partial_z) w_s - \frac{1}{2} \nabla \bcdot (a\nabla v_s) \Big] dz',v_n^* \Big)_{L^2}\Big|\\
    &\leq \cancel{\Big|\Big( v_s \bcdot \nabla_H w(v_n^*) + w_s \partial_z w(v_n^*), w(v_n^*)  \Big)_{L^2}\Big|} + \Big( \frac{1}{2} \nabla \bcdot (a\nabla v_n^*) , w(v_n^*) \Big)_{L^2} \\
    &- \Big(\alpha (-\Delta)^{\gamma_r} w(v_n^*) , w(v_n^*) \Big)_{L^2} +  C (\|U_n^*\|_H+1) \|U_n^*\|_V\\
    &\leq \underbrace{\frac{1}{2} \Big|\Big( a\nabla v_n^*, \nabla w(v_n^*) \Big)_{L^2}\Big|}_{\leq C \|U_n^*\|_H \|w(v_n^*)\|_{H^2}} - \alpha \|w(v_n^*)\|_{H^{\gamma_r}}^2+ C_\xi (\|U_n^*\|_H^2+1) + \xi \|U_n^*\|_V^2.
\end{align*}
\paragraph{First term of the total covariation noise contribution}
\begin{align*}
    \Big\|\Big[G_\sigma^n(U_n^*) - &\mathcal{M}^{\sigma,U_n^*}\Big]^*U_n^*\Big\|_H^2 = \sum_{k=0}^\infty \Big[\cancel{(w(v_n^*), \phi_k \bcdot \nabla_3 w(v_n^*))_{L^2}} + (w(v_n^*), \phi_k \bcdot \nabla_3 w_s)_{L^2} \Big]\\
    &+ \sum_{k=0}^\infty \Big[(w(v_n^*), A \phi_k^z )_{L^2} - \cancel{(v_n^*,[\phi_k \bcdot \nabla_3 v_n^*])_{L^2}} - (v_n^*, \phi_k \bcdot \nabla_3 v_s + A \phi_k^H + \Gamma \phi_k^H)_{L^2} \Big]^2 \\
    &\leq \sum_{k=0}^\infty C \|\phi_k\|_{H^3}^2 \|v_n^*\|_H^2 \leq C \|U_n^*\|_H^2.
\end{align*}
\paragraph{Second term of the total covariation noise contribution}
\begin{align*}
    \frac{1}{2} \Big\|G_\sigma^n(U_n^*)- \mathcal{M}^{\sigma,U_n^*}\Big\|_{\mathcal{L}_2(\mathcal{W},H)}^2  &= \frac{1}{2} \sum_{k=0}^\infty \Big\|\mathbf{P}[\phi_k \bcdot \nabla_3 v_n^* + \phi_k \bcdot \nabla_3 v_s] + A \phi_k^H + \Gamma \phi_k^H\Big\|_H^2 + \Big\|\mathbf{P} \nabla_H \pi_k \Big\|_H^2\\
    &+ \sum_{k=0}^\infty (\phi_k \bcdot \nabla_3 (v_n^*+v_s) + A \phi_k^H + \Gamma \phi_k^H, \mathbf{P} \nabla_H \pi_k)_{L^2}\\
    \leq \frac{1}{2} (a_3\nabla_3 v_n^*,\nabla_3 v_n^*)_{L^2} + &C_\xi(1+\|U_n^*\|_H^2) + \xi (\|U_n^*\|_V^2 + \|w(v_n^*)\|_{H^{\gamma_r}}^2) + \underbrace{\frac{1}{2} \Big( a\nabla v_n^*, \nabla w(v_n^*) \Big)_{L^2}}_{\leq C \|U_n^*\|_H \|w(v_n^*)\|_{H^2}}.
\end{align*}
Again, we conclude by using similar arguments as in subsection \ref{subsec-energy-1} (Step 3).
\CQFD

\paragraph*{Acknowledgements:}
The authors acknowledge the support of the ERC EU project 856408-STUOD and  benefit from the support of the French government ``Investissements d'Avenir'' program ANR-11-LABX-0020-01. AD acknowledges the support of the ANR project ADA.

\noindent
\modif{Additionally, the authors thank the reviewer for their careful reading and insightful comments.}

\bibliography{sn-bibliography}% common bib file
%% if required, the content of .bbl file can be included here once bbl is generated
%%\input sn-article.bbl

\end{document}